# Curl Constraint-Preserving Reconstruction and the Guidance it gives for Mimetic Scheme Design

By

Dinshaw S. Balsara (dbalsara@nd.edu),

Roger Käppeli (roger.kaepelli@sam.math.ethz.ch)

Walter Boscheri (walter.boscheri@unife.it) and

Michael Dumbser (michael.dumbser@unitn.it)

**Abstract**

Several important PDE systems, like magnetohydrodynamics and computational electrodynamics, are known to support involutions where the divergence of a vector field evolves in divergence-free or divergence constraint-preserving fashion. Recently, new classes of PDE systems have emerged for hyperelasticity, compressible multiphase flows, so-called first order reductions of the Einstein field equations, or a novel first order hyperbolic reformulation of Schrödinger's equation, to name a few, where the involution in the PDE supports curl-free or curl constraint-preserving evolution of a vector field. We study the problem of curl constraint-preserving reconstruction as it pertains to the design of mimetic finite volume (FV) WENO-like schemes for PDEs that support a curl-preserving involution. (Some insights into Discontinuous Galerkin (DG) schemes are also drawn, though that is not the prime focus of this paper.) This is done for two and three dimensional structured mesh problems where we deliver closed form expressions for the reconstruction. The importance of multidimensional Riemann solvers in facilitating the design of such schemes is also documented. In two dimensions, a von Neumann analysis of structure-preserving WENO-like schemes that mimetically satisfy the curl constraints, is also presented. It shows the tremendous value of higher order WENO-like schemes in minimizing dissipation and dispersion for this class of problems. Numerical results are also presented to show that the edge-centered curl-preserving (ECCP) schemes meet their design accuracy. This paper is the first paper that invents non-linearly hybridized curl-preserving reconstruction and integrates it with higher order Godunov philosophy. By its very design, this paper is therefore intended to be forward-looking and to set the stage for future work on curl involution-constrained PDEs.

## I.1) Introduction

There has been a lot of emerging interest in mimetic scheme design. These are schemes that preserve structures in the solution that arise from involutions in the governing PDEs. In other words, the PDE itself has some extra symmetries that result in certain features of the solution remaining invariant; and we want the numerical scheme to mimic that.

The simplest example of such involution-constrained PDEs consists of the magnetohydrodynamic (MHD) equations, where Faraday's law ensures divergence-free evolution of the magnetic induction vector field. Another prominent example consists of computational electrodynamics (CED) – the numerical solution of Maxwell's equations – where the divergences of the magnetic induction vector field and the electric displacement vector field are held zero as long as radiation does not interact with a conductor. Numerous papers have been written on these topics, where it has been realized that the divergence-preserving reconstruction of vector fields is an important building block for scheme design and adaptive mesh refinement (AMR) (Balsara and Spicer [3], Balsara [5], [6], [8], Balsara and Dumbser [14], Xu et al. [57], Balsara and Käppeli [20], [24], Balsara *et al*. [10], [18], [22], [23], [25], [26] and Hazra *et al*. [41]). To get fully constraint-preserving, mimetic, time-evolution it was realized that certain update variables have to be collocated at certain favored locations on a mesh. A weighted essentially non-oscillatory, i.e. WENO-like, reconstruction strategy that preserves the divergence of the vector fields was found to be very valuable in extending these methods to high order. We call such methods WENO-like because they draw on many insights from WENO schemes, however these insights are applied at the faces of the mesh while a divergence constraint-preserving reconstruction is used to obtain the vector field within the volume of the mesh. (It should also be noted that two-dimensional WENO reconstruction is used in the faces as part of the three-dimensional divergence constraint-preserving reconstruction, which is why we call the reconstruction WENO-like.) To update such variables in a properly upwinded fashion, it is also crucially important to invoke a multidimensional Riemann solver at the edges of the computational mesh (Balsara, [9], [10], [13], Balsara, Dumbser and Abgrall [12], Balsara and Dumbser [15], Balsara *et al*. [17], Balsara and Nkonga [21]). The multidimensional Riemann solver is, therefore, the other important building block of such divergence-constrained schemes.

In obtaining highly accurate globally constraint-preserving Discontinuous Galerkin-like (DG-like) schemes for MHD and CED, Balsara and Käppeli [20], [24] showed that both these building blocks were crucially important. They showed that if one attempts to bypass either of these building blocks, it will result in an unstable DG-like scheme. We call these schemes DG-like because they evolve all the face-based modes of the vector field so as to ensure globally divergence constraint-preserving evolution of the vector field; however, they are not like classical DG schemes because the modes are not defined on the volumes. We, therefore, see that a study of the involution-preserving reconstruction can provide substantial insights into scheme design. A review of globally divergence constraint-preserving DG schemes for CED is also available in Balsara and Simpson [27] which collects all the ideas together in an easily accessible format in one place.

While MHD and CED are relatively well-studied PDEs with a divergence constraint, a new class of PDEs has recently emerged and their involution constraints are equally interesting. We are referring to PDEs that support curl-free (or curl-preserving) evolution of vector fields. Indeed, the evolution is curl-free in these systems only as long the source terms in the governing equations are zero. Numerous PDEs of great practical interest fall in this category. Many of the hyperbolic systems resulting from the Godunov-Peshkov-Romenski (GPR) formulation for hyperelasticity and compressible multiphase flow with and without surface tension have such curl-preserving update equations (Godunov and Romenski [39], Romenski [47], Romenski *et al*. [48], Peshkov and Romenski [45], Dumbser *et al*. [33], [34], Schmidmayer *et al*. [49]). The equations of General Relativity when cast in the FO-CCZ4 formulation also have such a structure (Alic *et al*. [1], [2], 2012, Brown *et al*. [31], Dumbser *et al*. [35], Dumbser, *et al*. [36]). Similarly, it has recently become possible to recast Schrödinger's equation in first order hyperbolic form, and the time-evolution of this very important equation also has curl-preserving constraints (Dhaouadi *et al*. [32]). As with the divergence-preserving reconstruction, the curl-preserving reconstruction also plays an important role in guiding scheme design. The goal of this paper is to show how curl-preserving reconstruction of vector fields can be carried out and why it is so important in the design of curl-constraint-preserving schemes. We restrict our focus to structured meshes, since the treatment of unstructured meshes will be the topic of another paper.

In this work we take on the task of designing a WENO-like globally curl constraint-preserving reconstruction. This means that the curl of a vector field, evaluated over any closed loop, is always either zero or equal to a specified divergence-free vector field. (As we shall see in Sections II and III, one-dimensional WENO reconstruction is indeed used in the edges of the mesh as part of the three-dimensional curl constraint-preserving reconstruction which is why we think it is appropriate the call the reconstruction WENO-like. Note though that the constrained reconstruction has many further attributes that go beyond the basic WENO reconstruction.) It may even prove advantageous to refer to PDEs that keep the curl exactly zero as curl-free; whereas PDEs that only preserve the curl constraint in certain limits can be referred to as curl-preserving. Some families of involutionary PDEs, like the FO-CCZ4 formulation of the equations of general relativity, can guarantee that certain vector fields remain curl-free for all time. Other involutionary PDE systems, like the hyperbolic formulation of thermal conduction and viscosity and elastic-plastic transition, have vector fields that are only curl-free when the source term is zero. But important interactions with matter, like the use of thermal conduction, or viscosity or elastic-plastic deformation require the operation of non-zero stiff source terms in those PDE systems.

Having a globally constraint-preserving reconstruction in hand is very useful for computational problems for a very important reason:- When coupled with a three-dimensional Riemann solver (see Balsara [16]) it enables us to define a curl-preserving scheme over a single control volume. In other words, the fluid variables can be zone-centered and the curl-free vector field can share that same control volume. The primal curl-constraint-preserving vector field variables of such a scheme reside in the edges of that control volume. Specifically, for a Cartesian mesh, the x-components of such a vector field are collocated at the x-edges; the y-components of such a vector field are collocated at the y-edges and the z-components of such a vector field are collocated at the z-edges. The three-dimensional Riemann solver, invoked at the vertices of the three-dimensional mesh, then yields the curl constraint-preserving update.

### I.2) Introduction to a Sample Curl-Preserving PDE and Motivation for Curl-Preserving Schemes

There are several PDE systems that have curl-free, and curl-preserving, vector fields that arise from involutions in the differential equations. In fact, the entropy-consistent GPR formulation seems to churn out such involutionary PDEs with amazing regularity. But numerical implementations of other valuable PDE systems, like the numerical solution of first order reductions of the Einstein equations or the Schrödinger equation, also result in such constraints. Let us take a simple example involving a fluid with thermal conduction in the GPR formulation. Let us denote the density by $\rho$, the fluid velocity by $\mathbf{v}$, the fluid pressure by "$P$", the fluid temperature by "$T$", the internal thermal energy density by "$e$", the total energy density by $E \equiv e + \rho \mathbf{v}^2/2$, the thermal impulse by a vector $\mathbf{J}$, the heat flux by a vector $\mathbf{q}$ and the thermal stress by the second rank tensor $\boldsymbol{\sigma}$. The equations for a fluid with thermal conduction can be written as

$$\frac{\partial \rho}{\partial t} + \nabla \cdot (\rho \mathbf{v}) = 0 \tag{1.1a}$$

$$\frac{\partial (\rho \mathbf{v})}{\partial t} + \nabla \cdot (\rho \mathbf{v} \otimes \mathbf{v} + P\mathbf{I} + \boldsymbol{\sigma}) = 0 \tag{1.1b}$$

$$\frac{\partial E}{\partial t} + \nabla \cdot ((E+P)\mathbf{v} + \mathbf{v} \cdot \boldsymbol{\sigma} + \mathbf{q}) = 0 \tag{1.1c}$$

$$\frac{\partial \mathbf{J}}{\partial t} + \nabla (\mathbf{J} \cdot \mathbf{v} + T) - \mathbf{v} \times (\nabla \times \mathbf{J}) = -\frac{\rho T}{\tau} \mathbf{J} \tag{1.1d}$$

The identity matrix is denoted by $\mathbf{I}$ in the above equations. The first three of the four equations in eqn. (1.1) above reveal themselves to be the equations for mass, momentum and energy conservation for a fluid, with additional contributions from the thermal conduction vector, $\mathbf{q}$, and the thermal stress tensor, $\boldsymbol{\sigma}$. The fourth equation in eqn. (1.1) is a novel contribution from the GPR formulation, see (Romenski [47]). We see that it will be strongly dependent on the magnitude of the relaxation time $\tau$. When the relaxation term becomes stiff, i.e. the relaxation time is short, the heat conduction will behave asymptotically like the classical Fourier law for parabolic heat conduction. When the relaxation time is very large, the source term becomes irrelevant and the heat conduction will be described by purely hyperbolic heat waves or phonons, propagating with a characteristic speed that is called the second sound. The beauty of the above equations stems from the fact that they constitute a first order hyperbolic system with a source term that may indeed become stiff in certain limits. Therefore, all of the well-developed technologies that have been developed for solving hyperbolic PDE systems with stiff source terms can indeed be brought to bear on the numerical solution of the above PDE system. Furthermore, the solution method does not require the treatment of a parabolic sub-system, which can be computationally expensive. As already stated before, a formal asymptotic analysis of eqn. (1.1) shows that the above equations retrieve the Navier Stokes equations with the traditional Fourier law of heat conduction in the stiff limit when the relaxation time $\tau$ tends to zero. To complete our description of the above system, we also mention the constitutive relation for the thermal stress tensor $\sigma_{ij} = \rho c_h^2 J_i J_j$ and the other constitutive relation for the thermal conduction vector $q_i = \rho T c_h^2 J_i$ . Here $c_h$ denotes the hyperbolic speed of heat waves, i.e. the second sound.

Now let us focus on the last equation in eqn. (1.1). Let us consider the limit where the relaxation time is very large, so that the source term is irrelevant. Since the vector field **J** starts off curl-free, it is easy to see that it remains curl-free by considering the remaining two parts of that equation. The first part of the update equation, given by $\nabla(\mathbf{J} \cdot \mathbf{v} + T)$, is just the gradient of a scalar. Since the curl of a gradient is zero, the first term will not contribute to the curl if none is present initially. The second part of the update equation, given by $\mathbf{v} \times (\nabla \times \mathbf{J})$ , will also be zero if the vector **J** is initially curl-free. We see, therefore, that the vector field **J** stays curl-free if it is initially curl-free in the limit of very large relaxation time. Of course, when the relaxation time cannot be ignored, the curl of the vector field does indeed evolve in response to the presence of the stiff source term $-\rho T \, \mathbf{J}/\tau$ . It is important to realize that if the fourth equation in eqn. (1.1) does not have a consistent discretization then the curl of the vector field **J** will only be specified by the accuracy of the numerical method. As a result, even for regions of the flow that should have no thermal conduction, there will indeed continue to be some small amount of thermal conduction. This affects the fidelity of the method and its results. We, therefore, see the importance of a consistent, curl-preserving discretization and evolution strategy. Let us address that next.

The fourth equation in eqn. (1.1) contains the involution and therefore deserves further attention. Because the evolution of a curl-free vector field **J** is only governed by the gradient term $\nabla(\mathbf{J} \cdot \mathbf{v} + T)$ we must pick a mimetic discretization that ensures this curl-free evolution. A good conceptual model for a curl is the altitude in a mountainous region. It does not matter which closed curve one takes in that mountainous region, as long as the curve is closed, the total change in elevation will be zero. This is the model that we keep at the back of our mind when studying this problem. The closed curve could be the edges of a rectangular mesh in 2D. For a 3D Cartesian mesh, we have closed curves in all the faces of a cuboidal element. Along each of those faces, the circulation of the vector field $\mathbf{J}$ must be zero in all the situations where the vector field is required to evolve in a curl-free fashion. This is only guaranteed if the components of $\mathbf{J}$ are collocated at the edges of the mesh and $(\mathbf{J} \cdot \mathbf{v} + T)$ is collocated at the corners of the mesh. But realize that we are solving a hyperbolic system, as a result $(\mathbf{J} \cdot \mathbf{v} + T)$ will have to be obtained consistently with multidimensional upwinding at the corners of the mesh. There already exists a 3D Riemann solver that does this (Balsara [16]).

Let us now press on with our study of the last equation in eqn. (1.1). Let us focus on the term $\mathbf{v} \times (\nabla \times \mathbf{J})$. When the curl is zero, it is irrelevant. However, when the curl is non-zero, it does affect the time rate of change of the component of $\mathbf{J}$ that is aligned with the edges of the mesh. How can we get the measure of the curl of a vector field? In three-dimensions, we can only do that by reconstructing the vector field in a three-dimensional fashion. (Likewise, of course, in two dimensions!) In other words, we need to start with the components of $\mathbf{J}$ in the edges that surround each volume element of the mesh and obtain from it a consistent value of $\mathbf{J}$ within the volume element. This should be done in a way that reflects, in some appropriate fashion, the curl that is already present in the faces of that mesh element. This is the problem of reconstructing a vector field consistent with its constraints. We, therefore, see that we will have to pay special attention in this paper to curl-free and curl-preserving reconstruction of vector fields.

Let us now take the curl of the third equation in eqn. (1.1). Let us also make the definition $\mathbf{R} \equiv \nabla \times \mathbf{J}$, where the vector field $\mathbf{R}$ is referred to as the Burger's vector field. It is easy to see from its very definition that $\nabla \cdot \mathbf{R} = 0$; i.e. the Burger's vector is divergence-free. Now let us take the curl of the fourth equation in eqn. (1.1). We get (see Peshkov *et al.* [46])

$$\frac{\partial \mathbf{R}}{\partial t} - \nabla \times (\mathbf{v} \times \mathbf{R}) + \nabla \times \left( \frac{\rho T}{\tau} \mathbf{J} \right) = 0 \qquad (1.2)$$

We see immediately that eqn. (1.2) guarantees divergence-free time-evolution for the Burger's vector field. (Those who are familiar with the induction equation for MHD will also see the great parallels between eqn (1.2), evaluated in the limit where $\tau \to \infty$, and the MHD induction equation.)

Now notice that the fourth equation of the set of equations in eqn. (1.1) gives an evolutionary equation for the time-evolution of the curl of **J** -- see eqn. (1.2)! To preserve the mimetic nature of the time-evolution, the components of the curl-constrained vector field should be collocated at the edges of the mesh. We realize from eqn. (1.2) that curl-preserving reconstruction of a vector field should always pay attention to the curl of the vector field. The curl of the reconstructed vector field should be divergence-free. Therefore, we see that when designing a DG-like scheme for curl constraint-preserving PDEs we will get an additional equation for the evolution of the curl of the vector field **J** that is of interest. In keeping with the DG philosophy, the higher moments of the curl of this vector field should have its components collocated at the faces of the mesh and evolved in divergence-free fashion. For a finite volume (FV) WENO-like scheme, of course, we have to reconstruct the higher moments of the Burger's vector field **R**. Furthermore, these higher moments must guarantee that $\nabla \cdot \mathbf{R} = 0$. The zeroth moment of the Burger's vector does not need to be reconstructed because it is always given to us by a discrete application of the definition, $\mathbf{R} \equiv \nabla \times \mathbf{J}$, in the faces of the zone of interest. Therefore, for WENO-like schemes we do not need to simultaneously evolve eqn. (1.2).

**I.3) Introduction : Plan of the Paper**

The rest of the paper follows the ensuing plan. In Section II we show how the curl-preserving reconstruction can be carried out at all locations of a two-dimensional Cartesian mesh; this will include second to fourth order reconstructions. Section III extends these ideas to three-dimensional Cartesian meshes. Those two sections also have demonstrations that (when they are coupled to a multidimensional Riemann solver) the two-dimensional and three-dimensional schemes are multidimensionally upwinded and, therefore, stable. Section IV shows results of a von Neumann stability analysis of curl constraint-preserving WENO-like schemes. Section V shows some results from a couple of model problems where the lack of curl-preserving reconstruction is shown to have obvious deleterious effects. Section VI shows some applications of the GPR system in eqn. (1.1) to some further test problems. Section VII presents some conclusions.

## II) Curl-Preserving Reconstruction on a Two-dimensional Cartesian Mesh

It is easiest to get introduced to this subject in two dimensions, especially on a structured mesh. We consider this problem in five easy Sub-sections. In Sub-section II.1 we present a first order accurate reconstruction of a curl-preserving vector field. In Sub-section II.2 we present a second order accurate reconstruction of a curl-preserving vector field. In Sub-sections II.3 and II.4 we present third and fourth order extensions. Sub-section II.5 shows that the curl-free reconstruction, when combined with a two-dimensional Riemann solver, produces a properly upwinded numerical scheme. Each Sub-section is designed to give us a new and important insight.

**Place Fig. 1 here**

### II.1) Curl-Preserving Reconstruction on a Two-dimensional Cartesian Mesh at First Order

Let us consider what is entailed in a first order reconstruction. In keeping with the spirit of a first order finite volume scheme for fluid flow, it means that each edge of a rectangular/square zone has a component of the vector field along the direction of the edge. In simplest form, and for a unit square zone with extent $(x,y) \in [-1/2, 1/2] \times [-1/2, 1/2]$, this is shown in Fig. 1. Any rectangular zone can be mapped to such a square zone, so our results are perfectly general. Fig. 1 shows the collocation of vector components along the edges of a two-dimensional control volume. As evaluated over the edges of the square element, the discrete circulation is fully specified. (The mean value and its linear variation are shown along each edge in Fig. 1, in anticipation of a second order accurate reconstruction scheme. However, in this Sub-section we ignore the linear variation.) The reconstruction problem for a curl-free reconstruction consists of obtaining a polynomial-based vector field that is globally curl-free within this two-dimensional control volume. The reconstruction problem for a curl-preserving reconstruction consists of obtaining a polynomial-based vector field that matches the specified mean circulation in the zone.

From Fig. 1 we see that the bottom and top x-edges have x-components of the vector field that are given by $V_x^1$ and $V_x^2$ respectively. Likewise, the left and right y-edges have y-components of the vector field that are given by $V_y^1$ and $V_y^2$ respectively. A polynomial that holds over the entire unit square and matches the specified values at the edges is given by

$$V^x(x,y) = V_x^1\left(\frac{1}{2}-y\right) + V_x^2\left(\frac{1}{2}+y\right) \quad ; \quad V^y(x,y) = V_y^1\left(\frac{1}{2}-x\right) + V_y^2\left(\frac{1}{2}+x\right) \tag{2.1}$$

By taking the curl of the above vector field, we get

$$(\nabla \times \mathbf{V})_z = \frac{\partial V^y}{\partial x} - \frac{\partial V^x}{\partial y} = \left[V_x^1 - V_x^2 + V_y^2 - V_y^1\right] \tag{2.2}$$

We see that the curl, evaluated as a differential expression, gives back the discrete circulation of the vector field over the unit square shown in Fig. 1. If the discrete circulation is curl-free then it will evaluate to zero and our vector field in eqn. (2.1) will also be curl-free – i.e. the curl evaluated at each local point in the unit square is exactly zero. If the discrete circulation is not curl-free then the differential form in eqn. (2.2) matches the exact value of the discrete circulation at all locations of the unit square, which is reasonable. Observe too that $V^x(x,y)$ only has linear variation in the y-direction while $V^y(x,y)$ only has linear variation in the x-direction in eqn. (2.1) with the result that the reconstruction in eqn. (2.1) is only first order accurate.

It is also worthwhile to observe that if any three of the four components given by $V_x^1$, $V_x^2$, $V_y^2$ and $V_y^1$ are specified, and if we are told that the vector field is curl-free, then the fourth component is automatically satisfied. This is a small observation for now, but it will be expanded on in subsequent sections.

## II.2) Curl-Preserving Reconstruction on a Two-dimensional Cartesian Mesh at Second Order

Now let us consider second order extensions. In the spirit of van Leer [55] and Kolgan [43], this is tantamount to endowing each of the edges with a piecewise linear variation. From Fig. 1 we see that the bottom and top x-edges have x-components of the vector field that are now endowed with undivided differences in the x-direction given by $(\Delta_x V_x^1)$ and $(\Delta_x V_x^2)$ respectively. In a finite volume setting, these undivided differences can be obtained by using a TVD or WENO scheme applied one-dimensionally along the x-edges. To retain second order accuracy, the one-dimensional reconstruction of the first moments should be obtained from a WENO scheme that is at least second order accurate. Similarly, from Fig. 1 we see that the left and right y-edges have y-components of the vector field that are now endowed with undivided differences in the y-direction given by $\left(\Delta_y V_y^1\right)$ and $\left(\Delta_y V_y^2\right)$ respectively. In a finite volume setting, these undivided differences can be obtained by using a TVD or WENO scheme applied one-dimensionally along the y-edges. (For a second order DG scheme, these undivided differences will indeed become evolutionary modes.) Let us say that we follow exactly the same game-plan as in eqn. (2.1) and write

$$V^x(x,y) = \left(V_x^1 + \left(\Delta_x V_x^1\right)x\right)\left(\frac{1}{2} - y\right) + \left(V_x^2 + \left(\Delta_x V_x^2\right)x\right)\left(\frac{1}{2} + y\right) \;;$$
$$V^y(x,y) = \left(V_y^1 + \left(\Delta_y V_y^1\right)y\right)\left(\frac{1}{2} - x\right) + \left(V_y^2 + \left(\Delta_y V_y^2\right)y\right)\left(\frac{1}{2} + x\right)$$
(2.3)

To see the problem with eqn. (2.3), let us take its curl. We get

$$(\nabla \times \mathbf{V})_z = \frac{\partial V^y}{\partial x} - \frac{\partial V^x}{\partial y}$$
$$= \left[V_x^1 - V_x^2 + V_y^2 - V_y^1\right] + \left[\left(\Delta_x V_x^1\right) - \left(\Delta_x V_x^2\right)\right]x + \left[\left(\Delta_y V_y^2\right) - \left(\Delta_y V_y^1\right)\right]y$$
(2.4)

We see that even if the original vector field had a discrete circulation that was zero over the square shown in Fig. 1, the resulting curl evaluated at all points within the square will not be zero. This is because in general $\left(\Delta_x V_x^2\right) \neq \left(\Delta_x V_x^1\right)$ and $\left(\Delta_y V_y^2\right) \neq \left(\Delta_y V_y^1\right)$, so the linear variations in the x- and y-directions in eqn. (2.4) will not be zero. Therefore, eqn. (2.3) is not curl-preserving.

Having seen that a naïve attack on the problem yields nothing of value, let us renew our effort. We take inspiration from the divergence-free reconstruction of two-dimensional vector fields that was discussed in Sub-section III.1 of Balsara [5] and realize that when we are dealing with a constrained vector field, the components couple. In other words, the vector field is an entire entity and we cannot take the individual components as disjoint entities. Therefore, the x-component of the vector field will couple to the y-component of the vector field so as to preserve the constraints. Since we have already realized that curl-free and curl-preserving reconstruction are just two sides of the same coin, we focus on the former problem first. Let us write our vector field as

$$V^x(x,y) = \left[V_x^1 + (\Delta_x V_x^1)x\right]\left(\frac{1}{2} - y\right) + \left[V_x^2 + (\Delta_x V_x^2)x\right]\left(\frac{1}{2} + y\right) + a_{yy}(1 - 4y^2) \ ;$$

$$V^y(x,y) = \left[V_y^1 + (\Delta_y V_y^1)y\right]\left(\frac{1}{2} - x\right) + \left[V_y^2 + (\Delta_y V_y^2)y\right]\left(\frac{1}{2} + x\right) + b_{xx}(1 - 4x^2)$$

(2.5)

Notice that all the terms that are needed for obtaining second order accuracy are already present in eqn. (2.5). The $a_{yy}(1-4y^2)$ term is designed to go to zero at $y = \pm 1/2$, i.e. at the upper and lower x-edges of the mesh shown in Fig. 1. This ensures that at the abutting x-edges of a full two-dimensional mesh we have one and the same value for $V^x(x,y)$. Similarly, the $b_{xx}(1-4x^2)$ term is designed to go to zero at $x = \pm 1/2$, i.e. at the right and left y-edges of the mesh shown in Fig. 1. This ensures that the abutting y-edges of a full two-dimensional mesh have one and the same value of $V^y(x,y)$. The above two sentences ensure that the reconstruction strategy is globally curl-free or curl-preserving. We can now express the curl of the above vector field as

$$(\nabla \times \mathbf{V})_z = \frac{\partial V^y}{\partial x} - \frac{\partial V^x}{\partial y}$$

$$= \left[V_x^1 - V_x^2 + V_y^2 - V_y^1\right] + \left[(\Delta_x V_x^1) - (\Delta_x V_x^2) - 8b_{xx}\right]x + \left[(\Delta_y V_y^2) - (\Delta_y V_y^1) + 8a_{yy}\right]y$$

(2.6)

We see that the first square bracket in the above equation still expresses the discrete circulation, which is exactly zero for a curl-free vector field. The second and third square brackets in the above equation can be made zero by setting

$$b_{xx} = \frac{1}{8}\left[(\Delta_x V_x^1) - (\Delta_x V_x^2)\right] \quad ; \quad a_{yy} = \frac{1}{8}\left[(\Delta_y V_y^1) - (\Delta_y V_y^2)\right]$$

(2.7)

Notice that from a finite difference point of view, the coefficients $a_{yy}$ and $b_{xx}$ are just higher order derivatives of the undivided differences, and as a result, the second order accuracy of eqn. (2.5) is not affected by the inclusion of these additional terms.

The analogies with divergence-free reconstruction in Balsara [5] are also worth drawing. In both cases, the first order term is an expression of the discrete constraint applied to the boundaries of the element. The inclusion of higher order terms requires additional coefficients to ensure that the differential form of the constraint is exactly satisfied at all locations within the element.

Now that we have thoroughly discussed all the nuances of a curl-free reconstruction of a vector field, we are in a position to discuss how the idea goes over to a curl-preserving reconstruction. Recall that the fourth equation in eqn. (1.1) indeed has an evolutionary equation for the curl, see eqn. (1.2). Notice that in a curl-preserving reconstruction, Fig. 1 shows us that the discrete circulation in the square zone is given by $\left[ V_x^1 - V_x^2 + V_y^2 - V_y^1 \right]$. For a second order DG scheme, eqn. (1.2) would provide the time-evolving higher moments of the curl. For a FV scheme, we can reconstruct such a quantity for all zones of the two-dimensional mesh. Using neighboring elements, we can obtain a TVD-based or WENO-based piecewise linear, finite volume reconstruction of the circulation. Such a reconstruction should match the discrete circulation in the target zone. As a result, for the zone shown in Fig. 1, we can write the piecewise linear circulation as

$$R^z(x,y) = \left[ V_x^1 - V_x^2 + V_y^2 - V_y^1 \right] + \left( \Delta_x R^z \right) x + \left( \Delta_y R^z \right) y \tag{2.8}$$

Matching eqns. (2.6) and (2.8) we get

$$b_{xx} = \frac{1}{8}\left[ -\left( \Delta_x R^z \right) + \left( \Delta_x V_x^1 \right) - \left( \Delta_x V_x^2 \right) \right] \quad ; \quad a_{yy} = \frac{1}{8}\left[ \left( \Delta_y R^z \right) + \left( \Delta_y V_y^1 \right) - \left( \Delta_y V_y^2 \right) \right] \tag{2.9}$$

Comparing eqns. (2.7) and (2.9) we now observe an obvious correspondence between curl-free and curl-preserving reconstruction. Specifying one is tantamount to specifying the other.

It is also worth pointing out that eqns. (1.1) and (1.2) show us that the curl of the vector field explicitly participates in the time-update. Therefore, it is useful to provide explicit expressions not just for the vector field but also for its curl, as was done in eqn. (2.6). Such expressions prove to be quite valuable for making numerical implementations.

When a one-dimensional TVD or WENO reconstruction is used to obtain $\left( \Delta_x V_x^1 \right)$, $\left( \Delta_x V_x^2 \right)$, $\left( \Delta_y V_y^1 \right)$ and $\left( \Delta_y V_y^2 \right)$ in eqn. (2.5), we realize that we are automatically building in non-linear hybridization into the multidimensional curl-preserving reconstruction. It is for this reason that all the resulting schemes presented in this paper inherit all the good non-linear stabilization properties that are an integral part of higher order Godunov schemes. Consequently, by its very construction, the non-linearly hybridized curl-preserving reconstruction developed here integrates very well with higher order Godunov philosophy.

**II.3) Curl-Preserving Reconstruction on a Two-dimensional Cartesian Mesh at Third Order**

We extend the results from the previous Sub-section to the third order case here. In addition to being useful for scheme design, this is useful for analytic work on WENO-like or DG-like schemes, and their von Neumann stability analysis. Notice first off that $V^x(x,y)$ in eqn. (2.5) has constant, $x$, $y$, $xy$ and $y^2$ terms. Therefore, to become a truly third order reconstruction, it minimally needs an $x^2$–dependent term, which will indeed be added along the x-edges. Similarly, $V^y(x,y)$ in eqn. (2.5) has constant, $x$, $y$, $xy$ and $x^2$ terms. Therefore, to become a truly third order reconstruction, it minimally needs a $y^2$–dependent term, which will indeed be added along the y-edges. Such a way of thinking shows us how each reconstruction of the curl constraint-preserving vector field at a certain order illuminates the way to the reconstruction at the next higher order. At least on a structured mesh, where the polynomial terms can proliferate, this is the systematic strategy that one should pursue.

It is important to be emphatic about a point of detail that we develop in this paragraph. One may think that it is unreasonable to claim that the $y^2$ mode is present in $V^x(x,y)$ in eqn. (2.5) because that mode comes purely from the constraint-satisfaction at second order. Similarly, one may think that it is unreasonable to claim that the $x^2$ mode is present in $V^y(x,y)$ in eqn. (2.5) because that mode also comes from constraint-satisfaction at second order. However, indeed those modes are truly present because this is the very idea behind a constrained vector field. The constraint basically tells us that if a mode in $V^x(x,y)$ or $V^y(x,y)$ is needed in order to satisfy the curl-free (or curl-preserving) constraint, then it is indeed truly satisfied. It does not matter that it is satisfied by variation in the other vector component, because the curl-free (or curl-preserving) vector field is just a single entity. None of the components of the curl-free vector field are entire in themselves, they only exist as parts of a whole! An entirely analogous observation has been found to be true over and over again in divergence constraint-preserving reconstruction for MHD and CED (Balsara [5], [6], [8], Balsara *et al*. [22], [23]).

Let us now extend the curl-free reconstruction to third order. At the bottom and top x-edges of the square shown in Fig. 1 we now add piecewise quadratic modes that we denote by $\left(\Delta_{xx}V_x^1\right)$ and $\left(\Delta_{xx}V_x^2\right)$ respectively. In a finite volume setting, these quadratic modes can be obtained by using a WENO scheme applied one-dimensionally along the x-edges. The linear modes are also provided by the same WENO scheme. To retain third order accuracy, the one-dimensional reconstruction of the linear and quadratic modes should be obtained from a WENO scheme that is at least third order accurate. Similarly, at the left and right y-edges of the square shown in Fig. 1 we now add piecewise quadratic modes that we denote by $\left(\Delta_{yy}V_y^1\right)$ and $\left(\Delta_{yy}V_y^2\right)$ respectively. In a finite volume setting, these quadratic modes can be obtained by using a WENO scheme applied one-dimensionally along the y-edges; likewise for the linear modes. The linear modes are also provided by the same WENO scheme. We use a sequence of orthogonal Legendre polynomials because the higher order polynomials all average to zero and the polynomial sequence retains a nice orthogonality property. It is important to notice that the inclusion of a quadratic $x^2$-dependent term along each of the x-edges in $V^x(x,y)$ will also trigger additional modes of the form $x^2 y$. Likewise, the inclusion of a quadratic $y^2$-dependent term along each of the y-edges in $V^y(x,y)$ will also trigger additional modes of the form $y^2 x$. To compensate for the effect of these terms on the curl, some higher order polynomial terms have to be added. We now have at third order

$$V^x(x,y) = \left[V_x^1 + \left(\Delta_x V_x^1\right)x + \left(\Delta_{xx}V_x^1\right)\left(x^2 - 1/12\right)\right]\left(\frac{1}{2} - y\right)$$
$$+ \left[V_x^2 + \left(\Delta_x V_x^2\right)x + \left(\Delta_{xx}V_x^2\right)\left(x^2 - 1/12\right)\right]\left(\frac{1}{2} + y\right)$$
$$+ a_{yy}\left(1 - 4y^2\right) + a_{yyy}y\left(1 - 4y^2\right) + a_{xyy}x\left(1 - 4y^2\right) \ ;$$

$$V^y(x,y) = \left[V_y^1 + \left(\Delta_y V_y^1\right)y + \left(\Delta_{yy}V_y^1\right)\left(y^2 - 1/12\right)\right]\left(\frac{1}{2} - x\right)$$
$$+ \left[V_y^2 + \left(\Delta_y V_y^2\right)y + \left(\Delta_{yy}V_y^2\right)\left(y^2 - 1/12\right)\right]\left(\frac{1}{2} + x\right)$$
$$+ b_{xx}\left(1 - 4x^2\right) + b_{xxx}x\left(1 - 4x^2\right) + b_{xxy}y\left(1 - 4x^2\right)$$

(2.10)

The $a_{xyy}$ and $b_{xxy}$ are not mandatory for order property preservation, but we shall show shortly that they are needed in the construction of a WENO-like or DG-like scheme for curl constraint-preserving vector fields. We can now write out the curl of the above vector field as

$$(\nabla \times \mathbf{V})_z = \frac{\partial V^y}{\partial x} - \frac{\partial V^x}{\partial y}$$

$$= \left[V_x^1 - V_x^2 + V_y^2 - V_y^1\right] + \left[(\Delta_x V_x^1) - (\Delta_x V_x^2) - 8b_{xx}\right]x + \left[(\Delta_y V_y^2) - (\Delta_y V_y^1) + 8a_{yy}\right]y$$

$$+ \left[(\Delta_{xx} V_x^1) - (\Delta_{xx} V_x^2) - 12 b_{xxx}\right](x^2 - 1/12) + \left[-(\Delta_{yy} V_y^1) + (\Delta_{yy} V_y^2) + 12 a_{yyy}\right](y^2 - 1/12) + 8(a_{xyy} - b_{xxy})xy$$

(2.11)

As with eqn. (2.8), we can now reconstruct the discrete circulation up to and including quadratic variation over each zone, and write the result as

$$R^z(x, y) = \left[V_x^1 - V_x^2 + V_y^2 - V_y^1\right] + (\Delta_x R^z)x + (\Delta_y R^z)y$$
$$+ (\Delta_{xx} R^z)(x^2 - 1/12) + (\Delta_{yy} R^z)(y^2 - 1/12) + (\Delta_{xy} R^z)xy$$

(2.12)

Please note that for a WENO-like scheme, the higher moments in eqn. (2.12) are reconstructed; whereas for a DG-like scheme the higher moments in eqn. (2.12) are evolved. Equating like terms in eqns. (2.11) and (2.12) we get

$$b_{xxx} = \frac{1}{12}\left[-(\Delta_{xx} R^z) + (\Delta_{xx} V_x^1) - (\Delta_{xx} V_x^2)\right] \quad ; \quad a_{yyy} = \frac{1}{12}\left[(\Delta_{yy} R^z) + (\Delta_{yy} V_y^1) - (\Delta_{yy} V_y^2)\right] ;$$

$$b_{xx} = \frac{1}{8}\left[-(\Delta_x R^z) + (\Delta_x V_x^1) - (\Delta_x V_x^2)\right] \quad ; \quad a_{yy} = \frac{1}{8}\left[(\Delta_y R^z) + (\Delta_y V_y^1) - (\Delta_y V_y^2)\right] ;$$

$$a_{xyy} = -b_{xxy} = \frac{1}{16}(\Delta_{xy} R^z)$$

(2.13)

This gives us the third order curl-free or curl-preserving reconstruction on a two-dimensional Cartesian mesh. To get a curl-free reconstruction, just set all the coefficients in eqn. (2.12) to zero. We can now also notice that a third order accurate DG-like scheme which evolves all the modes of the circulation in eqn. (2.12) will indeed evolve a value for $(\Delta_{xy} R^z)$. As a result, the terms $a_{xyy}$ and $b_{xxy}$ in eqn. (2.10) were needed for matching all the modes of a third order accurate DG-like scheme which evolves the primal vector field in eqn. (2.10) as well as its curl in eqn. (2.12). If the vector field in eqn. (2.10) can be guaranteed to be curl-free then the terms $a_{xyy}$ and $b_{xxy}$ in eqn. (2.10) are not needed.

Now notice that the vector field in eqn. (2.10) only needs to be up to third order accurate, i.e. it only needs to retain all the quadratic terms that arise in a two-dimensional Taylor series expansion. Therefore, when dealing with a finite volume scheme, its curl only needs to be second order accurate. In other words, for a third order accurate FV scheme, if we had set the coefficients $\left(\Delta_{xx} R^z\right)$, $\left(\Delta_{yy} R^z\right)$ and $\left(\Delta_{xy} R^z\right)$ to zero, eqn. (2.12) would still have been second order accurate. Therefore, in a FV scheme, it would have been acceptable, and third order accurate, to have set $a_{xyy} = b_{xxy} = 0$ in eqn. (2.10). For a DG-like scheme, of course, all the modes in eqn. (2.12) are needed. Also notice that from a finite difference point of view, the coefficients $a_{yy}$, $a_{yyy}$, $b_{xx}$ and $b_{xxx}$ in eqn. (2.13) are just higher order derivatives of the undivided differences, and as a result, the third order accuracy of eqn. (2.10) is not affected by the inclusion of these additional terms.

## II.4) Curl-Preserving Reconstruction on a Two-dimensional Cartesian Mesh at Fourth Order

Let us now make a fourth order extension. We use our idea of systematically thinking about the terms that are present in the third order reconstruction and using them to inform our choices at fourth order. Notice, first off, that $V^x(x, y)$ in eqn. (2.10) has constant, $x$, $y$, $x^2$, $xy$, $y^2$, $y^3$ and $x^2 y$ terms. To that, along each x-edge, we will indeed add an $x^3$–dependent term. However, to have full fourth order accurate reconstruction, we will still need an $xy^2$ term, which must indeed be added with a zone-centered collocation! In other words, by enriching the moments along each x-edge we simply cannot obtain a term with $xy^2$ variation, so we have to include it at a location where all the moments have validity, namely at the zone center. Furthermore, notice that $V^y(x, y)$ in eqn. (2.10) has constant, $x$, $y$, $x^2$, $xy$, $y^2$, $x^3$ and $xy^2$ terms. To that, along each y-edge, we will indeed add a $y^3$– dependent term. However, to have full fourth order accurate reconstruction, we will still need an $x^2 y$ term, which must indeed be added with a volume-centered collocation! As before, by enriching the moments along each y-edge we simply cannot obtain a term with $x^2 y$ variation, so we have to include it at a location where all the moments have validity, namely at the zone center. We now see the value of our systematic, order-by-order approach because it has highlighted for us the fourth order terms that are supplied by enriching the basis space along the edges and the additional modes that have to be supplied volumetrically. (A similar subdivision occurs in divergence-free reconstruction for CED and MHD where we already know that at fourth order and beyond, many of the modes are face-centered, but some are volume-centered. There too, we realized that we could not enrich the space of spatial modes to obtain all the terms that are needed in a fourth order accurate Taylor series expansion. As a result, some of the modes had to be volume-centered; see Balsara and Käppeli [24] and Hazra *et al*. [41].)

Let us now extend the curl-free reconstruction to fourth order. At the bottom and top x-edges of the square shown in Fig. 1 we now add piecewise cubic modes that we denote by $\left(\Delta_{xxx} V_x^1\right)$ and $\left(\Delta_{xxx} V_x^2\right)$ respectively. Similarly, at the left and right y-edges of the square shown in Fig. 1 we now add piecewise cubic modes that we denote by $\left(\Delta_{yyy} V_y^1\right)$ and $\left(\Delta_{yyy} V_y^2\right)$ respectively. We use a sequence of orthogonal Legendre polynomials, as before. The inclusion of a cubic $x^3$–dependent term along each of the x-edges in $V^x(x, y)$ will also trigger additional modes of the form $x^3 y$. Likewise, the inclusion of a cubic $y^3$-dependent term along each of the y-edges in $V^y(x, y)$ will also trigger additional modes of the form $y^3 x$. To compensate for the effect of these terms on the curl, some higher order polynomial terms have to be added. The analysis from the previous paragraph allows us to write the fourth order curl constraint-preserving vector field as

$$V^x(x,y) = \left[V_x^1 + (\Delta_x V_x^1)x + (\Delta_{xx}V_x^1)(x^2 - 1/12) + (\Delta_{xxx}V_x^1)(x^3 - 3x/20)\right]\left(\frac{1}{2} - y\right)$$

$$+ \left[V_x^2 + (\Delta_x V_x^2)x + (\Delta_{xx}V_x^2)(x^2 - 1/12) + (\Delta_{xxx}V_x^2)(x^3 - 3x/20)\right]\left(\frac{1}{2} + y\right) + a_{xyy}x(1-4y^2)$$

$$+ a_{yy}(1-4y^2) + a_{yyy}y(1-4y^2) + a_{yyyy}y^2(1-4y^2) + a_{xyyy}xy(1-4y^2) + a_{xxyy}(x^2 - 1/12)(1-4y^2)$$

$$+ a_{xyyyy}x\left((y^4 - 3y^2/14 + 3/560) - 1/70\right) + a_{xxxyy}(x^3 - 3x/20)(1 - 4y^2) \; ;$$

$$V^y(x,y) = \left[V_y^1 + (\Delta_y V_y^1)y + (\Delta_{yy}V_y^1)(y^2 - 1/12) + (\Delta_{yyy}V_y^1)(y^3 - 3y/20)\right]\left(\frac{1}{2} - x\right)$$

$$+ \left[V_y^2 + (\Delta_y V_y^2)y + (\Delta_{yy}V_y^2)(y^2 - 1/12) + (\Delta_{yyy}V_y^2)(y^3 - 3y/20)\right]\left(\frac{1}{2} + x\right) + b_{xxy}y(1-4x^2)$$

$$+ b_{xx}(1-4x^2) + b_{xxx}x(1-4x^2) + b_{xxxx}x^2(1-4x^2) + b_{xxxy}xy(1-4x^2) + b_{xxyy}(y^2 - 1/12)(1 - 4x^2)$$

$$+ b_{xxxxy}y\left((x^4 - 3x^2/14 + 3/560) - 1/70\right) + b_{xxyyy}(y^3 - 3y/20)(1 - 4x^2)$$

(2.14)

Note that the modes $a_{xyy}$ and $b_{xxy}$ correspond to zone-centered modes that carry the $xy^2$ and $x^2 y$ variation. There are 16 coefficients in the above equations and we need a strategy for fixing them up. The curl now becomes

$$(\nabla \times \mathbf{V})_z = \frac{\partial V^y}{\partial x} - \frac{\partial V^x}{\partial y} = \left[V_x^1 - V_x^2 + V_y^2 - V_y^1\right] + \left[(\Delta_x V_x^1) - (\Delta_x V_x^2) - 8b_{xx} - 2b_{xxxx}/5\right]x$$

$$+ \left[(\Delta_y V_y^2) - (\Delta_y V_y^1) + 8a_{yy} + 2a_{yyyy}/5\right]y$$

$$+ \left[(\Delta_{xx}V_x^1) - (\Delta_{xx}V_x^2) - 12b_{xxx}\right](x^2 - 1/12) + \left[-(\Delta_{yy}V_y^1) + (\Delta_{yy}V_y^2) + 12a_{yyy}\right](y^2 - 1/12)$$

$$+ \left[8(a_{xyy} - b_{xxy}) - 6(a_{xyyyy} - b_{xxxxy})/35\right]xy$$

$$+ \left[(\Delta_{xxx}V_x^1) - (\Delta_{xxx}V_x^2) - 16b_{xxxx}\right](x^3 - 3x/20) + \left[-(\Delta_{yyy}V_y^1) + (\Delta_{yyy}V_y^2) + 16a_{yyyy}\right](y^3 - 3y/20)$$

$$+ \left[8a_{xxyy} - 12b_{xxxy}\right]y(x^2 - 1/12) + \left[-8b_{xxyy} + 12a_{xyyy}\right]x(y^2 - 1/12)$$

$$+ \left[8a_{xxxyy} + 4b_{xxxxy}\right](x^3 - 3x/20)y + \left[-8b_{xxyyy} - 4a_{xyyyy}\right](y^3 - 3y/20)x$$

(2.15)

As with eqn. (2.12), we can now reconstruct the discrete circulation up to and including quadratic variation over each zone and including the minimal number of cubic modes that arise in eqn. (2.15), and write the result as

$$\begin{aligned}
R^z(x,y) &= \left[V_x^1 - V_x^2 + V_y^2 - V_y^1\right] + \left(\Delta_x R^z\right)x + \left(\Delta_y R^z\right)y \\
&\quad + \left(\Delta_{xx} R^z\right)\left(x^2 - 1/12\right) + \left(\Delta_{yy} R^z\right)\left(y^2 - 1/12\right) + \left(\Delta_{xy} R^z\right)xy \\
&\quad + \left(\Delta_{xxx} R^z\right)\left(x^3 - 3x/20\right) + \left(\Delta_{yyy} R^z\right)\left(y^3 - 3y/20\right) + \left(\Delta_{xxy} R^z\right)y\left(x^2 - 1/12\right) \\
&\quad + \left(\Delta_{xyy} R^z\right)x\left(y^2 - 1/12\right)
\end{aligned} \quad (2.16)$$

Equating like terms in eqns. (2.15) and (2.16) gives us 11 equations. Therefore, we see that at fourth (and higher) orders the curl constraints, by themselves, do not give us sufficient information for uniquely fixing up the 16 coefficients in eqn. (2.14). However, we expect a solution which has 5 free parameters that can be fixed by using some other logic, which we will describe shortly. The parametric solution can be written quite easily if we introduce a parameter "$\alpha$" defined by $\alpha = a_{xxyy} + b_{xxyy}$. We can then write the entire solution in terms of the 5 parameters $a_{xyy}$, $b_{xxy}$, $a_{xxyy}$, $b_{xxyy}$ and $\alpha$ as follows:-

$$a_{xxxy} = \frac{1}{24}\left[12\alpha + 280 a_{xyy} - 280 b_{xxy} - 35\left(\Delta_{xy} R^z\right)\right] \quad ; \quad b_{xyyy} = \frac{1}{24}\left[12\alpha - 280 a_{xyy} + 280 b_{xxy} + 35\left(\Delta_{xy} R^z\right)\right] \quad ;$$

$$a_{xyyy} = -2 b_{xxyy} \quad ; \quad b_{xxxy} = -2 a_{xxyy} \quad ;$$

$$a_{yyyy} = \frac{1}{16}\left[\left(\Delta_{yyy} R^z\right) + \left(\Delta_{yyy} V_y^1\right) - \left(\Delta_{yyy} V_y^2\right)\right] \quad ; \quad b_{xxxx} = \frac{1}{16}\left[-\left(\Delta_{xxx} R^z\right) + \left(\Delta_{xxx} V_x^1\right) - \left(\Delta_{xxx} V_x^2\right)\right] \quad ;$$

$$a_{xyy} = \frac{1}{12}\left[8 b_{xxyy} + \left(\Delta_{yy} R^z\right)\right] \quad ; \quad b_{xxy} = \frac{1}{12}\left[8 a_{xxyy} - \left(\Delta_{xy} R^z\right)\right] \quad ;$$

$$a_{yyy} = \frac{1}{12}\left[\left(\Delta_{yy} R^z\right) + \left(\Delta_{yy} V_y^1\right) - \left(\Delta_{yy} V_y^2\right)\right] \quad ; \quad b_{xxx} = \frac{1}{12}\left[-\left(\Delta_{xx} R^z\right) + \left(\Delta_{xx} V_x^1\right) - \left(\Delta_{xx} V_x^2\right)\right] \quad ;$$

$$a_{yy} = \frac{1}{320}\left[40\left(\Delta_y V_y^1\right) - 40\left(\Delta_y V_y^2\right) - \left(\Delta_{yyy} V_y^1\right) + \left(\Delta_{yyy} V_y^2\right) + 40\left(\Delta_y R^z\right) - \left(\Delta_{yyy} R^z\right)\right] \quad ;$$

$$b_{xx} = \frac{1}{320}\left[40\left(\Delta_x V_x^1\right) - 40\left(\Delta_x V_x^2\right) - \left(\Delta_{xxx} V_x^1\right) + \left(\Delta_{xxx} V_x^2\right) - 40\left(\Delta_x R^z\right) + \left(\Delta_{xxx} R^z\right)\right]$$

$$(2.17)$$

Up to this point in the narrative, we have left the 5 parameters $a_{xyy}$, $b_{xxy}$, $a_{xxyy}$, $b_{xxyy}$ and $\alpha$ undetermined. They will only be fixed after we make the following consideration. The vector field in eqn. (2.14) includes modes that reside on the edges of the mesh and modes that are zone-centered. It is, therefore, interesting to ask how both kinds of modes (edge-centered and zone-centered) can be accommodated seamlessly in a WENO-like or DG-like scheme? Indeed, we get a valuable new insight by addressing this question. Realize that the vector field in eqn. (2.14) can also be decomposed in terms of orthogonal Legendre polynomials as follows:-

$$V^x(x,y) = \left[(V_x^1 + V_x^2)/2 + ((\Delta_y R^z) + (\Delta_y V_y^1) - (\Delta_y V_y^2))/12\right]$$
$$+ \left[((\Delta_x V_x^1) + (\Delta_x V_x^2))/2 + (\Delta_{xy} R^z)/24 + \alpha/70 + (a_{xyy} + b_{xxy})/3\right] x$$
$$+ \left[-V_x^1 + V_x^2 + ((\Delta_{yy} R^z) + (\Delta_{yy} V_y^1) - (\Delta_{yy} V_y^2))/30\right] y$$
$$+ \left[2a_{xxyy}/3 + ((\Delta_{xx} V_x^1) + (\Delta_{xx} V_x^2))/2\right](x^2 - 1/12)$$
$$+ \left[\begin{array}{l}((\Delta_{yyy} R^z) + (\Delta_{yyy} V_y^1) - (\Delta_{yyy} V_y^2))/112 \\ -(40(\Delta_y V_y^1) - 40(\Delta_y V_y^2) - (\Delta_{yyy} V_y^1) + (\Delta_{yyy} V_y^2) + 40(\Delta_y R^z) - (\Delta_{yyy} R^z))/80\end{array}\right](y^2 - 1/12)$$
$$+ \left[-(\Delta_x V_x^1) + (\Delta_x V_x^2) + (8b_{xxy} + (\Delta_{xy} R^z))/30\right] xy$$
$$+ \left[((\Delta_{xxx} V_x^1) + (\Delta_{xxx} V_x^2))/2 + (12\alpha + 280a_{xyy} - 280b_{xxy} - 35(\Delta_{xy} R^z))/36\right](x^3 - 3x/20)$$
$$+ \left[-((\Delta_{yy} R^z) + (\Delta_{yy} V_y^1) - (\Delta_{yy} V_y^2))/3\right](y^3 - 3y/20)$$
$$+ \left[-(\Delta_{xx} V_x^1) + (\Delta_{xx} V_x^2)\right](x^2 - 1/12) y + \left[-4a_{xyy}\right](y^2 - 1/12) x$$
$$+ \left[-((\Delta_{yyy} R^z) + (\Delta_{yyy} V_y^1) - (\Delta_{yyy} V_y^2))/4\right](y^4 - 3y^2/14 + 3/560) + \left[-(\Delta_{xxx} V_x^1) + (\Delta_{xxx} V_x^2)\right](x^3 - 3x/20) y$$
$$+ \left[-(8b_{xxy} + (\Delta_{xy} R^z))/3\right](y^3 - 3y/20) x + \left[-4a_{xxyy}\right](x^2 - 1/12)(y^2 - 1/12)$$
$$+ \left[-(12\alpha - 280a_{xyy} + 280b_{xxy} + 35(\Delta_{xy} R^z))/12\right] x(y^4 - 3y^2/14 + 3/560)$$
$$+ \left[-(12\alpha + 280a_{xyy} - 280b_{xxy} - 35(\Delta_{xy} R^z))/6\right](x^3 - 3x/20)(y^2 - 1/12)$$

(2.18a)

$$V^y(x,y) = \left[(V_y^1 + V_y^2)/2 + (-(\Delta_x R^z) + (\Delta_x V_x^1) - (\Delta_x V_x^2))/12\right]$$

$$+ \left[-V_y^1 + V_y^2 + (-(\Delta_{xx} R^z) + (\Delta_{xx} V_x^1) - (\Delta_{xx} V_x^2))/30\right] x$$

$$+ \left[((\Delta_y V_y^1) + (\Delta_y V_y^2))/2 - (\Delta_{xy} R^z)/24 + \alpha/70 + (a_{xyy} + b_{xxy})/3\right] y$$

$$+ \begin{bmatrix} (-(\Delta_{xxx} R^z) + (\Delta_{xxx} V_x^1) - (\Delta_{xxx} V_x^2))/122 \\ -(40(\Delta_x V_x^1) - 40(\Delta_x V_x^2) - (\Delta_{xxx} V_x^1) + (\Delta_{xxx} V_x^2) - 40(\Delta_x R^z) + (\Delta_{xxx} R^z))/80 \end{bmatrix} (x^2 - 1/12)$$

$$+ \left[2b_{xxyy}/3 + ((\Delta_{yy} V_y^1) + (\Delta_{yy} V_y^2))/2\right](y^2 - 1/12) + \left[-(\Delta_y V_y^1) + (\Delta_y V_y^2) + (8a_{xyy} - (\Delta_{xy} R^z))/30\right] xy$$

$$+ \left[-(-(\Delta_{xx} R^z) + (\Delta_{xx} V_x^1) - (\Delta_{xx} V_x^2))/3\right](x^3 - 3x/20)$$

$$+ \left[((\Delta_{yyy} V_y^1) + (\Delta_{yyy} V_y^2))/2 + (12\alpha - 280 a_{xyy} + 280 b_{xxy} + 35(\Delta_{xy} R^z))/36\right](y^3 - 3y/20)$$

$$+ \left[-4b_{xxy}\right](x^2 - 1/12) y + \left[-(\Delta_{yy} V_y^1) + (\Delta_{yy} V_y^2)\right](y^2 - 1/12) x$$

$$+ \left[-(-(\Delta_{xxx} R^z) + (\Delta_{xxx} V_x^1) - (\Delta_{xxx} V_x^2))/4\right](x^4 - 3x^2/14 + 3/560) + \left[-(8a_{xyy} - (\Delta_{xy} R^z))/3\right](x^3 - 3x/20) y$$

$$+ \left[-(\Delta_{yyy} V_y^1) + (\Delta_{yyy} V_y^2)\right](y^3 - 3y/20) x + \left[-4b_{xxyy}\right](x^2 - 1/12)(y^2 - 1/12)$$

$$+ \left[-(12\alpha + 280 a_{xyy} - 280 b_{xxy} - 35(\Delta_{xy} R^z))/12\right] y(x^4 - 3x^2/14 + 3/560)$$

$$+ \left[-(12\alpha - 280 a_{xyy} + 280 b_{xxy} + 35(\Delta_{xy} R^z))/6\right](x^2 - 1/12)(y^3 - 3y/20)$$

(2.18b)

The above two equations give us all the modes that are present in the curl-preserving reconstruction of the vector field. Notice that the 2nd, 4th, 6th, 7th and 10th terms in $V^x(x, y)$ as well as the 3rd, 5th, 6th, 8th and 9th terms in $V^y(x, y)$ are only specified parametrically in terms of the 5 parameters $a_{xyy}$, $b_{xxy}$, $a_{xxyy}$, $b_{xxyy}$ and $\alpha$. In other words, the edge-centered vector field components do not fix those modes. Some additional information should be gleaned from the zone-centered formulation of the hyperbolic system. (Interestingly, this is also true for fourth and higher order divergence-preserving reconstruction in CED.)

First, notice that even without specifying $a_{xyy}$, $b_{xxy}$, $a_{xxyy}$, $b_{xxyy}$ and $\alpha$ we can still use eqn. (2.18) to write the zone-averaged values for the two vector components as:-

$$\langle V^x(x,y)\rangle = (V_x^1 + V_x^2)/2 + ((\Delta_y V_y^1) - (\Delta_y V_y^2) + (\Delta_y R^z))/12 \ ;$$
$$\langle V^y(x,y)\rangle = (V_y^1 + V_y^2)/2 + ((\Delta_x V_x^1) - (\Delta_x V_x^2) - (\Delta_x R^z))/12 \quad (2.19)$$

Therefore, the zone-averages of the vector field are indeed fully specified by the edge values and the higher order moments in the edges.

Next, notice that the curl-constrained, edge-based vector components have not fully pinned down the 2$^{nd}$, 4$^{th}$, 6$^{th}$, 7$^{th}$ and 10$^{th}$ terms in $V^x(x,y)$; nor have they fully pinned down the 3$^{rd}$, 5$^{th}$, 6$^{th}$, 8$^{th}$ and 9$^{th}$ terms in $V^y(x,y)$. Therefore, owing to the mass matrix being diagonal for the expansion in eqns. (2.18a) and (2.18b), we can indeed use a traditional, finite-volume-based DG scheme to evolve the first 10 moments in those equations. Such a DG scheme would be evolved as an auxiliary scheme to the edge-based DG-like scheme and is intended to help pin down the moments that cannot be specified exclusively by the edge-based DG-like scheme. Let us denote those modes with a superscript of "FV". Since we seek the maximum concordance between the zone-centered DG scheme and the edge-centered DG-like scheme. This is done by requiring that the following 10 linear equations should have their residuals minimized in a least squares sense:-

$$((\Delta_x V_x^1) + (\Delta_x V_x^2))/2 + (\Delta_{xy} R^z)/24 + \alpha/70 + (a_{yy} + b_{xy})/3 = (\Delta_x V_x^{FV}) \ ;$$
$$2a_{xyy}/3 + ((\Delta_{xx} V_x^1) + (\Delta_{xx} V_x^2))/2 = (\Delta_{xx} V_x^{FV}) \ ;$$
$$-(\Delta_y V_x^1) + (\Delta_y V_x^2) + (8b_{xyy} + (\Delta_{yy} R^z))/30 = (\Delta_{xy} V_x^{FV}) \ ;$$
$$((\Delta_{xxx} V_x^1) + (\Delta_{xxx} V_x^2))/2 + (12\alpha + 280 a_{xyy} - 280 b_{xxy} - 35(\Delta_{xy} R^z))/36 = (\Delta_{xxx} V_x^{FV}) \ ;$$
$$-4a_{yy} = (\Delta_{yy} V_x^{FV}) \ ;$$
$$((\Delta_y V_y^1) + (\Delta_y V_y^2))/2 - (\Delta_{xy} R^z)/24 + \alpha/70 + (a_{yy} + b_{xy})/3 = (\Delta_y V_y^{FV}) \ ; \quad (2.20)$$
$$2b_{xxy}/3 + ((\Delta_{yy} V_y^1) + (\Delta_{yy} V_y^2))/2 = (\Delta_{yy} V_y^{FV}) \ ;$$
$$-(\Delta_x V_y^1) + (\Delta_x V_y^2) + (8a_{xyy} - (\Delta_{xx} R^z))/30 = (\Delta_{xy} V_y^{FV}) \ ;$$
$$((\Delta_{yyy} V_y^1) + (\Delta_{yyy} V_y^2))/2 + (12\alpha - 280 a_{xyy} + 280 b_{xxy} + 35(\Delta_{xy} R^z))/36 = (\Delta_{yyy} V_y^{FV}) \ ;$$
$$-4b_{xy} = (\Delta_{xy} V_y^{FV})$$

We see that there are only 5 free parameters ($a_{yy}$, $b_{xy}$, $a_{xyy}$, $b_{xxy}$ and $\alpha$) in the above 10 equations with the result that they cannot be satisfied exactly, therefore, the least squares minimization is our best option. The least squares procedure is justified because any finite volume higher order reconstruction is already quite accurate, therefore it preserves the curl constraint reasonably closely, though not exactly. The least squares procedure just does the extra little bit to keep the curl constraint exactly satisfied while bringing all the moments in eqn. (2.18a,b) as close as possible to their finite volume counterparts. Section V shows that the use of the least squares procedure does not affect the accuracy of the final fourth order scheme. Since the equations are linear in the 5 free parameters, the optimal parameters can be written explicitly as

$$\begin{aligned}
a_{xyy} = (\ & 27811175(\Delta_{xy}R^z) - 4799550(\Delta_x V_x^1) - 6945750(\Delta_{xxx}V_x^1) - 4799550(\Delta_y V_y^1) + 7357140(\Delta_{yyy}V_y^1) \\
& - 4799550(\Delta_x V_x^2) - 6945750(\Delta_{xxx}V_x^2) - 4799550(\Delta_y V_y^2) + 7357140(\Delta_{yyy}V_y^2) + 9599100(\Delta_x V_x^{FV}) \\
& + 13891500(\Delta_{xxx}V_x^{FV}) - 61378272(\Delta_{xyy}V_x^{FV}) - 54022500(\Delta_{xxy}V_y^{FV}) + 9599100(\Delta_y V_y^{FV}) \\
& - 14714280(\Delta_{yyy}V_y^{FV})\ )\ /\ 474401888\ ;
\end{aligned}$$

$$\begin{aligned}
b_{xxy} = (\ & -27811175(\Delta_{xy}R^z) - 4799550(\Delta_x V_x^1) + 7357140(\Delta_{xxx}V_x^1) - 4799550(\Delta_y V_y^1) - 6945750(\Delta_{yyy}V_y^1) \\
& -\ 4799550(\Delta_x V_x^2) + 7357140(\Delta_{xxx}V_x^2) - 4799550(\Delta_y V_y^2) - 6945750(\Delta_{yyy}V_y^2) + 9599100(\Delta_x V_x^{FV}) \\
& - 14714280(\Delta_{xxx}V_x^{FV}) - 54022500(\Delta_{xyy}V_x^{FV}) - 61378272(\Delta_{xxy}V_y^{FV}) + 9599100(\Delta_y V_y^{FV}) \\
& + 13891500(\Delta_{yyy}V_y^{FV})\ )\ /\ 474401888\ ;
\end{aligned}$$

$$\begin{aligned}
a_{xxyy} = (\ & 2(\Delta_{xxy}R^z) - 75(\Delta_{xx}V_x^1) + 60(\Delta_y V_y^1) - 75(\Delta_{xx}V_x^2) - 60(\Delta_y V_y^2) + 150(\Delta_{xx}V_x^{FV}) \\
& + 60(\Delta_{xy}V_y^{FV})\ )\ /\ 116\ ;
\end{aligned}$$

$$\begin{aligned}
b_{xxyy} = (& -2(\Delta_{xyy}R^z) + 60(\Delta_x V_x^1) - 75(\Delta_{yy}V_y^1) - 60(\Delta_x V_x^2) - 75(\Delta_{yy}V_y^2) + 60(\Delta_{xy}V_x^{FV}) \\
& + 150(\Delta_{yy}V_y^{FV})\ )\ /\ 116\ ;
\end{aligned}$$

$$\begin{aligned}
\alpha = -105\ (\ & 54(\Delta_x V_x^1) + 1295(\Delta_{xxx}V_x^1) + 54(\Delta_y V_y^1) + 1295(\Delta_{yyy}V_y^1) + 54(\Delta_x V_x^2) + 1295(\Delta_{xxx}V_x^2) \\
& + 54(\Delta_y V_y^2) + 1295(\Delta_{yyy}V_y^2) - 108(\Delta_x V_x^{FV}) - 2590(\Delta_{xxx}V_x^{FV}) - 18(\Delta_{xyy}V_x^{FV}) \\
& - 18(\Delta_{xxy}V_y^{FV}) - 108(\Delta_y V_y^{FV}) - 2590(\Delta_{yyy}V_y^{FV})\ )\ /\ 181624
\end{aligned}$$

(2.21)

Once eqn. (2.21) gives us the optimized parameters, they can be substituted in eqn. (2.17) to obtain all the coefficients in the reconstruction. This gives us a complete strategy for the vector field reconstruction in an edge-collocated, WENO-like or DG-like, curl-constraint preserving scheme.

The previous paragraph has shown us how we obtain an edge-collocated, DG-like, curl-constraint preserving reconstruction. We now show that an edge-collocated, WENO-like, curl-constraint preserving reconstruction can also be designed. Realize that eqns. (2.19) can be used to obtain zone-averaged values for the vector field. Furthermore, these zone-averaged values can be obtained without regard to the 5 free parameters that are optimized via eqn. (2.20). The values from eqn. (2.19) can be used to carry out a conventional, finite-volume, fourth order accurate WENO reconstruction and such a reconstruction provides all the moments that are required in the right-hand sides of eqn. (2.20). Therefore, one does not need to evolve any auxiliary system when constructing a fourth order, curl-constraint preserving WENO-like scheme. In Section IV we show that a fourth order WENO-like curl-constraint preserving scheme that is based on this section does indeed meet its design accuracy. This completes our description of curl-preserving reconstruction at fourth order in two dimensions.

## II.5) Combining Curl-Free Reconstruction and the Two-Dimensional Riemann Solver to Obtain a Multidimensionally Upwinded, Globally Curl-Free Scheme

When studying one-dimensional advection, it is indeed a very instructive to realize that one-dimensional upwinding from a one-dimensional Riemann solver yields a stable parabolized scheme. The transition to a higher order scheme for advection is then easy to justify with the inclusion of TVD or WENO limiters, as was shown by van Leer [55], [56], Jiang and Shu [42], Balsara and Shu [4], Balsara, Garain and Shu [19]. As the order of accuracy of the reconstruction is increased, it is easy to see that the dissipative terms will become smaller. It is very desirable to show that an analogous plan exists for PDEs with a curl involution constraint on a vector field. Such a demonstration has to be multidimensional because the curl operator is only meaningful in two or more dimensions.

Let us begin by obtaining a very important insight from divergence-preserving vector fields. In their study of the induction equation for globally divergence-free MHD, Balsara and Käppeli [20] were able to show that multidimensional divergence-free reconstruction of the magnetic field at first order, coupled with a two-dimensional Riemann solver, also results in a stable scheme with the correct, multidimensionally parabolized, dissipation; see Section 4 of the above-mentioned paper. Because of the presence of these parabolic terms, the multidimensional Riemann solver *always plays a stabilizing* role in the induction equation. This showed us that the transition to higher order, by applying WENO-type or DG-type limiting, would indeed succeed for the induction equation. Notice that Balsara and Käppeli [20] have made a non-trivial demonstration, because the components of a divergence-free vector field are collocated at the faces of the mesh, which is quite different from the zone-centered collocation that is used in DG schemes for conservation laws.

The previous paragraph has shown us that divergence-free evolution of vector fields has been developed to a point where it has a solid footing. It is very desirable to show that vector fields that evolve in a curl-free fashion have a similar assurance. This will again be a non-trivial demonstration because the components of a curl-free vector field are indeed collocated at the edges of the mesh. This is indeed different from the zone-centered collocation that is used for conservation laws as well as the face-centered collocation that is used for mimetic schemes that support globally divergence-free evolution. (Curl constraint-preserving evolution is then just a matter of adding source terms to the curl-free evolution equations; so we do not consider that in this Sub-section.) Let us consider the simplest two-dimensional equations that give us curl-free evolution. They can be written as

$$\frac{\partial J^x}{\partial t} + \frac{\partial \left( v^x J^x + v^y J^y \right)}{\partial x} + v^y \left( \frac{\partial J^x}{\partial y} - \frac{\partial J^y}{\partial x} \right) = 0 \; ;$$

$$\frac{\partial J^y}{\partial t} + \frac{\partial \left( v^x J^x + v^y J^y \right)}{\partial y} + v^x \left( \frac{\partial J^y}{\partial x} - \frac{\partial J^x}{\partial y} \right) = 0$$

(2.22)

Here, for the purposes of our theoretical study in this Sub-section, we will take the velocity components $\left( v^x, v^y \right)$ to be constant. Let us view the update equations shown above in two space and one time direction. We also consider a uniform Cartesian mesh in the two spatial dimensions with zones of size $\Delta x$ and $\Delta y$ in the x- and y-directions. Let the timestep be of size $\Delta t$. The first equation in eqn. (2.22) can be integrated along the x- and t-directions of the three-dimensional space-time mesh (2 space + 1 time dimensions). Since $J_x$ is collocated at the x-edges of the mesh, the x-directional integration should, of course, coincide with the x-edges of the Cartesian mesh. As long as the curl is exactly zero, the $v^y \left( \partial J^x / \partial y - \partial J^y / \partial x \right)$ term in that equation will be exactly zero; and the $\partial \left( v^x J^x + v^y J^y \right) / \partial x$ term is just a gradient applied in the x-direction. Similarly, the second equation in eqn. (2.22) can be integrated along the y- and t-directions of a three-dimensional space-time mesh. Since $J_y$ is collocated at the y-edges of the mesh, the y-directional integration should, of course, coincide with the y-edges of the two-dimensional Cartesian mesh. As long as the curl is exactly zero, the $v^x \left( \partial J^y / \partial x - \partial J^x / \partial y \right)$ term in that equation will be exactly zero; and the $\partial \left( v^x J^x + v^y J^y \right) / \partial y$ term is just a gradient applied in the y-direction. It is then easy to see that eqn. (2.22) ensures that if the *same, properly upwinded, vertex-centered*, values of the potential, defined in two dimensions by $\phi \equiv \left( v^x J^x + v^y J^y \right)$, are used at the corners of the mesh in order to update the x-edge centered $J^x$ and the y-edge centered $J^y$ then the update will be globally curl-free. (This is very similar to globally divergence-free update of the induction equation in MHD which requires the *same, properly upwinded, edge-centered*, electric fields to be used for the update of the face-centered components of the magnetic induction vector.)

**Place Fig. 2 here**

Let us now consider constant velocity components with $v^x > 0$ and $v^y > 0$, just to keep our initial discussion simple. Eqn. (2.22) becomes

$$\frac{\partial J^x}{\partial t} + v^x \frac{\partial J^x}{\partial x} + v^y \frac{\partial J^x}{\partial y} = 0 \; ;$$
$$\frac{\partial J^y}{\partial t} + v^x \frac{\partial J^y}{\partial x} + v^y \frac{\partial J^y}{\partial y} = 0 \quad (2.23)$$

We see, therefore, that "properly upwinded" in this context means that each zone in Fig. 2 will contribute its $J^x$ and $J^y$ values to its top right corner in Fig. 2. We require that the first order accurate curl-free reconstruction from eqn. (2.1) should be used in the zones $(i, j)$, $(i-1, j)$, $(i-1, j-1)$ and $(i, j-1)$ in Fig. 2. Fig. 2 also catalogues the values of the potential $\phi \equiv \left( v^x J^x + v^y J^y \right)$ at all the vertices of zone $(i, j)$. We are interested in the time update of the components $J^{x;n}_{i, j-1/2}$ and $J^{y;n}_{i-1/2, j}$ at the lower x-edge and lower y-edge of the zone $(i, j)$ respectively. The subscripts in $J^{x;n}_{i, j-1/2}$ and $J^{y;n}_{i-1/2, j}$ indicate spatial collocation points on the mesh; the superscript of "$n$" denotes the $n^{th}$ timestep. The components of the vector field **J** are also indicated by "$x$" or

"y" in the superscripts. We can write the discretized, first order in space and time, update equations as

$$\frac{J^{x;n+1}_{i,j-1/2} - J^{x;n}_{i,j-1/2}}{\Delta t} + \frac{\phi^n_{i+1/2,j-1/2} - \phi^n_{i-1/2,j-1/2}}{\Delta x} = 0 \ ;$$
$$\frac{J^{y;n+1}_{i-1/2,j} - J^{y;n}_{i-1/2,j}}{\Delta t} + \frac{\phi^n_{i-1/2,j+1/2} - \phi^n_{i-1/2,j-1/2}}{\Delta x} = 0$$
(2.24)

Recall that we assume a uniform Cartesian mesh with zone sizes $\Delta x$ and $\Delta y$ in the x- and y-directions and a timestep of size $\Delta t$. By substituting the upwinded potentials from Fig. 2, we get

$$\frac{J^{x;n+1}_{i,j-1/2} - J^{x;n}_{i,j-1/2}}{\Delta t} + \frac{\left[\left(v^x J^{x;n}_{i,j-1/2} + v^y J^{y;n}_{i+1/2,j-1}\right) - \left(v^x J^{x;n}_{i-1,j-1/2} + v^y J^{y;n}_{i-1/2,j-1}\right)\right]}{\Delta x} = 0 \ ;$$
$$\frac{J^{y;n+1}_{i-1/2,j} - J^{y;n}_{i-1/2,j}}{\Delta t} + \frac{\left[\left(v^x J^{x;n}_{i-1,j+1/2} + v^y J^{y;n}_{i-1/2,j}\right) - \left(v^x J^{x;n}_{i-1,j-1/2} + v^y J^{y;n}_{i-1/2,j-1}\right)\right]}{\Delta y} = 0$$
(2.25)

We now write the above equations in a format that allows us to see the velocity-dependence more clearly as follows

$$\frac{J^{x;n+1}_{i,j-1/2} - J^{x;n}_{i,j-1/2}}{\Delta t} + \frac{v^x}{\Delta x}\left(J^{x;n}_{i,j-1/2} - J^{x;n}_{i-1,j-1/2}\right) + \frac{v^y}{\Delta x}\left(J^{y;n}_{i+1/2,j-1} - J^{y;n}_{i-1/2,j-1}\right) = 0 \ ;$$
$$\frac{J^{y;n+1}_{i-1/2,j} - J^{y;n}_{i-1/2,j}}{\Delta t} + \frac{v^x}{\Delta y}\left(J^{x;n}_{i-1,j+1/2} - J^{x;n}_{i-1,j-1/2}\right) + \frac{v^y}{\Delta y}\left(J^{y;n}_{i-1/2,j} - J^{y;n}_{i-1/2,j-1}\right) = 0$$
(2.26)

The above equations still do not look like discretized versions of eqn. (2.23). As with the divergence-free evolution of vector fields, the concordance will only be established if the discrete circulations in the zones $(i, j-1)$ and $(i-1, j)$ are used in the first and second equations of eqn. (2.26). If we utilize the fact that the discrete circulations in those two zones are indeed zero, we can make the transcriptions

$$\frac{1}{\Delta x}\left(J^{y;n}_{i+1/2,j-1} - J^{y;n}_{i-1/2,j-1}\right) \to \frac{1}{\Delta y}\left(J^{x;n}_{i,j-1/2} - J^{x;n}_{i,j-3/2}\right) \ ;$$
$$\frac{1}{\Delta y}\left(J^{x;n}_{i-1,j+1/2} - J^{x;n}_{i-1,j-1/2}\right) \to \frac{1}{\Delta x}\left(J^{y;n}_{i-1/2,j} - J^{y;n}_{i-3/2,j}\right)$$
(2.27)

Inserting the transcriptions from eqn. (2.27) into eqn. (2.26), we now get

$$\frac{J_{i,j-1/2}^{x;n+1} - J_{i,j-1/2}^{x;n}}{\Delta t} + \frac{v^x}{\Delta x}\left(J_{i,j-1/2}^{x;n} - J_{i-1,j-1/2}^{x;n}\right) + \frac{v^y}{\Delta y}\left(J_{i,j-1/2}^{x;n} - J_{i,j-3/2}^{x;n}\right) = 0 \; ;$$

$$\frac{J_{i-1/2,j}^{y;n+1} - J_{i-1/2,j}^{y;n}}{\Delta t} + \frac{v^x}{\Delta x}\left(J_{i-1/2,j}^{y;n} - J_{i-3/2,j}^{y;n}\right) + \frac{v^y}{\Delta y}\left(J_{i-1/2,j}^{y;n} - J_{i-1/2,j-1}^{y;n}\right) = 0$$

(2.28)

It is now easy to see that the first equation in eqn. (2.28) is just a first order upwinded approximation for the time update of $J_{i,j-1/2}^{x;n}$. Please compare that equation to the first equation in eqn. (2.23). Also please envision it as an upwind scheme for $J_{i,j-1/2}^{x;n}$ applied to the red control volume in Fig. 2. It is also easy to see that the second equation in eqn. (2.28) is just a first order upwinded approximation for the time update of $J_{i-1/2,j}^{y;n}$. Please compare that equation to the second equation in eqn. (2.23). Also please envision it as an upwind scheme for $J_{i-1/2,j}^{y;n}$ applied to the blue control volume in Fig. 2. The curl-free reconstruction which couples all components of the vector **J**, along with the use of a unique upwinded potential $\phi$ at the vertices of the mesh, makes the time-evolution curl-free. It is very important to understand the role of the curl-free reconstruction which couples all the edge-collocated components of the vector **J**. It is similarly very important to understand the role of the uniquely defined, multidimensionally upwinded, potential $\phi$ at each of the vertices of the mesh which couples all the update equations for all the components of **J**. The two innovations work together to yield the globally curl-preserving scheme. We therefore see with the help of Fig. 2 that the curl-free reconstruction, along with multidimensional upwinding applied to the vertices of the mesh, gives us a globally curl-free update strategy for our model equations, i.e. eqn. (2.22). Realize too that the multidimensional Riemann solver is indeed designed to automate that multidimensional upwinding in the general case of a system of PDEs. Therefore, we realize that the curl-free reconstruction (and the edge-centered collocation of vector components that it entails), along with the application of a multidimensional Riemann solver, indeed gives us a stable, globally curl-free, mimetic update strategy for our curl-free model equations.

The discussion in the previous paragraph was restricted to first order accuracy. To keep the discussion extremely accessible, we also drew on our notional understanding of multidimensional upwinding, and we kept all the velocities positive. In general, we will want to use a full multidimensional Riemann solver which can accommodate to waves propagating in a general hyperbolic system in any direction. Eqns. (12), (13) and (14) of Balsara [13] show how such a multidimensional Riemann solver can be designed for systems in the general case. The multidimensional Riemann solver from Balsara [13] works on Cartesian meshes. It was extended to unstructured meshes in Balsara and Dumbser [15]. For our model equations in eqn. (2.22), we do not have the support of a general underlying PDE system, however, we can design a multidimensional Riemann solver (using eqns. (12), (13) and (14) of Balsara [13]) in the locally Lax-Friedrichs (LLF) approximation. This enables us to write our potentials at the vertices of Fig. 2 as a centered part as well as a multidimensionally dissipative part. We will then write discrete evolution equations and show that we get a centered update along with a parabolic contribution at first order. This makes it easier for us to understand that as the spatial order of accuracy of our curl-free reconstruction is improved, and as the accuracy of our time-update is improved, the dissipative parts will become progressively smaller. In other words, the same exercise that was presented in Section 4 of Balsara and Käppeli [20] for mimetic, globally divergence-free schemes is now replicated in the ensuing paragraphs for mimetic, globally curl-free schemes.

We still consider the case where the velocity components $(v^x, v^y)$ are constant but now they can have any sign. We write the potential at the vertex $(i-1/2, j-1/2)$ as a centered term plus a diffusive term. We have

$$\begin{aligned}\phi^n_{i-1/2, j-1/2} &= \frac{1}{2} v^x \left( J^{x;n}_{i, j-1/2} + J^{x;n}_{i-1, j-1/2} \right) - \frac{1}{2} |v^x| \left( J^{x;n}_{i, j-1/2} - J^{x;n}_{i-1, j-1/2} \right) \\ &+ \frac{1}{2} v^y \left( J^{y;n}_{i-1/2, j} + J^{y;n}_{i-1/2, j-1} \right) - \frac{1}{2} |v^y| \left( J^{y;n}_{i-1/2, j} - J^{y;n}_{i-1/2, j-1} \right)\end{aligned} \quad (2.29)$$

Analogous expressions can be written for $\phi^n_{i+1/2, j-1/2}$ and $\phi^n_{i-1/2, j+1/2}$ by appropriate shifting of the indices. Making transcriptions that are analogous to the ones in eqn. (2.27) we obtain

$$\frac{J^{x;n+1}_{i,j-1/2} - J^{x;n}_{i,j-1/2}}{\Delta t} + \frac{v^x}{2\Delta x}\left(J^{x;n}_{i+1,j-1/2} - J^{x;n}_{i-1,j-1/2}\right) + \frac{v^y}{2\Delta y}\left(J^{x;n}_{i,j+1/2} - J^{x;n}_{i,j-3/2}\right) =$$

$$\frac{|v^x|}{2\Delta x}\left(J^{x;n}_{i+1,j-1/2} - 2J^{x;n}_{i,j-1/2} + J^{x;n}_{i-1,j-1/2}\right) + \frac{|v^y|}{2\Delta y}\left(J^{x;n}_{i,j+1/2} - 2J^{x;n}_{i,j-1/2} + J^{x;n}_{i,j-3/2}\right) \ ;$$

$$\frac{J^{y;n+1}_{i-1/2,j} - J^{y;n}_{i-1/2,j}}{\Delta t} + \frac{v^x}{2\Delta x}\left(J^{y;n}_{i+1/2,j} - J^{y;n}_{i-3/2,j}\right) + \frac{v^y}{2\Delta y}\left(J^{y;n}_{i-1/2,j+1} - J^{y;n}_{i-1/2,j-1}\right) =$$

$$\frac{|v^x|}{2\Delta x}\left(J^{y;n}_{i+1/2,j} - 2J^{y;n}_{i-1/2,j} + J^{y;n}_{i-3/2,j}\right) + \frac{|v^y|}{2\Delta y}\left(J^{y;n}_{i-1/2,j+1} - 2J^{y;n}_{i-1/2,j} + J^{y;n}_{i-1/2,j-1}\right)$$

(2.30)

The above equations clearly show us the centered terms and the dissipation terms in the first order scheme. We see that the dissipation terms have the correct scaling to perfectly stabilize the scheme. The use of higher order curl-free reconstruction and higher order time stepping, in conjunction with the multidimensional Riemann solver, will reduce the dissipation. We now have our assurance that higher order, mimetic, globally curl-free schemes will be stable.

### III) Curl-Preserving Reconstruction on a Three-dimensional Cartesian Mesh

The prior exercise at reconstructing curl-free or curl-preserving vector fields in two dimensions has left us with two very valuable insights. We list them below:-

**1)** We see that one should start at the lowest order and systematically build up to higher orders. This is because the polynomial terms at each order give us insight into which modes are needed at the next higher order.

**2)** We also obtained the insight that some additional polynomial contributions will be needed to ensure curl-free vector fields at all locations within a zone. Because of the nature of the curl operator, and our need to only use the lowest order additional polynomials so as to retain a modicum of stability, we found that one component of the vector field usually takes on contributions that help cancel extra terms that arise in another component of the vector field.

Armed with these insights, we now extend our studies to three-dimensional Cartesian meshes.

Sub-section III.1 deals with the first order reconstruction on 3D Cartesian meshes. Sub-section III.2 extends this to second order. Sub-section III.3 presents the third order case. Sub-section III.4 shows that the curl-free reconstruction, when combined with a three-dimensional Riemann solver, produces a properly upwinded numerical scheme. Sub-section III.5 gives a step-by-step implementation strategy.

**III.1) Curl-Preserving Reconstruction on a Three-dimensional Cartesian Mesh at First Order**

Fig. 3 shows the collocation of vector components along the edges of the control volume. Within each of the two x-faces, the two y-faces and the two z-faces, the discrete circulation (evaluated over those faces) is exactly zero if we are dealing with a curl-free PDE. If the PDE is only curl-preserving, the discrete circulation is easily evaluated at each face, as we soon show. The mean value of the vector field's parallel component and its linear variation are shown along each edge, in anticipation of a second order accurate reconstruction scheme. The reconstruction problem consists of obtaining a polynomial that is globally curl-free/curl-preserving within this control volume. The discrete circulation within each of the six faces of the zone shown in Fig. 3 gives us the six conditions for the mean values. For simplicity, let us consider curl-free evolution in the next three equations. At the top and bottom x-faces we have

$$V_y^3 + V_z^4 - V_y^4 - V_z^2 = 0 \quad ; \quad V_y^1 + V_z^3 - V_y^2 - V_z^1 = 0 \tag{3.1}$$

At the top and bottom y-faces we have

$$V_z^3 + V_x^4 - V_z^4 - V_x^2 = 0 \quad ; \quad V_z^1 + V_x^3 - V_z^2 - V_x^1 = 0 \tag{3.2}$$

At the top and bottom z-faces we have

$$V_x^3 + V_y^4 - V_x^4 - V_y^2 = 0 \quad ; \quad V_x^1 + V_y^3 - V_x^2 - V_y^1 = 0 \tag{3.3}$$

In general, the second and higher moments in Fig. 3 are only obtained with some level of approximation. So we cannot guarantee that an analogous set of equations hold for the slopes.

**Place Fig. 3 here**

It is very important to make one further observation about a curl-free vector field. The six discrete curl conditions in eqns. (3.1), (3.2) and (3.3) at the six faces of the mesh in Fig. 3 are not independent. The circulation condition at any one face can always be written in terms of the circulation conditions at the other five faces. Let us consider the z-faces and say that we know the condition $V_y^3 - V_x^2 - V_y^1 + V_x^1 = 0$ at the bottom face. We wish to show that the condition at the top face can be obtained in terms of the condition at the bottom face and the use of eqns. (3.1) and (3.2). To that end, realize that $V_y^3$ appears only in the first discrete circulation condition in eqn. (3.1) and can be written as $V_y^3 = -V_z^4 + V_y^4 + V_z^2$. Similarly, $V_y^1$ appears only in the second discrete circulation condition in eqn. (3.1) and can be written as $-V_y^1 = +V_z^3 - V_y^2 - V_z^1$. Likewise, $V_x^2$ only appears in the first discrete circulation condition in eqn. (3.2) and can be written as $-V_x^2 = V_z^4 - V_x^4 - V_z^3$. Furthermore, $V_x^1$ only appears in the second discrete circulation condition in eqn. (3.2) and can be written as $V_x^1 = -V_z^2 + V_x^3 + V_z^1$. Adding the above four equations immediately retrieves the first discrete circulation condition in eqn. (3.3). We have started with the second discrete circulation condition in eqn. (3.3) and shown that it immediately proves the first circulation condition in eqn. (3.3) if the four other circulation conditions in eqns. (3.1) and (3.2) can be assumed true. Therefore, on a mesh with six faces, only five of them are truly, mutually independent.

Consider the following reconstruction on the unit cube spanning $[-1/2, 1/2]^3$. We write the three components of the first order accurate curl-preserving vector field as

$$V^x(x,y,z) = V_x^1\left(y - \frac{1}{2}\right)\left(z - \frac{1}{2}\right) + V_x^2\left(y + \frac{1}{2}\right)\left(\frac{1}{2} - z\right)$$
$$+ V_x^3\left(\frac{1}{2} - y\right)\left(z + \frac{1}{2}\right) + V_x^4\left(y + \frac{1}{2}\right)\left(z + \frac{1}{2}\right)$$

(3.4)

and

$$V^y(x,y,z) = V_y^1\left(z - \frac{1}{2}\right)\left(x - \frac{1}{2}\right) + V_y^2\left(z + \frac{1}{2}\right)\left(\frac{1}{2} - x\right)$$
$$+ V_y^3\left(\frac{1}{2} - z\right)\left(x + \frac{1}{2}\right) + V_y^4\left(z + \frac{1}{2}\right)\left(x + \frac{1}{2}\right)$$

(3.5)

and

$$V^z(x,y,z) = V_z^1\left(x-\frac{1}{2}\right)\left(y-\frac{1}{2}\right) + V_z^2\left(x+\frac{1}{2}\right)\left(\frac{1}{2}-y\right)$$
$$+ V_z^3\left(\frac{1}{2}-x\right)\left(y+\frac{1}{2}\right) + V_z^4\left(x+\frac{1}{2}\right)\left(y+\frac{1}{2}\right)$$
(3.6)

In the above three equations, the polynomial space is so chosen that the edge values of the vector field components are exactly retrieved at the edges of the zone shown in Fig. 3. As a result, the vector field will be globally constraint-preserving. We evaluate the curl of the above vector field in the next paragraph.

Evaluating $(\nabla \times \mathbf{V})_x$ we get

$$(\nabla \times \mathbf{V})_x = \frac{\partial V^z}{\partial y} - \frac{\partial V^y}{\partial z} = \frac{1}{2}\left(V_y^1 - V_y^2 + V_y^3 - V_y^4 - V_z^1 - V_z^2 + V_z^3 + V_z^4\right)$$
$$+ \left(-V_y^1 + V_y^2 + V_y^3 - V_y^4 + V_z^1 - V_z^2 - V_z^3 + V_z^4\right)x$$
(3.7)

Set $x = \pm 1/2$ in the above equation to see that it retrieves the discrete circulation equations in eqn. (3.1) at the top and bottom x-faces respectively. If the vector field is circulation-free then this guarantees that $(\nabla \times \mathbf{V})_x$ is strictly zero in the above equation. Evaluating $(\nabla \times \mathbf{V})_y$ we get

$$(\nabla \times \mathbf{V})_y = \frac{\partial V^x}{\partial z} - \frac{\partial V^z}{\partial x} = \frac{1}{2}\left(V_z^1 - V_z^2 + V_z^3 - V_z^4 - V_x^1 - V_x^2 + V_x^3 + V_x^4\right)$$
$$+ \left(-V_z^1 + V_z^2 + V_z^3 - V_z^4 + V_x^1 - V_x^2 - V_x^3 + V_x^4\right)y$$
(3.8)

Set $y = \pm 1/2$ in the above equation to see that it retrieves the discrete circulation equations in eqn. (3.2) at the top and bottom y-faces respectively. If the vector field is circulation-free then this guarantees that $(\nabla \times \mathbf{V})_y$ is strictly zero in the above equation. Evaluating $(\nabla \times \mathbf{V})_z$ we get

$$(\nabla \times \mathbf{V})_z = \frac{\partial V^y}{\partial x} - \frac{\partial V^x}{\partial y} = \frac{1}{2}\left(V_x^1 - V_x^2 + V_x^3 - V_x^4 - V_y^1 - V_y^2 + V_y^3 + V_y^4\right)$$
$$+ \left(-V_x^1 + V_x^2 + V_x^3 - V_x^4 + V_y^1 - V_y^2 - V_y^3 + V_y^4\right)z$$
(3.9)

Set $z = \pm 1/2$ in the above equation to see that it retrieves the discrete circulation equations in eqn. (3.3) at the top and bottom z-faces respectively. If the vector field is circulation-free then this guarantees that $(\nabla \times \mathbf{V})_z$ is strictly zero in the above equation.

From eqns. (3.7), (3.8) and (3.9) we also see that $\nabla \cdot (\nabla \times \mathbf{V}) = \partial_x (\nabla \times \mathbf{V})_x + \partial_y (\nabla \times \mathbf{V})_y + \partial_z (\nabla \times \mathbf{V})_z = 0$. In other words, the discrete divergence of the discrete curl is also exactly zero. This also tells us that a good zone-centered approximation of $(\nabla \times \mathbf{V})_x$ is given by the first term in eqn. (3.7). Likewise, a good zone-centered approximation of $(\nabla \times \mathbf{V})_y$ is given by the first term in eqn. (3.8). Similarly, a good zone-centered approximation of $(\nabla \times \mathbf{V})_z$ is given by the first term in eqn. (3.9). We will see in the next section that while these approximations are available, they are indeed sub-optimal because better approximations of the circulation are directly available in the faces of the mesh.

Now notice that eqns. (3.4), (3.5) and (3.6) are only first order accurate. This is because eqn. (3.4) lacks any linear variation in the x-direction; eqn. (3.5) lacks any linear variation in the y-direction and eqn. (3.6) lacks any linear variation in the z-direction. However, for the first order accurate case, the equations are exactly curl-free or curl-preserving. Besides, the first order curl-free reconstruction reflects the six discrete circulations evaluated over the six faces of the control volume using the edges of the same control volume. (This mirrors the known fact that at first order, the discrete divergence-preserving reconstruction reflects the one discrete divergence evaluated over the control volume using the faces of the same control volume.) Notice too that while there is only one divergence condition evaluated over a 3D control volume in a divergence-preserving scheme, there are six curl conditions evaluated over a 3D control volume in a curl-preserving scheme. This makes curl-free reconstruction more complicated than divergence-free reconstruction, especially in three dimensions.

Eqns. (3.7), (3.8) and (3.9) give us yet another insight if we compare them to eqn. (2.2). Notice that eqn. (2.2) is a single equation for the discrete circulation. However, because of the linear variation in the x-direction, eqn. (3.7) is an expression of the discrete circulation at the top and bottom x-faces. Similarly, because of the linear variation in the y-direction, eqn. (3.8) is an expression of the discrete circulation at the top and bottom y-faces. Likewise, because of the linear variation in the z-direction, eqn. (3.9) is an expression of the discrete circulation at the top and bottom z-faces. We, therefore, understand that three-dimensional curl-preserving reconstruction can be quite different from two-dimensional curl-preserving reconstruction. In three-dimensional curl-preserving reconstruction, even when we make a higher order reconstruction of the curl/circulation vector in the zone of interest, the two modes that are present in each of eqns. (3.7), (3.8) and (3.9) must be retained.

**III.2) Curl-Preserving Reconstruction on a Three-dimensional Cartesian Mesh at Second Order**

Notice that eqn. (3.4) already has a constant part and $y$, $z$ and $yz$ variation. To attain full second order accuracy, we need to add a linear x-directional variation to $V^x(x,y,z)$. This would be added to the x-edges of the zone shown in Fig. 3. The inclusion of such an x-variation will also trigger $xy$, $xz$ and $xyz$ terms in $V^x(x,y,z)$. Similarly, notice that eqn. (3.5) already has a constant part and $x$, $z$ and $xz$ variation. To attain full second order accuracy, we need to add a linear y-directional variation to $V^y(x,y,z)$. This would be added to the y-edges of the zone shown in Fig. 3. The inclusion of such a y-variation will also trigger $xy$, $yz$ and $xyz$ terms in $V^y(x,y,z)$. Likewise, notice that eqn. (3.6) already has a constant part and $x$, $y$ and $xy$ variation. To attain full second order accuracy, we need to add a linear z-directional variation to $V^z(x,y,z)$. This would be added to the z-edges of the zone shown in Fig. 3. The inclusion of such a z-variation will also trigger $xz$, $yz$ and $xyz$ terms in $V^z(x,y,z)$. The inclusion of all these terms also causes the curl operator to acquire additional moments and to ensure curl-free reconstruction (or to make sure that the curl-preserving reconstruction has the appropriate moments) we need to add some complementing terms.

Notice that in the ensuing three equations, the polynomials are chosen with such factors that they do not affect the vector components, or their linear variation, in the edges. In other words, if the ensuing polynomial for $V^x(x,y,z)$ is evaluated at any x-edge of Fig. 3, we indeed retrieve only the constant and linear variations in that x-edge, and nothing but that variation. Similarly, if $V^y(x,y,z)$ is evaluated at any y-edge of Fig. 3, we indeed retrieve only the constant and linear variations in that y-edge, and nothing but that variation. Similarly, if $V^z(x,y,z)$ is evaluated at any z-edge of Fig. 3, we indeed retrieve only the constant and linear variations in that z-edge, and nothing but that variation. We write the three components of the second order accurate curl-preserving vector field as

$$V^x(x,y,z) = \left[V_x^1 + (\Delta_x V_x^1)x\right]\left(y - \frac{1}{2}\right)\left(z - \frac{1}{2}\right) + \left[V_x^2 + (\Delta_x V_x^2)x\right]\left(y + \frac{1}{2}\right)\left(\frac{1}{2} - z\right)$$
$$+ \left[V_x^3 + (\Delta_x V_x^3)x\right]\left(\frac{1}{2} - y\right)\left(z + \frac{1}{2}\right) + \left[V_x^4 + (\Delta_x V_x^4)x\right]\left(y + \frac{1}{2}\right)\left(z + \frac{1}{2}\right) \quad (3.10)$$
$$+ a_{yy}(1 - 4y^2) + a_{zz}(1 - 4z^2) + a_{yyz} z(1 - 4y^2) + a_{yzz} y(1 - 4z^2)$$

and

$$V^y(x,y,z) = \left[V_y^1 + (\Delta_y V_y^1)y\right]\left(z - \frac{1}{2}\right)\left(x - \frac{1}{2}\right) + \left[V_y^2 + (\Delta_y V_y^2)y\right]\left(z + \frac{1}{2}\right)\left(\frac{1}{2} - x\right)$$
$$+ \left[V_y^3 + (\Delta_y V_y^3)y\right]\left(\frac{1}{2} - z\right)\left(x + \frac{1}{2}\right) + \left[V_y^4 + (\Delta_y V_y^4)y\right]\left(z + \frac{1}{2}\right)\left(x + \frac{1}{2}\right) \quad (3.11)$$
$$+ b_{xx}(1 - 4x^2) + b_{zz}(1 - 4z^2) + b_{xxz} z(1 - 4x^2) + b_{xzz} x(1 - 4z^2)$$

and

$$V^z(x,y,z) = \left[V_z^1 + (\Delta_z V_z^1)z\right]\left(x - \frac{1}{2}\right)\left(y - \frac{1}{2}\right) + \left[V_z^2 + (\Delta_z V_z^2)z\right]\left(x + \frac{1}{2}\right)\left(\frac{1}{2} - y\right)$$
$$+ \left[V_z^3 + (\Delta_z V_z^3)z\right]\left(\frac{1}{2} - x\right)\left(y + \frac{1}{2}\right) + \left[V_z^4 + (\Delta_z V_z^4)z\right]\left(x + \frac{1}{2}\right)\left(y + \frac{1}{2}\right) \quad (3.12)$$
$$+ c_{xx}(1 - 4x^2) + c_{yy}(1 - 4y^2) + c_{xxy} y(1 - 4x^2) + c_{xyy} x(1 - 4y^2)$$

Now let us study the curl-preserving constraints that are put on the above vector field. The object of this study is to derive insights that help us to fix up the 12 coefficients – $a_{yy}, a_{zz}, a_{yyz}, a_{yzz}$, $b_{xx}, b_{zz}, b_{xxz}, b_{xzz}, c_{xx}, c_{yy}, c_{xxy}, c_{xyy}$ – in eqns. (3.10), (3.11) and (3.12) above.

From the Introduction, we have seen that every curl-preserving, edge-centered vector field has to be such that its curl matches a divergence-free vector field. This divergence-free vector field has components that are face-centered. Therefore, at the top and bottom x-faces of Fig. 3 we have the mean circulations given by $\left(V_y^3 + V_z^4 - V_y^4 - V_z^2\right)$ and $\left(V_y^1 + V_z^3 - V_y^2 - V_z^1\right)$ respectively. Using two-dimensional, second order, finite-volume WENO reconstruction that is restricted to the x-faces of the mesh, we can write the facial variation of the x-component of the circulation in the top and bottom x-faces as

$$R^{x+}(y,z) = \left(V_y^3 + V_z^4 - V_y^4 - V_z^2\right) + R_y^{x+} y + R_z^{x+} z \quad ;$$
$$R^{x-}(y,z) = \left(V_y^1 + V_z^3 - V_y^2 - V_z^1\right) + R_y^{x-} y + R_z^{x-} z$$

(3.13a)

Similarly, at the top and bottom y-faces of Fig. 3 we have the mean circulations given by $\left(V_z^3 + V_x^4 - V_z^4 - V_x^2\right)$ and $\left(V_z^1 + V_x^3 - V_z^2 - V_x^1\right)$ respectively. Using two-dimensional, second order, finite-volume WENO reconstruction that is restricted to the y-faces of the mesh, we can write the facial variation of the y-component of the circulation in the top and bottom y-faces as

$$R^{y+}(x,z) = \left(V_z^3 + V_x^4 - V_z^4 - V_x^2\right) + R_x^{y+} x + R_z^{y+} z \quad ;$$
$$R^{y-}(x,z) = \left(V_z^1 + V_x^3 - V_z^2 - V_x^1\right) + R_x^{y-} x + R_z^{y-} z$$

(3.13b)

Likewise, at the top and bottom z-faces of Fig. 3 we have the mean circulations given by $\left(V_x^3 + V_y^4 - V_x^4 - V_y^2\right)$ and $\left(V_x^1 + V_y^3 - V_x^2 - V_y^1\right)$ respectively. Using two-dimensional, second order, finite-volume WENO reconstruction that is restricted to the z-faces of the mesh, we can write the facial variation of the z-component of the circulation in the top and bottom z-faces as

$$R^{z+}(x,y) = \left(V_x^3 + V_y^4 - V_x^4 - V_y^2\right) + R_x^{z+} x + R_y^{z+} y \quad ;$$
$$R^{z-}(x,y) = \left(V_x^1 + V_y^3 - V_x^2 - V_y^1\right) + R_x^{z-} x + R_y^{z-} y$$

(3.13c)

Notice that WENO has given us the facial variations of the x-, y- and z-components of the circulations in the x-, y- and z-faces of the mesh. But we still do not have expressions for the circulation vector field over the volume of the zone that we are considering.

We can obtain expressions for the x-, y- and z-components of the circulation vector field by using the divergence-free reconstruction strategy catalogued in Balsara *et al.* [22] (see Section 3 of that paper). The x-component of the volumetrically reconstructed circulation can be written as

$$R^x(x,y,z) = \alpha_0 + \alpha_x x + \alpha_y y + \alpha_z z + \alpha_{xx}\left(x^2 - 1/12\right) + \alpha_{xy} xy + \alpha_{xz} xz \tag{3.14a}$$

The y-component of the volumetrically reconstructed circulation can be written as

$$R^y(x,y,z) = \beta_0 + \beta_x x + \beta_y y + \beta_z z + \beta_{xy} xy + \beta_{yy}\left(y^2 - 1/12\right) + \beta_{yz} yz \tag{3.14b}$$

The z-component of the volumetrically reconstructed circulation can be written as

$$R^z(x,y,z) = \gamma_0 + \gamma_x x + \gamma_y y + \gamma_z z + \gamma_{xz} xz + \gamma_{yz} yz + \gamma_{zz}\left(z^2 - 1/12\right) \tag{3.14c}$$

Using the divergence-free constraint, Section 3 of Balsara *et al.* [22] gives us an explicit procedure for matching the 21 coefficients in eqns. (3.14a), (3.14b) and (3.14c) to the coefficients in eqns. (3.13a), (3.13b) and (3.13c). We do not repeat the procedure here because all the formulae have already been explicitly catalogued in the literature. Realize that we now have divergence-free expressions for the circulation vector field over the volume of the zone that we are considering.

With the divergence-free expressions for the circulation vector field in hand, we can now match the curl of eqns. (3.10), (3.11) and (3.12) to eqns. (3.14a), (3.14b) and (3.14c). Operationally, we write the three equations $\partial_y V^z(x,y,z) - \partial_z V^y(x,y,z) = R^x(x,y,z)$, $\partial_z V^x(x,y,z) - \partial_x V^z(x,y,z) = R^y(x,y,z)$ and $\partial_x V^y(x,y,z) - \partial_y V^x(x,y,z) = R^z(x,y,z)$. We then match all non-zero moments across the three previous equations. This gives us a sparse linear system which we can solve for the 12 free coefficients in eqns. (3.10), (3.11) and (3.12). With a little ingenuity, the coefficients in eqn. (3.10) are found to be

$$\begin{aligned}
a_{yy} &= \frac{1}{16}\left[R_y^{z-} + R_y^{z+} + \left(\Delta_y V_y^1\right) + \left(\Delta_y V_y^2\right) - \left(\Delta_y V_y^3\right) - \left(\Delta_y V_y^4\right)\right] ; \\
a_{zz} &= \frac{1}{16}\left[-R_z^{y-} - R_z^{y+} + \left(\Delta_z V_z^1\right) - \left(\Delta_z V_z^2\right) + \left(\Delta_z V_z^3\right) - \left(\Delta_z V_z^4\right)\right] ; \\
a_{yyz} &= \frac{1}{8}\left[-R_y^{z-} + R_y^{z+} - \left(\Delta_y V_y^1\right) + \left(\Delta_y V_y^2\right) + \left(\Delta_y V_y^3\right) - \left(\Delta_y V_y^4\right)\right] ; \\
a_{yzz} &= \frac{1}{8}\left[R_z^{y-} - R_z^{y+} - \left(\Delta_z V_z^1\right) + \left(\Delta_z V_z^2\right) + \left(\Delta_z V_z^3\right) - \left(\Delta_z V_z^4\right)\right]
\end{aligned} \tag{3.15}$$

The coefficients in eqn. (3.11) are found to be

$$b_{xx} = \frac{1}{16}\left[-R_x^{z-} - R_x^{z+} + \left(\Delta_x V_x^1\right) - \left(\Delta_x V_x^2\right) + \left(\Delta_x V_x^3\right) - \left(\Delta_x V_x^4\right)\right] ;$$

$$b_{zz} = \frac{1}{16}\left[R_z^{x-} + R_z^{x+} + \left(\Delta_z V_z^1\right) + \left(\Delta_z V_z^2\right) - \left(\Delta_z V_z^3\right) - \left(\Delta_z V_z^4\right)\right] ;$$

$$b_{xxz} = \frac{1}{8}\left[R_x^{z-} - R_x^{z+} - \left(\Delta_x V_x^1\right) + \left(\Delta_x V_x^2\right) + \left(\Delta_x V_x^3\right) - \left(\Delta_x V_x^4\right)\right] ;$$

$$b_{xzz} = \frac{1}{8}\left[-R_z^{x-} + R_z^{x+} - \left(\Delta_z V_z^1\right) + \left(\Delta_z V_z^2\right) + \left(\Delta_z V_z^3\right) - \left(\Delta_z V_z^4\right)\right]$$

(3.16)

The coefficients in eqn. (3.12) are found to be

$$c_{xx} = \frac{1}{16}\left[R_x^{y-} + R_x^{y+} + \left(\Delta_x V_x^1\right) + \left(\Delta_x V_x^2\right) - \left(\Delta_x V_x^3\right) - \left(\Delta_x V_x^4\right)\right] ;$$

$$c_{yy} = \frac{1}{16}\left[-R_y^{x-} - R_y^{x+} + \left(\Delta_y V_y^1\right) - \left(\Delta_y V_y^2\right) + \left(\Delta_y V_y^3\right) - \left(\Delta_y V_y^4\right)\right] ;$$

$$c_{xxy} = \frac{1}{8}\left[-R_x^{y-} + R_x^{y+} - \left(\Delta_x V_x^1\right) + \left(\Delta_x V_x^2\right) + \left(\Delta_x V_x^3\right) - \left(\Delta_x V_x^4\right)\right] ;$$

$$c_{xyy} = \frac{1}{8}\left[R_y^{x-} - R_y^{x+} - \left(\Delta_y V_y^1\right) + \left(\Delta_y V_y^2\right) + \left(\Delta_y V_y^3\right) - \left(\Delta_y V_y^4\right)\right]$$

(3.17)

It is easy to see the finite difference-like structure in eqns. (3.15), (3.16) and (3.17). We see that they truly represent higher derivatives that can be derived from the arrangement of gradients at the edges of the control volume in Fig. 3. This also tells us that the effect of these higher derivative terms in eqns. (3.10), (3.11) and (3.12) is to slightly modify the higher moments in those equations so as to restore curl-free or curl-preserving behavior. It is also important to note that, in spite of this modification, the reconstructed vector field will indeed exactly match the mean values and the linear gradients of the vector field at the twelve edges of the zone shown in Fig. 3. As a result, the vector field will be globally constraint-preserving. Note too that the modifications will be slight owing to the fact that the modifying coefficients represent higher derivatives. Therefore, the reconstructed vector field in eqns. (3.10), (3.11) and (3.12) is suitable for the construction of globally curl constraint-preserving schemes.

Because the second order accurate vector field only has to be specified up to linear terms, its discrete curl only needs to be specified up to constant terms when we are considering a second order finite volume scheme. But realize from an examination of eqn. (1.2) and how it arises from the last equation in eqn. (1.1) that a second order DG scheme will also have an evolutionary equation for the curl of the original vector field; and the latter equation also needs to be evolved with second order accuracy. (In two-dimensions, this consideration was immaterial; but in three dimensions it is a meaningful consideration.) This completes our study of curl constraint-preserving reconstruction of vector fields at second order of accuracy.

### III.3) Curl-Preserving Reconstruction on a Three-dimensional Cartesian Mesh at Third Order

To understand the nuances that enter the reconstruction procedure at third order, it is helpful to divide our study into two stages in the following two paragraphs. In the first stage, it is helpful to study curl-free reconstruction of vector fields and the restrictions it places on the various terms in the reconstruction. Only in the second stage will we study curl constraint-preserving reconstruction of vector fields. This two-stage sub-division is also useful because several very useful PDE systems, like general relativity, only require a curl-free reconstruction of vector fields and do not need the computationally heavier details of a curl-preserving reconstruction of vector fields.

In this first stage of our study, let us begin by considering a curl-free vector field that is represented on a three-dimensional Cartesian mesh with third order of accuracy. The component $V^x(x,y,z)$ in eqn. (3.10) has the following modes:- It has a constant mode and it also has modes for $x, y, z, xy, xz, yz, y^2, z^2, y^2z, yz^2, xyz$. Notice that the only mode that is missing in $V^x(x,y,z)$ for obtaining full third order accuracy is the $x^2$ mode. This is good news, because it means that along the four x-edges of Fig. 3 we have to reconstruct/evolve the second moment of the x-component of the vector field along those edges. Let us label those four modes $\left(\Delta_{xx}V_x^1\right)$, $\left(\Delta_{xx}V_x^2\right)$, $\left(\Delta_{xx}V_x^3\right)$ and $\left(\Delta_{xx}V_x^4\right)$ with an obvious extension of the notation. (This extension of notation also applies to the y- and z-directions. ) When an $x^2$ mode is added along each of the four x-edges, it also triggers the additional presence of $x^2y, x^2z, x^2yz$ variation in $V^x(x,y,z)$. Similar considerations apply to the other directions so that we realize that along the four y-edges of Fig. 3 we have to provide four $y^2$-dependent modes to $V^y(x,y,z)$; this also triggers the additional presence of $y^2x, y^2z, y^2xz$ variation in $V^y(x,y,z)$. Furthermore, along the four z-edges of Fig. 3 we have to provide four $z^2$-dependent modes to $V^z(x,y,z)$; this also triggers the additional presence of $z^2 x, z^2 y, z^2 xy$ variation in $V^z(x,y,z)$. To balance all the terms that have to be added in the other two components, we have to add $y^3, z^3, y^3 z, z^3 y$ variations in $V^x(x,y,z)$. Symmetry considerations help us to realize that we will have to add $x^3, z^3, x^3 z, z^3 x$ variations in $V^y(x,y,z)$. Similarly, we will have to add $x^3, y^3, x^3 y, y^3 x$ variations in $V^z(x,y,z)$. These additional terms are the minimum number of terms needed for ensuring the curl-free aspect of the vector field.

In this second stage of our study, we realize that the evolutionary equations for a curl-preserving vector field will have a corresponding evolutionary equation for its curl vector. For instance, see eqn. (1.2) which results from taking the curl of the last equation in eqn. (1.1). Even with all the modes that we have argued for in the previous paragraph, an evaluation of the curl of the original vector field will not have all the terms that are needed for matching all the modes that are specified in the evolutionary equation for the curl. This is especially true for a third order accurate curl-preserving DG-like scheme. In the three equations that follow, such terms are denoted by a capital "A" in $V^x(x,y,z)$; they are denoted by a capital "B" in $V^y(x,y,z)$; and they are denoted by a capital "C" in $V^z(x,y,z)$.

As in the second order case, the polynomials are chosen with such factors that they do not affect the vectors in the edges. This was deemed essential for globally curl-preserving reconstruction. In other words, if the ensuing polynomial for $V^x(x,y,z)$ is evaluated at any x-edge of Fig. 3, we indeed retrieve only the constant, linear and quadratic variations in that x-edge, and nothing but that variation. Similarly, if $V^y(x,y,z)$ is evaluated at any y-edge of Fig. 3, we indeed retrieve only the constant, linear and quadratic variations in that y-edge, and nothing but that variation. Similarly, if $V^z(x,y,z)$ is evaluated at any z-edge of Fig. 3, we indeed retrieve only the constant, linear and quadratic variation in that z-edge, and nothing but that variation. We write

$$V^x(x,y,z) = \left[V_x^1 + \left(\Delta_x V_x^1\right)x + \left(\Delta_{xx} V_x^1\right)\left(x^2 - 1/12\right)\right]\left(y - \frac{1}{2}\right)\left(z - \frac{1}{2}\right)$$

$$+ \left[V_x^2 + \left(\Delta_x V_x^2\right)x + \left(\Delta_{xx} V_x^2\right)\left(x^2 - 1/12\right)\right]\left(y + \frac{1}{2}\right)\left(\frac{1}{2} - z\right)$$

$$+ \left[V_x^3 + \left(\Delta_x V_x^3\right)x + \left(\Delta_{xx} V_x^3\right)\left(x^2 - 1/12\right)\right]\left(\frac{1}{2} - y\right)\left(z + \frac{1}{2}\right)$$

$$+ \left[V_x^4 + \left(\Delta_x V_x^4\right)x + \left(\Delta_{xx} V_x^4\right)\left(x^2 - 1/12\right)\right]\left(y + \frac{1}{2}\right)\left(z + \frac{1}{2}\right) \quad (3.18)$$

$$+ a_{yy}\left(1 - 4y^2\right) + a_{zz}\left(1 - 4z^2\right) + a_{yyz} z\left(1 - 4y^2\right) + a_{yzz} y\left(1 - 4z^2\right)$$

$$+ a_{yyy} y\left(1 - 4y^2\right) + a_{zzz} z\left(1 - 4z^2\right) + a_{yyyz} yz\left(1 - 4y^2\right) + a_{yzzz} yz\left(1 - 4z^2\right)$$

$$+ A_{xyy} x\left(1 - 4y^2\right) + A_{xzz} x\left(1 - 4z^2\right) + A_{xyyz} xz\left(1 - 4y^2\right) + A_{xyzz} xy\left(1 - 4z^2\right)$$

and

$$V^y(x,y,z) = \left[V_y^1 + \left(\Delta_y V_y^1\right)y + \left(\Delta_{yy} V_y^1\right)\left(y^2 - 1/12\right)\right]\left(z - \frac{1}{2}\right)\left(x - \frac{1}{2}\right)$$

$$+ \left[V_y^2 + \left(\Delta_y V_y^2\right)y + \left(\Delta_{yy} V_y^2\right)\left(y^2 - 1/12\right)\right]\left(z + \frac{1}{2}\right)\left(\frac{1}{2} - x\right)$$

$$+ \left[V_y^3 + \left(\Delta_y V_y^3\right)y + \left(\Delta_{yy} V_y^3\right)\left(y^2 - 1/12\right)\right]\left(\frac{1}{2} - z\right)\left(x + \frac{1}{2}\right)$$

$$+ \left[V_y^4 + \left(\Delta_y V_y^4\right)y + \left(\Delta_{yy} V_y^4\right)\left(y^2 - 1/12\right)\right]\left(z + \frac{1}{2}\right)\left(x + \frac{1}{2}\right) \quad (3.19)$$

$$+ b_{xx}\left(1 - 4x^2\right) + b_{zz}\left(1 - 4z^2\right) + b_{xxz} z\left(1 - 4x^2\right) + b_{xzz} x\left(1 - 4z^2\right)$$

$$+ b_{xxx} x\left(1 - 4x^2\right) + b_{zzz} z\left(1 - 4z^2\right) + b_{xxxz} xz\left(1 - 4x^2\right) + b_{xzzz} xz\left(1 - 4z^2\right)$$

$$+ B_{xxy} y\left(1 - 4x^2\right) + B_{yzz} y\left(1 - 4z^2\right) + B_{xxyz} yz\left(1 - 4x^2\right) + B_{xyzz} xy\left(1 - 4z^2\right)$$

and

$$V^z(x,y,z) = \left[V_z^1 + (\Delta_z V_z^1)z + (\Delta_{zz} V_z^1)(z^2 - 1/12)\right]\left(x - \frac{1}{2}\right)\left(y - \frac{1}{2}\right)$$

$$+ \left[V_z^2 + (\Delta_z V_z^2)z + (\Delta_{zz} V_z^2)(z^2 - 1/12)\right]\left(x + \frac{1}{2}\right)\left(\frac{1}{2} - y\right)$$

$$+ \left[V_z^3 + (\Delta_z V_z^3)z + (\Delta_{zz} V_z^3)(z^2 - 1/12)\right]\left(\frac{1}{2} - x\right)\left(y + \frac{1}{2}\right)$$

$$+ \left[V_z^4 + (\Delta_z V_z^4)z + (\Delta_{zz} V_z^4)(z^2 - 1/12)\right]\left(x + \frac{1}{2}\right)\left(y + \frac{1}{2}\right) \quad (3.20)$$

$$+ c_{xx}(1 - 4x^2) + c_{yy}(1 - 4y^2) + c_{xxy} y(1 - 4x^2) + c_{xyy} x(1 - 4y^2)$$

$$+ c_{xxx} x(1 - 4x^2) + c_{yyy} y(1 - 4y^2) + c_{xxxy} xy(1 - 4x^2) + c_{xyyy} xy(1 - 4y^2)$$

$$+ C_{xxz} z(1 - 4x^2) + C_{yyz} z(1 - 4y^2) + C_{xxyz} yz(1 - 4x^2) + C_{xyyz} xz(1 - 4y^2)$$

Now let us study the curl-preserving constraints that are imposed on the above vector field. The object of this study is to derive insights that help us to fix up the 36 coefficients that are as yet unspecified in eqns. (3.18), (3.19) and (3.20) above.

From the Introduction, we have seen that every curl-preserving, edge-centered vector field has to be such that its curl matches a divergence-free vector field. This divergence-free vector field has components that are face-centered. Therefore, at the top and bottom x-faces of Fig. 3 we have the mean circulations given by $(V_y^3 + V_z^4 - V_y^4 - V_z^2)$ and $(V_y^1 + V_z^3 - V_y^2 - V_z^1)$ respectively. Using two-dimensional, third order, finite-volume WENO reconstruction that is restricted to the x-faces of the mesh, we can write the facial variation of the x-component of the circulation in the top and bottom x-faces as

$$R^{x+}(y,z) = (V_y^3 + V_z^4 - V_y^4 - V_z^2) + R_y^{x+} y + R_z^{x+} z + R_{yy}^{x+}(y^2 - 1/12) + R_{zz}^{x+}(z^2 - 1/12) + R_{yz}^{x+} yz \ ;$$
$$R^{x-}(y,z) = (V_y^1 + V_z^3 - V_y^2 - V_z^1) + R_y^{x-} y + R_z^{x-} z + R_{yy}^{x-}(y^2 - 1/12) + R_{zz}^{x-}(z^2 - 1/12) + R_{yz}^{x-} yz$$
(3.21a)

Similarly, at the top and bottom y-faces of Fig. 3 we have the mean circulations given by $(V_z^3 + V_x^4 - V_z^4 - V_x^2)$ and $(V_z^1 + V_x^3 - V_z^2 - V_x^1)$ respectively. Using two-dimensional, third order, finite-volume WENO reconstruction that is restricted to the y-faces of the mesh, we can write the facial variation of the y-component of the circulation in the top and bottom y-faces as

$$R^{y+}(x,z) = \left(V_z^3 + V_x^4 - V_z^4 - V_x^2\right) + R_x^{y+}x + R_z^{y+}z + R_{xx}^{y+}\left(x^2 - 1/12\right) + R_{zz}^{y+}\left(z^2 - 1/12\right) + R_{xz}^{y+}xz ;$$
$$R^{y-}(x,z) = \left(V_z^1 + V_x^3 - V_z^2 - V_x^1\right) + R_x^{y-}x + R_z^{y-}z + R_{xx}^{y-}\left(x^2 - 1/12\right) + R_{zz}^{y-}\left(z^2 - 1/12\right) + R_{xz}^{y-}xz$$
(3.21b)

Likewise, at the top and bottom z-faces of Fig. 3 we have the mean circulations given by $\left(V_x^3 + V_y^4 - V_x^4 - V_y^2\right)$ and $\left(V_x^1 + V_y^3 - V_x^2 - V_y^1\right)$ respectively. Using two-dimensional, third order, finite-volume WENO reconstruction that is restricted to the z-faces of the mesh, we can write the facial variation of the z-component of the circulation in the top and bottom z-faces as

$$R^{z+}(x,y) = \left(V_x^3 + V_y^4 - V_x^4 - V_y^2\right) + R_x^{z+}x + R_y^{z+}y + R_{xx}^{z+}\left(x^2 - 1/12\right) + R_{yy}^{z+}\left(y^2 - 1/12\right) + R_{xy}^{z+}xy ;$$
$$R^{z-}(x,y) = \left(V_x^1 + V_y^3 - V_x^2 - V_y^1\right) + R_x^{z-}x + R_y^{z-}y + R_{xx}^{z-}\left(x^2 - 1/12\right) + R_{yy}^{z-}\left(y^2 - 1/12\right) + R_{xy}^{z-}xy$$
(3.21c)

Notice that WENO has given us the facial variations of the x-, y- and z-components of the circulations in the x-, y- and z-faces of the mesh. But we still do not have expressions for the circulation vector field over the volume of the zone that we are considering.

We can obtain expressions for the x-, y- and z-components of the circulation vector field by using the divergence-free reconstruction strategy catalogued in Balsara *et al*. [23] (see Section 3 of that paper). The x-component of the volumetrically reconstructed circulation can be written as

$$\begin{aligned}R^x(x,y,z) &= \alpha_0 + \alpha_x x + \alpha_y y + \alpha_z z + \alpha_{xx}\left(x^2 - 1/12\right) + \alpha_{yy}\left(y^2 - 1/12\right) + \alpha_{zz}\left(z^2 - 1/12\right)\\&+ \alpha_{xy}xy + \alpha_{yz}yz + \alpha_{xz}xz + \alpha_{xxx}\left(x^3 - 3x/20\right) + \alpha_{xxy}\left(x^2 - 1/12\right)y + \alpha_{xxz}\left(x^2 - 1/12\right)z\\&+ \alpha_{xyy}x\left(y^2 - 1/12\right) + \alpha_{xzz}x\left(z^2 - 1/12\right) + \alpha_{xyz}xyz\end{aligned}$$
(3.22a)

The y-component of the volumetrically reconstructed circulation can be written as

$$\begin{aligned}R^y(x,y,z) &= \beta_0 + \beta_x x + \beta_y y + \beta_z z + \beta_{xx}\left(x^2 - 1/12\right) + \beta_{yy}\left(y^2 - 1/12\right) + \beta_{zz}\left(z^2 - 1/12\right)\\&+ \beta_{xy}xy + \beta_{yz}yz + \beta_{xz}xz + \beta_{yyy}\left(y^3 - 3y/20\right) + \beta_{xyy}x\left(y^2 - 1/12\right) + \beta_{yyz}\left(y^2 - 1/12\right)z\\&+ \beta_{xxy}\left(x^2 - 1/12\right)y + \beta_{yzz}y\left(z^2 - 1/12\right) + \beta_{xyz}xyz\end{aligned}$$
(3.22b)

The z-component of the volumetrically reconstructed circulation can be written as

$$R^z(x,y,z) = \gamma_0 + \gamma_x x + \gamma_y y + \gamma_z z + \gamma_{xx}(x^2 - 1/12) + \gamma_{yy}(y^2 - 1/12) + \gamma_{zz}(z^2 - 1/12)$$
$$+ \gamma_{xy} xy + \gamma_{yz} yz + \gamma_{xz} xz + \gamma_{zzz}(z^3 - 3z/20) + \gamma_{xzz} x(z^2 - 1/12) + \gamma_{yzz} y(z^2 - 1/12) \qquad (3.22c)$$
$$+ \gamma_{xxz}(x^2 - 1/12)z + \gamma_{yyz}(y^2 - 1/12)z + \gamma_{xyz} xyz$$

Using the divergence-free constraint, Section 3 of Balsara *et al.* [23] gives us an explicit procedure for matching the 45 coefficients in eqns. (3.22a), (3.22b) and (3.22c) to the coefficients in eqns. (3.21a), (3.21b) and (3.21c). We do not repeat the procedure here because all the formulae have already been explicitly catalogued in the literature. Realize that we now have divergence-free expressions for the circulation vector field over the volume of the zone that we are considering.

With the divergence-free expressions for the circulation vector field in hand, we can now match the curl of eqns. (3.18), (3.19) and (3.20) to eqns. (3.22a), (3.22b) and (3.22c). This gives us a sparse linear system which we can solve for the 36 free coefficients in eqns. (3.18), (3.19) and (3.20). With a little ingenuity, the coefficients in eqn. (3.18) are found to be

$$a_{yy} = \frac{1}{16}\left[R_y^{z-} + R_y^{z+} + (\Delta_y V_y^1) + (\Delta_y V_y^2) - (\Delta_y V_y^3) - (\Delta_y V_y^4)\right] ;$$

$$a_{zz} = \frac{1}{16}\left[-R_z^{y-} - R_z^{y+} + (\Delta_z V_z^1) - (\Delta_z V_z^2) + (\Delta_z V_z^3) - (\Delta_z V_z^4)\right] ;$$

$$a_{yyy} = \frac{1}{24}\left[R_{yy}^{z-} + R_{yy}^{z+} + (\Delta_{yy} V_y^1) + (\Delta_{yy} V_y^2) - (\Delta_{yy} V_y^3) - (\Delta_{yy} V_y^4)\right] ;$$

$$a_{zzz} = \frac{1}{24}\left[-R_{zz}^{y-} - R_{zz}^{y+} + (\Delta_{zz} V_z^1) - (\Delta_{zz} V_z^2) + (\Delta_{zz} V_z^3) - (\Delta_{zz} V_z^4)\right] ;$$

$$a_{yyz} = \frac{1}{8}\left[-R_y^{z-} + R_y^{z+} - (\Delta_y V_y^1) + (\Delta_y V_y^2) + (\Delta_y V_y^3) - (\Delta_y V_y^4)\right] ;$$

$$a_{yzz} = \frac{1}{8}\left[R_z^{y-} - R_z^{y+} - (\Delta_z V_z^1) + (\Delta_z V_z^2) + (\Delta_z V_z^3) - (\Delta_z V_z^4)\right] ;$$

$$A_{xzz} = -\frac{1}{32}\left[R_{xz}^{y-} + R_{xz}^{y+}\right] ; \quad A_{xyy} = \frac{1}{32}\left[R_{xy}^{z-} + R_{xy}^{z+}\right] ;$$

$$a_{yyyz} = \frac{1}{12}\left[-R_{yy}^{z-} + R_{yy}^{z+} - (\Delta_{yy} V_y^1) + (\Delta_{yy} V_y^2) + (\Delta_{yy} V_y^3) - (\Delta_{yy} V_y^4)\right] ;$$

$$a_{yzzz} = \frac{1}{12}\left[R_{zz}^{y-} - R_{zz}^{y+} - (\Delta_{zz} V_z^1) + (\Delta_{zz} V_z^2) + (\Delta_{zz} V_z^3) - (\Delta_{zz} V_z^4)\right] ; \qquad (3.23)$$

$$A_{xyyz} = \frac{1}{16}\left[-R_{xy}^{z-} + R_{xy}^{z+}\right] ; \quad A_{xyzz} = \frac{1}{16}\left[R_{xz}^{y-} - R_{xz}^{y+}\right]$$

The coefficients in eqn. (3.19) are found to be

$$b_{xx} = \frac{1}{16}\left[-R_x^{z-} - R_x^{z+} + \left(\Delta_x V_x^1\right) - \left(\Delta_x V_x^2\right) + \left(\Delta_x V_x^3\right) - \left(\Delta_x V_x^4\right)\right] ;$$

$$b_{zz} = \frac{1}{16}\left[R_z^{x-} + R_z^{x+} + \left(\Delta_z V_z^1\right) + \left(\Delta_z V_z^2\right) - \left(\Delta_z V_z^3\right) - \left(\Delta_z V_z^4\right)\right] ;$$

$$b_{xxx} = \frac{1}{24}\left[-R_{xx}^{z-} - R_{xx}^{z+} + \left(\Delta_{xx} V_x^1\right) - \left(\Delta_{xx} V_x^2\right) + \left(\Delta_{xx} V_x^3\right) - \left(\Delta_{xx} V_x^4\right)\right] ;$$

$$b_{zzz} = \frac{1}{24}\left[R_{zz}^{x-} + R_{zz}^{x+} + \left(\Delta_{zz} V_z^1\right) + \left(\Delta_{zz} V_z^2\right) - \left(\Delta_{zz} V_z^3\right) - \left(\Delta_{zz} V_z^4\right)\right] ;$$

$$b_{xxz} = \frac{1}{8}\left[R_x^{z-} - R_x^{z+} - \left(\Delta_x V_x^1\right) + \left(\Delta_x V_x^2\right) + \left(\Delta_x V_x^3\right) - \left(\Delta_x V_x^4\right)\right] ;$$

$$b_{xzz} = \frac{1}{8}\left[-R_z^{x-} + R_z^{x+} - \left(\Delta_z V_z^1\right) + \left(\Delta_z V_z^2\right) + \left(\Delta_z V_z^3\right) - \left(\Delta_z V_z^4\right)\right] ;$$

$$B_{yzz} = \frac{1}{32}\left[R_{yz}^{x-} + R_{yz}^{x+}\right] ; \quad B_{xxy} = -\frac{1}{32}\left[R_{xy}^{z-} + R_{xy}^{z+}\right] ;$$

$$b_{xxxz} = \frac{1}{12}\left[R_{xx}^{z-} - R_{xx}^{z+} - \left(\Delta_{xx} V_x^1\right) + \left(\Delta_{xx} V_x^2\right) + \left(\Delta_{xx} V_x^3\right) - \left(\Delta_{xx} V_x^4\right)\right] ;$$

$$b_{xzzz} = \frac{1}{12}\left[-R_{zz}^{x-} + R_{zz}^{x+} - \left(\Delta_{zz} V_z^1\right) + \left(\Delta_{zz} V_z^2\right) + \left(\Delta_{zz} V_z^3\right) - \left(\Delta_{zz} V_z^4\right)\right] ;$$

$$B_{xxyz} = \frac{1}{16}\left[R_{xy}^{z-} - R_{xy}^{z+}\right] ; \quad B_{xyzz} = \frac{1}{16}\left[-R_{yz}^{x-} + R_{yz}^{x+}\right]$$

(3.24)

The coefficients in eqn. (3.20) are found to be

$$c_{xx} = \frac{1}{16}\left[ R_x^{y-} + R_x^{y+} + \left(\Delta_x V_x^1\right) + \left(\Delta_x V_x^2\right) - \left(\Delta_x V_x^3\right) - \left(\Delta_x V_x^4\right)\right] \ ;$$

$$c_{yy} = \frac{1}{16}\left[ -R_y^{x-} - R_y^{x+} + \left(\Delta_y V_y^1\right) - \left(\Delta_y V_y^2\right) + \left(\Delta_y V_y^3\right) - \left(\Delta_y V_y^4\right)\right] \ ;$$

$$c_{xxx} = \frac{1}{24}\left[ R_{xx}^{y-} + R_{xx}^{y+} + \left(\Delta_{xx} V_x^1\right) + \left(\Delta_{xx} V_x^2\right) - \left(\Delta_{xx} V_x^3\right) - \left(\Delta_{xx} V_x^4\right)\right] \ ;$$

$$c_{zzz} = \frac{1}{24}\left[ -R_{yy}^{x-} - R_{yy}^{x+} + \left(\Delta_{yy} V_y^1\right) - \left(\Delta_{yy} V_y^2\right) + \left(\Delta_{yy} V_y^3\right) - \left(\Delta_{yy} V_y^4\right)\right] \ ;$$

$$c_{xxy} = \frac{1}{8}\left[ -R_x^{y-} + R_x^{y+} - \left(\Delta_x V_x^1\right) + \left(\Delta_x V_x^2\right) + \left(\Delta_x V_x^3\right) - \left(\Delta_x V_x^4\right)\right] \ ;$$

$$c_{xyy} = \frac{1}{8}\left[ R_y^{x-} - R_y^{x+} - \left(\Delta_y V_y^1\right) + \left(\Delta_y V_y^2\right) + \left(\Delta_y V_y^3\right) - \left(\Delta_y V_y^4\right)\right] \ ;$$

$$C_{yyz} = -\frac{1}{32}\left[ R_{yz}^{x-} + R_{yz}^{x+} \right] \ ; \quad C_{xxz} = \frac{1}{32}\left[ R_{xz}^{y-} + R_{xz}^{y+} \right] \ ;$$

$$c_{xxxy} = \frac{1}{12}\left[ -R_{xx}^{y-} + R_{xx}^{y+} - \left(\Delta_{xx} V_x^1\right) + \left(\Delta_{xx} V_x^2\right) + \left(\Delta_{xx} V_x^3\right) - \left(\Delta_{xx} V_x^4\right)\right] \ ;$$

$$c_{xyyy} = \frac{1}{12}\left[ R_{yy}^{x-} - R_{yy}^{x+} - \left(\Delta_{yy} V_y^1\right) + \left(\Delta_{yy} V_y^2\right) + \left(\Delta_{yy} V_y^3\right) - \left(\Delta_{yy} V_y^4\right)\right] \ ;$$

$$C_{xxyz} = \frac{1}{16}\left[ -R_{xz}^{y-} + R_{xz}^{y+} \right] \ ; \quad C_{xyyz} = \frac{1}{16}\left[ R_{yz}^{x-} - R_{yz}^{x+} \right]$$

(3.25)

As in the second order case, it is easy to see the finite difference-like structure in eqns. (3.23), (3.24) and (3.25). We see that they truly represent higher derivatives that can be derived from the arrangement of gradients (and higher moments) at the edges of the control volume in Fig. 3. Because the third order accurate vector field only has to be specified up to quadratic terms, we might be able to get away with specifying its discrete curl only up to linear terms when we are considering a third order WENO-like scheme. But realize from an examination of eqn. (1.2) and how it arises from the last equation in eqn. (1.1) that a third order DG scheme will also have an evolutionary equation for the curl of the original vector field; and the latter equation also needs to be evolved with third order accuracy. This completes our study of curl constraint-preserving reconstruction of vector fields at third order of accuracy.

**III.4) Combining Curl-Free Reconstruction and the Three-Dimensional Riemann Solver to Obtain a Multidimensionally Upwinded, Globally Curl-Free Scheme**

By now it is very evident that curl-free reconstruction of vector fields requires a collocation of vector components at the edges of the mesh. This is true in two and three dimensions. In Sub-section II.5 we showed that in two dimensions, the use of curl-free reconstruction, in conjunction with a two-dimensional Riemann solver that provides the two-dimensional upwinding, can indeed result in a properly upwinded, globally curl-free scheme. Since the curious reader might wonder whether there is an analogous extension to three dimensions, we provide such an extension here. To keep the discussion generally applicable to any possible time stepping strategy, we present it within the context of a semi-discrete formulation in time.

Let us focus on the model system that is the three dimensional extension of the one we studied in eqn. (2.22). We have

$$\frac{\partial J^x}{\partial t} + \frac{\partial \left( v^x J^x + v^y J^y + v^z J^z \right)}{\partial x} + v^y \left( \frac{\partial J^x}{\partial y} - \frac{\partial J^y}{\partial x} \right) + v^z \left( \frac{\partial J^x}{\partial z} - \frac{\partial J^z}{\partial x} \right) = 0 \;;$$

$$\frac{\partial J^y}{\partial t} + \frac{\partial \left( v^x J^x + v^y J^y + v^z J^z \right)}{\partial y} + v^x \left( \frac{\partial J^y}{\partial x} - \frac{\partial J^x}{\partial y} \right) + v^z \left( \frac{\partial J^y}{\partial z} - \frac{\partial J^z}{\partial y} \right) = 0 \;; \quad (3.26)$$

$$\frac{\partial J^z}{\partial t} + \frac{\partial \left( v^x J^x + v^y J^y + v^z J^z \right)}{\partial z} + v^x \left( \frac{\partial J^z}{\partial x} - \frac{\partial J^x}{\partial z} \right) + v^y \left( \frac{\partial J^z}{\partial y} - \frac{\partial J^y}{\partial z} \right) = 0$$

Analogous to the demonstration in Sub-section II.5, it is easy to show that if a curl-free vector field is initially provided, and if a curl-free reconstruction is used, then the semi-discrete form of the above equations ought to evolve the vector field in curl-free fashion. Additionally, we require that the potential be collocated at the vertices of the mesh. The above system is of great practical interest because it arises naturally as part of the first order CCZ4 hyperbolic system that has to be solved for numerical general relativity; in fact, in the general relativistic context the equations do not have source terms. As in Sub-section II.5, the emphasis in this Sub-section is on curl-free evolution of the three-dimensional vector field because the inclusion of source terms on the right hand sides of the above equations can easily turn them into curl-preserving equations. For the simple case where the velocity vector $\left( v^x, v^y, v^z \right)$ is a constant, eqn. (3.26) can be simplified to yield

$$\frac{\partial J^x}{\partial t} + v^x \frac{\partial J^x}{\partial x} + v^y \frac{\partial J^x}{\partial y} + v^z \frac{\partial J^x}{\partial z} = 0 \ ;$$

$$\frac{\partial J^y}{\partial t} + v^x \frac{\partial J^y}{\partial x} + v^y \frac{\partial J^y}{\partial y} + v^z \frac{\partial J^y}{\partial z} = 0 \ ; \qquad (3.27)$$

$$\frac{\partial J^z}{\partial t} + v^x \frac{\partial J^z}{\partial x} + v^y \frac{\partial J^z}{\partial y} + v^z \frac{\partial J^z}{\partial z} = 0 \ ;$$

While the above equations suggest that we are simply advecting a vector field in three dimensions, please note that eqn. (3.26) enjoins us to keep the vector field globally curl-free. This requires the edge-centered collocation of the components of the vector field. Eqn. (3.26) also tells us that the potentials "$v^x J^x + v^y J^y + v^z J^z$" should be uniquely specified at the vertices of the mesh in a multidimensionally upwinded fashion. To obtain stable advection, we need to show that the update methodology is also properly upwinded. We demonstrate all these facets within the context of a first order semi-discrete scheme.

**Place Fig. 4 here**

Fig. 4 shows the collocation of curl-free vector components along the edges of a three-dimensional zone. The zone center is indexed by *(i,j,k)* and the edges are indexed suitably, consistent with the zone center's indexing. As in the two-dimensional case, the potentials are defined by $\phi \equiv v^x J^x + v^y J^y + v^z J^z$ and they are collocated at the vertices of the mesh. As long as the same potential at a vertex is used for the update of all the vector components in all the edges that meet at that vertex, the update will be globally curl-free. To keep the discussion simple, we take all the velocity components to be constant and positive. All the zones of the Cartesian mesh are also taken to be uniform with mesh sizes $\Delta x$, $\Delta y$ and $\Delta z$ in the x-, y- and z-directions. The upwinded potentials at two of the vertices of the mesh are also shown. The potentials at other vertices can be obtained by suitable shifts in the indexing. The purpose of this figure is to make it easy for us to understand how a curl-free reconstruction that is based on edge-centered vector components, in conjunction with a three dimensional Riemann solver, can give us a stable, globally curl-free scheme.

In the ensuing discussion, we shall focus on only the first of the equations in eqns. (3.26) and (3.27) because manipulations that are identical to the ones shown below can be made for the other two components of the vector field **J** . Since the vector field is curl-free, the terms $v^y \left( \partial J^x / \partial y - \partial J^y / \partial x \right)$ and $v^z \left( \partial J^x / \partial z - \partial J^z / \partial x \right)$ do not contribute. As a result, the semi-discrete update for $J^x_{i,j+1/2,k+1/2}$ in Fig. 4 becomes

$$\frac{\partial J^x_{i,j+1/2,k+1/2}}{\partial t} + \frac{\phi_{i+1/2,j+1/2,k+1/2} - \phi_{i-1/2,j+1/2,k+1/2}}{\Delta x} = 0 \tag{3.28}$$

Using the potentials $\phi_{i+1/2,j+1/2,k+1/2} = v^x J^x_{i,j+1/2,k+1/2} + v^y J^y_{i+1/2,j,k+1/2} + v^z J^z_{i+1/2,j+1/2,k}$ and $\phi_{i-1/2,j+1/2,k+1/2} = v^x J^x_{i-1,j+1/2,k+1/2} + v^y J^y_{i-1/2,j,k+1/2} + v^z J^z_{i-1/2,j+1/2,k}$ from Fig. 4, we get

$$\begin{aligned}\frac{\partial J^x_{i,j+1/2,k+1/2}}{\partial t} &+ \frac{v^x}{\Delta x} \left( J^x_{i,j+1/2,k+1/2} - J^x_{i-1,j+1/2,k+1/2} \right) + \frac{v^y}{\Delta x} \left( J^y_{i+1/2,j,k+1/2} - J^y_{i-1/2,j,k+1/2} \right) \\ &+ \frac{v^z}{\Delta x} \left( J^z_{i+1/2,j+1/2,k} - J^z_{i-1/2,j+1/2,k} \right) = 0\end{aligned} \tag{3.29}$$

It is easy to see that eqn. (3.29) scarcely resembles the first equation in eqn. (3.27). However, we now use the discrete circulations in faces $(i,j,k+1/2)$ and $(i,j+1/2,k)$ of Fig. 4 to make the transcription

$$\begin{aligned}\frac{1}{\Delta x} \left( J^y_{i+1/2,j,k+1/2} - J^y_{i-1/2,j,k+1/2} \right) &\to \frac{1}{\Delta y} \left( J^x_{i,j+1/2,k+1/2} - J^x_{i,j-1/2,k+1/2} \right) \\ \frac{1}{\Delta x} \left( J^z_{i+1/2,j+1/2,k} - J^z_{i-1/2,j+1/2,k} \right) &\to \frac{1}{\Delta z} \left( J^x_{i,j+1/2,k+1/2} - J^x_{i,j+1/2,k-1/2} \right)\end{aligned} \tag{3.30}$$

Putting the transcription from eqn. (3.30) in eqn. (3.29) we get

$$\begin{aligned}\frac{\partial J^x_{i,j+1/2,k+1/2}}{\partial t} &+ \frac{v^x}{\Delta x} \left( J^x_{i,j+1/2,k+1/2} - J^x_{i-1,j+1/2,k+1/2} \right) + \frac{v^y}{\Delta y} \left( J^x_{i,j+1/2,k+1/2} - J^x_{i,j-1/2,k+1/2} \right) \\ &+ \frac{v^z}{\Delta z} \left( J^x_{i,j+1/2,k+1/2} - J^x_{i,j+1/2,k-1/2} \right) = 0\end{aligned} \tag{3.31}$$

The concordance between eqn. (3.31) and the first equation in eqn. (3.27) is now very obvious. We see that at first order our globally curl-free scheme is indeed properly upwinded. For the update in eqn. (3.31) we should use an effective control volume of size $\Delta x \times \Delta y \times \Delta z$ that is centered on the center of the edge $(i, j+1/2, k+1/2)$ shown in Fig. 4. As we use higher order curl-free or curl-preserving reconstructions and higher order timestepping, we have the assurance that the numerical dissipation will be progressively reduced with increasing order but the scheme will remain stable. The demonstration might require a few steps, but the previous equation can also be written in a form that makes the centered and dissipation parts self-evident. We don't show the steps, but we show the final result. We have

$$\frac{\partial J^x_{i,j+1/2,k+1/2}}{\partial t} + \frac{v^x}{2\Delta x}\left(J^x_{i+1,j+1/2,k+1/2} - J^x_{i-1,j+1/2,k+1/2}\right) + \frac{v^y}{2\Delta y}\left(J^x_{i,j+3/2,k+1/2} - J^x_{i,j-1/2,k+1/2}\right)$$
$$+ \frac{v^z}{2\Delta z}\left(J^x_{i,j+1/2,k+3/2} - J^x_{i,j+1/2,k-1/2}\right) = \frac{|v^x|}{2\Delta x}\left(J^x_{i+1,j+1/2,k+1/2} - 2J^x_{i,j+1/2,k+1/2} + J^x_{i-1,j+1/2,k+1/2}\right) \quad (3.32)$$
$$+ \frac{|v^y|}{2\Delta y}\left(J^x_{i,j+3/2,k+1/2} - 2J^x_{i,j+1/2,k+1/2} + J^x_{i,j-1/2,k+1/2}\right) + \frac{|v^z|}{2\Delta z}\left(J^x_{i,j+1/2,k+3/2} - 2J^x_{i,j+1/2,k+1/2} + J^x_{i,j+1/2,k-1/2}\right)$$

We now see that a higher order extension can be made, as detailed in the steps given in the next sub-section.

### III.5) Stepwise Strategy for Implementing a Curl-Preserving Scheme

Thus we can identify the three essential steps for implementing a curl-preserving scheme:-

**1)** We make a higher order curl-free reconstruction of the sort that is presented in Sections II and III for any curl-constrained vector fields. Flow variables that have a zone-centered interpretation in a traditional higher order Godunov scheme can be reconstructed using well-known TVD or WENO methods. This gives us suitably high order spatial reconstruction in all instances.

**2)** We continue to have the support of one-dimensional Riemann solver technology. However, on a two-dimensional mesh, we also use the two-dimensional Riemann solver to give us the potential $\phi$ at the corners of the mesh, as shown in Fig. 2. On a three dimensional mesh, like the one shown in Fig. 4, we can use the three-dimensional Riemann solver from Balsara [16] to give us the potential $\phi$ at the vertices of the mesh. We should also incorporate any suitable integration of source terms. For example, it is evident that the terms $\mathbf{v} \times (\nabla \times \mathbf{J})$ and $-\frac{\rho T}{\tau} \mathbf{J}$ in eqn. (1.1d) should be averaged at an edge from the zones that come together at the edge of interest. This step and the previous step give us a single stage of a Runge-Kutta scheme.

**3)** A higher order SSP-RK time-stepping strategy can then be used to obtain higher order temporal accuracy. This combination of innovations will produce a stable, high order, multidimensionally upwinded, globally curl-free (or curl-preserving) scheme.

**IV) von Neumann Stability Analysis of Curl Constraint-Preserving WENO-like Schemes**

A von Neumann stability analysis for *classical, volume-centered, DG* schemes has been done (Liu *et al*. [44], Zhang and Shu [58]). In that analysis, the authors focused on the advection equation with a constant velocity. While the above authors showed a DG scheme, their methods can be adapted as well for WENO schemes. In a classical DG/WENO scheme, the primal variables are zone-centered and the objective is to satisfy a *telescoping density-conservation constraint*. (In other words, mass, momentum and energy densities are conserved in any subset of zones because of a telescoping application of mass, momentum and energy fluxes.) A similar von Neumann stability analysis for WENO-like and DG-like schemes for evolving the induction equation in MHD has also been carried out by Balsara and Käppeli [20]. The induction equation evolves a vector field in divergence-free fashion, and the analysis can be carried out with the simplification of a constant velocity in two dimensions. With that simplification, the induction equation also reduces to an advection of the two components of the vector field. In a *divergence-preserving, face-centered, WENO-like or DG-like* scheme for the induction equation, the primal variables are face-centered and the objective is to satisfy a *telescoping divergence-preserving constraint*. This choice of collocation also holds true for any divergence-constraint preserving PDE like Maxwell's equations (Balsara and Käppeli [24]). It is now easy to see that for the *curl-constraint preserving, edge-centered, WENO-like* schemes for treating eqn. (2.22), we have the objective to satisfy a *telescoping curl-preserving constraint*. When we restrict eqn. (2.22) so as to have a constant velocity, eqns. (2.23) show us that the model PDE again reduces to the advection of two curl-free vector field components. (The reader should note that it is extremely difficult to analyze a full 2D scheme with non-constant velocity and that is the only reason for choosing a constant velocity. All the schemes that are analyzed in this paper can be applied to practical problems that have non-constant and time-evolving velocities, as shown later in Sections V and VI.) In all such stability analyses it is traditional to simplify the spatial and temporal parts of the problem by using multi-stage Runge-Kutta timestepping. Therefore, we will use the SSP-RK timestepping schemes from (Shu and Osher [52], [53], Shu [54], Spiteri and Ruuth [50], [51], Gottlieb *et al*. [40]). The temporal order in our von Neumann stability analyses will always be matched to the spatial order of accuracy of the WENO-like scheme. (Because this paper is focused on WENO methods, we only develop the von Neumann stability analysis for curl-free WENO-like schemes for treating eqn. (2.22). In subsequent work, we would like to develop the same stability analysis for curl-free DG-like and PNPM-like schemes.)

We saw in Sub-sections II.5, III.4 and III.5 that the ingredients of a successful curl constraint-preserving scheme consist of :- 1) a higher order curl-free reconstruction, 2) a multidimensional Riemann solver and 3) a suitable high order timestepping strategy. It is possible to use these three building blocks to design WENO-like, PNPM-like and DG-like schemes of higher order Godunov type that preserve the global constraints. The full description of such a plan requires indeed a separate paper and such a paper is under construction (Balsara and Käppeli [28]). The description of DG-like schemes is more intricate and is not done here; however, we have indeed described WENO-like schemes in detail in this paper. Therefore, it makes sense to present the von Neumann stability analysis of edge-centered curl-preserving schemes here.

From the edge-centered primal variables, the reconstruction strategy described in Section II (and its associated Fig. 1) is used to reconstruct the entire vector field. For a von Neumann stability analysis, smooth flow variables can be assumed, so that we always use results from the suitably high order central stencil from the WENO reconstruction. The application of the multidimensional Riemann solver, along with the use of a suitably higher order SSP-RK timestepping scheme then completes the update strategy. In our von Neumann stability analysis we used a 2D Cartesian mesh with square zones. In such a von Neumann stability analysis, one posits a curl-free vector field in 2D with wave vector $\left(k_x, k_y\right)$ which has harmonic variation of the form $e^{i\left(k_x x + k_y y\right)}$ . For the model problem shown in eqn. (2.22) we then obtain the amplification factor as well as the phase of the entire globally curl-free WENO-like scheme. The exclusive use of the central stencil makes the entire scheme linear in the edge-centered variables with the result that a computer algebra system can be used to extract the entire amplification matrix.

Notice that each different choice of velocity components $(v_x, v_y)$ and each different choice of wave vector $(k_x, k_y)$ yields a different amplification factor and phase. What matters is that the amplification factor and phase depend on:- 1) The angle between the velocity vector and the x-direction of the mesh, 2) the relative angle between the wave vector and the velocity and 3) the ratio of the wavelength to the mesh size. We consider situations where the velocity vector makes angles of 0º, 15º, 30º and 45º relative to the x-direction of the mesh. This gives us a sufficiently interesting range of velocity directions. For each of those velocity directions we allow the wave vector to sweep over all possible angles between the direction of the velocity and the direction of the wave vector. Higher order WENO schemes usually display good resolving capabilities, so we display the amplitude and phase information when the wavelength spans 5 zones, 10 zones and 15 zones. We do this for second, third and fourth order accurate WENO-like schemes.

**Place Fig. 5 here**

We study the wave propagation characteristics for globally curl-free WENO-like schemes, which are second order accurate in space. The time-stepping was SSP-RK2 so that the entire scheme is spatially and temporally second order accurate. The von Neumann stability analysis shows us that this combination of spatial and temporal discretization yields a maximum CFL of 0.7071.

For the plots shown in Fig. 5, the CFL was set at 90% of its maximum value. Figs. 5a to 5d show one minus the absolute value of the amplification factor when the velocity vector makes angles of 0º, 15º, 30º and 45º relative to the x-direction of the mesh. Figs. 5e to 5h show the phase error, again for the same angles. The 2D wave vector can make any angle relative to the 2D direction of velocity propagation, with the result that the amplitude and phase information are shown with respect to the angle made between the velocity direction and the direction of the wave vector. In each plot, the blue curve refers to waves that span 5 cells per wavelength; the green curve refers to waves that span 10 cells per wavelength; the red curve refers to waves that span 15 waves per wavelength. Notice that the cases where the velocity vector makes angles of 0º and 45º relative to the x-direction of the mesh indeed show symmetrical wave propagation characteristics, as expected. When the velocity vector makes angles of 15º and 30º relative to the x-direction of the mesh, there is no symmetry between the velocity direction, the mesh direction and the direction of the wave vector, with the result that we don't expect to see symmetrical plots, and indeed we don't. We do, however, observe that when the waves span 10 cells per wavelength and 15 cells per wavelength the wave propagation becomes very close to isotropic and quite free of dissipation. This is a good sign that even our second order WENO-like scheme shows rather isotropic wave propagation with increasing wave length. In all instances, Figs. 5a to 5d show us that one minus the amplification factor is always positive or zero, indicating that the globally curl-free, second order, WENO-like scheme is indeed stable.

**Place Fig. 6 here**

We study the wave propagation characteristics for globally curl-free WENO-like schemes, which are spatially third order accurate. The time-stepping was SSP-RK3 so that the entire scheme is spatially and temporally third order accurate. The von Neumann stability analysis shows us that this combination of spatial and temporal discretization yields a maximum CFL of 1.1507. Note that a full application would involve zone-centered fluid-like variables which may need to be evolved with their own WENO reconstruction. That will impose its own CFL restriction, with the result that the smaller of the two CFL numbers has to be chosen.

For the plots shown in Fig. 6, the CFL was set at 90% of its maximum value. The first four panels in Fig. 6 show one minus the amplification factor, and the next four panels in Fig. 6 show the phase error. Comparing these results to their analogues in Fig. 5 we see that the third order scheme shows a significant improvement in the phase accuracy compared to the second order scheme. Owing to the built-in dissipation properties of SSP-RK3 time-stepping, the dissipation in Fig 6 is only competitive with that in Fig. 5. However, the improving trend is restored once we go to fourth order. This shows the benefit of resorting to a higher order WENO-like scheme. We also see that the wave propagation in Fig. 6 is much more isotropic compared to Fig. 5.

**Place Fig. 7 here**

We study the wave propagation characteristics for globally curl-free WENO-like schemes, which are spatially fourth order accurate. The time-stepping was SSP-RK(5,4) so that the entire scheme is spatially and temporally fourth order accurate. The von Neumann stability analysis shows us that this combination of spatial and temporal discretization yields a maximum CFL of 1.3040. As before, when using this scheme along with a formulation that requires zone-centered fluid-like variables, one has to choose the smaller of the permissible CFL numbers.

For the plots shown in Fig. 7, the CFL was set at 90% of its maximum value. The first four panels in Fig. 7 show one minus the amplification factor, and the next four panels in Fig. 7 show the phase error. Comparing these results to analogous results in Figs. 5 or 6, we see that the transition to fourth order of accuracy has made a very substantial improvement in amplification factor as well as the phase accuracy. As before, comparing Fig. 7 to Figs. 5 or 6 shows that the fourth order WENO-like scheme has not just vastly improved accuracy but also significantly improved isotropy of wave propagation. Especially when there are ten or more zones per wavelength, Fig. 7 shows that the dissipation and dispersion errors of a fourth order WENO-like scheme are roughly one order of magnitude lower than that of a second or third order scheme. This shows that in multidimensions, higher order WENO-like curl-free schemes provide not just vastly improved accuracy but also significantly improved isotropy of wave propagation.

Taken together, the results of this Section show that our curl constraint-preserving reconstruction, coupled with the multidimensional Riemann solver, provides a successful framework for the design of a very proficient class of mimetic WENO-like schemes for involution-constrained PDEs. Most importantly, the results presented point to a class of high order WENO - like mimetic schemes for involution-constrained PDEs that have superior amplitude preservation and phase accuracy even in multiple dimensions. Analysis of DG-like and PNPM-like schemes will, it is hoped, show even further possibilities for improvement.

### V) Numerical Results for Two Model Problems; and a Demonstration of Order Property

It behooves us to design and display a model problem where one can palpably witness the value of an exactly curl-preserving scheme. Moreover, since we have presented entire classes of such schemes with increasing order of accuracy, we would like to demonstrate the value of high order of accuracy. To that end, we present a test problem for a simple model PDE system, where the curl-free evolution of the vector field is of crucial importance.

Let us begin our discussion by considering the currently-available alternatives. Of course, a GLM-style cleaning procedure has been developed in Dumbser *et al.* [36], but it requires increasing signal speed of the cleaning equations and also adds many more vector fields than are originally necessary. In Dumbser *et al.* [37] and Boscheri *et al.* [30] a new exactly curl-free semi-implicit scheme was presented using a vertex-based *staggered* mesh, but it is limited to second order of accuracy and requires frequent interpolation of the velocity field and of the curl-free vector field **J** to different staggered locations on the mesh, which makes it more difficult to extend to adaptive mesh refinement (AMR) techniques.

By contrast, the present formulation preserves the same control volume for the fluid variables as well as the curl-constraint preserving vector field, making it suitable for an eventual future extension to AMR. The availability of higher order curl-preserving formulations also allows us to show another interesting facet that has gone unappreciated in the literature. It turns out that in certain important limits, the curl-preserving vector field satisfies an energy principle. The quadratic energy of the vector field should remain unchanged in time. A good scheme should preserve this quadratic energy as much as possible. We show in this section that our increasingly accurate curl-preserving schemes preserve the quadratic energy with increasing precision.

In this Section we illustrate the capabilities of our new high order accurate numerical method with the help of the toy system introduced in Dumbser *et al.* [36], [37]. The governing PDE system reads

$$\frac{\partial \rho}{\partial t} + \frac{\partial}{\partial x_i}(\rho v_i) = 0 \tag{5.1}$$

$$\frac{\partial \rho v_k}{\partial t} + \frac{\partial}{\partial x_i}\left(\rho v_i v_k + \delta_{ik} p + \rho c_0^2 J_i J_k\right) = 0 \tag{5.2}$$

$$\frac{\partial J_k}{\partial t} + \frac{\partial}{\partial x_k}(v_m J_m) + v_m\left(\frac{\partial J_k}{\partial x_m} - \frac{\partial J_m}{\partial x_k}\right) = 0 \tag{5.3}$$

where the Einstein convention implying summation over repeated indexes is adopted. The system of equations (5.1)-(5.3) describes the evolution of a scalar quantity ρ and two vector fields **v** and **J**, that in a fluid dynamic context could be interpreted as density, velocity and a kind of thermal impulse, respectively, while p represents the equivalent of a pressure, see Dumbser et al. [33]. As shown in Dumbser et al. [37] the system satisfies the extra energy conservation law

$$\frac{\partial \rho E}{\partial t} + \frac{\partial}{\partial x_k}\left(v_k(\rho E + p) + v_i \rho c_0^2 J_i J_k\right) = 0. \tag{5.4}$$

In this paper the model system is closed by the simple linear relation

$$p = \gamma^2 \rho \tag{5.5}$$

where $\gamma$ is a given constant, as well as $c_0$ in eqn. (5.2). The system (5.1)-(5.3) with (5.4) falls into the larger class of symmetric hyperbolic and thermodynamically compatible (SHTC) systems studied by Godunov and Romenski, see Godunov [38] and Romenski [47] and references therein.

Now, let us focus on the third PDE (5.3) and let $\chi = v_m J_m$ be a scalar quantity. Applying the Schwarz theorem, which implies the symmetry of second derivatives, i.e.

$$\frac{\partial}{\partial x_k}\frac{\partial}{\partial x_m}\chi - \frac{\partial}{\partial x_m}\frac{\partial}{\partial x_k}\chi = 0,$$

it follows that eqn. (6.3) maintains the linear involution constraint

$$C_{mk} = \left(\frac{\partial J_k}{\partial x_m} - \frac{\partial J_m}{\partial x_k}\right) = 0 \tag{5.6}$$

for all times if the field **J** was curl-free at the initial time.

It is therefore crucial to satisfy this constraint even at the discrete level, that is if $C_{mk} = 0$ at the initial time it must remain zero for all times.

### V.1) Model Problem: A Stationary curl-free solution

To verify that the novel curl constraint-preserving scheme is able to fulfill over time the involution constraint in eqn. (5.6) we propose to solve the following test problem, which is an exact stationary and smooth solution of the model system (5.1)-(5.3). Let the computational domain be the square $\Omega=[-5;5]^2$ and let the generic radial coordinate r satisfy $r^2=x^2+y^2$. The quantity **J** is defined as the *gradient* of a scalar potential $\phi$, so that the initial condition is ensured to be curl-free. The potential is given by

$$\varphi(r) = A \operatorname{erf}\left(\frac{r-R_0}{\sigma}\right), \tag{5.7}$$

with the parameters A, $R_0$ and $\sigma$. Then, the initial condition for the radial component $J_r(r) = \mathbf{J}(r) \cdot \mathbf{e}_r$ reads

$$J_r(r) = \frac{\partial \varphi}{\partial r} = \frac{2A}{\sqrt{\pi}\sigma} \exp\left(-\frac{(r-R_0)^2}{\sigma^2}\right), \tag{5.8}$$

while the angular component is set to $J_\phi = \mathbf{J} \cdot \mathbf{e}_\phi$. Here, $\mathbf{e}_r$ and $\mathbf{e}_\phi$ are the unit vectors in the radial and in the angular direction, respectively.

The initial condition in eqn. (5.8) guarantees that the involution constraint in eqn. (5.6) is satisfied at the initial time t=0. We furthermore impose **v**=0 at the initial time. From radial direction of eqn. (5.2) rewritten in polar coordinates it then follows that a stationary equilibrium is preserved if

$$\frac{d}{dr}\left(p(r) + \rho(r)c_0^2 J_r^2\right) + \frac{1}{r}\rho(r)c_0^2 J_r^2 = 0, \tag{5.9}$$

Solving the above equilibrium condition for the radial derivative of $p(\rho)$ and using eqn. (5.5) yields the following non-autonomous ODE for $\rho(r)$:

$$\frac{d\rho}{dr} = -\frac{\rho(r)J(r)c_0^2}{\gamma^2 + J^2 c_0^2}\left(2\frac{dJ}{dr} + \frac{J(r)}{r}\right), \qquad \rho(0) = \rho_0 \tag{5.10}$$

The ODE (5.10) can be solved numerically to obtain the initial condition for the density profile $\rho(r)$, which completes the setup of the initial condition of this test case. To this purpose we use a classical fourth order Runge-Kutta ODE solver with a very fine mesh spacing in radial direction so that the solution of the ODE can be considered as quasi exact. The parameters in this model problem are set to A=0.2, $R_0$=2, $\sigma$=0.5, $\rho_0$=2, $c_0$=2 and $\gamma$=2.

**V.2) Accuracy Analysis of the Numerical Scheme Using our Model Problem**

**Place Fig. 8 here**

In Figure 8 we show the numerical results for the vector component $J_x$ obtained on a mesh with 50 x 50 cells by running four different schemes: the semi-implicit second order accurate staggered curl-free (SCF) scheme proposed in (Boscheri *et al.* [30]) and the second, third and fourth order accurate edge centered curl-preserving (ECCP) reconstruction methods developed in Section II of this work. The final time of the simulation is very large, i.e. $t_{end}$=100, in order to show the behavior and the stability of the scheme for very long time computations. The less dissipative behavior achieved by the high order order reconstructions is clearly visible. We also notice that the SCF scheme is the most dissipative of the schemes shown in Fig. 8 because of the copious interpolation of the velocity field and the curl-free vector field **J** to different staggered locations on the mesh. We also mentioned that the quadratic energy of the vector field should, in principle, be preserved for this physical problem. All numerical schemes fall short of this ideal goal. In Fig. 9 we plot out the mesh-integrated quadratic energy, simply evaluated as $J^2 = J_x^2 + J_y^2$ as a function of time. We see that the higher order schemes preserve the quadratic energy much better than the lower order schemes. This is because the numerical viscosity is significantly reduced when a high order reconstruction technique is adopted. Finally, we would also like to point out that without a curl-preserving scheme the solution is spoiled after short times and the equilibrium is violated very soon leading to catastrophically unphysical results. At late times, the result of not treating the curl-preserving aspect of the PDE is indeed a code blowup.

**Place Fig. 9 here**

A numerical convergence study for this test problem is carried out by solving the model PDE system with the stationary initial condition discussed above until a final time of t=10 using a sequence of successively refined meshes. The results for our novel curl-free WENO schemes presented in this paper are shown in Table II for nominal orders of accuracy from two to four. In Fig. 9 we finally show the temporal evolution of the curl error of second, third and fourth order curl-free WENO schemes on a uniform Cartesian mesh composed of 32 x 32 elements. It can be clearly seen that in all cases and for all times, the curl errors remain at the level of machine precision, as expected.

TABLE II. Numerical convergence study from second to fourth order accuracy for the novel high order curl-preserving schemes that draw on the Curl-Preserving reconstruction from Section II presented in this paper. Errors for the variable $J_x$ are shown.

| Method | $N_x$ x $N_y$ | $L_1$ Error | $L_1$ Order | $L_\infty$ Error | $L_\infty$ Order |
|---|---|---|---|---|---|
| Curl-Preserving O2 | | | | | |
| | 64 x 64 | 2.4453E-3 | | 4.7175E-2 | |
| | 128 x 128 | 5.0732E-4 | **2.3** | 1.7351E-2 | **1.4** |
| | 256 x 256 | 9.2619E-5 | **2.4** | 7.4771E-3 | **1.2** |
| | 512 x 512 | 1.5713E-5 | **2.6** | 1.8191E-3 | **2.0** |
| Curl-Preserving O3 | | | | | |
| | 64 x 64 | 2.0530E-3 | | 3.8879E-2 | |
| | 128 x 128 | 4.6588E-4 | **2.1** | 1.0150E-2 | **1.9** |
| | 256 x 256 | 6.9384E-5 | **2.8** | 1.6952E-3 | **2.6** |
| | 512 x 512 | 9.0052E-6 | **3.0** | 2.1711E-4 | **3.0** |
| Curl-Preserving O4 | | | | | |
| | 64 x 64 | 3.7155E-4 | | 9.1343E-3 | |
| | 128 x 128 | 1.6224E-5 | **4.5** | 4.7515E-4 | **4.2** |
| | 256 x 256 | 1.6224E-5 | **4.3** | 2.4617E-5 | **4.3** |
| | 512 x 512 | 5.3898E-8 | **3.9** | 1.9865E-6 | **3.6** |

**Place Fig. 10 here**

Taken together, the results of this Sub-section show that our curl constraint-preserving reconstruction, coupled with the multidimensional Riemann solver, provides a successful framework for the design of mimetic finite volume schemes of increasing order of accuracy. Moreover, these schemes preserve the curl-constraint during very long time integrations. This constraint-preservation also contributes significantly to the enhanced stability of the scheme. The model problem that we provide here is also quite novel and it enables us to precisely document that the methods presented here do indeed meet their designed order of accuracy. The utility of mimetic schemes with high accuracy is also emphasized by the fact that additional quadratic energy terms are also preserved with superlative precision as one goes to higher order. As a result, we have presented high order mimetic finite volume-type schemes which have long time stability and excellent preservation of quadratic energy.

### V.4) Model Problem: Inhomogeneous Curl Involution

This test case aims at demonstrating the consistency of the time evolution related to the curl of the vector field ***J***. Let us consider a non-homogeneous curl involution in the toy model by adding a source term in eqn. (5.3) so that the evolutionary equation for the x- and y-components can be explicitly written as

$$\frac{\partial J_x}{\partial t} + \frac{\partial}{\partial x}\left(v_x J_x + v_y J_y\right) + v_y \left(\frac{\partial J_x}{\partial y} - \frac{\partial J_y}{\partial x}\right) = S_x \;;$$

$$\frac{\partial J_y}{\partial t} + \frac{\partial}{\partial y}\left(v_x J_x + v_y J_y\right) + v_x \left(\frac{\partial J_y}{\partial x} - \frac{\partial J_x}{\partial y}\right) = S_y \quad (5.11)$$

with $S_x = \sin x \sin y$ and $S_y = \cos x \cos y$

The initial conditions are given by setting the density field to unity and the components of the velocity field to be zero throughout the computational domain. The initial conditions for the thermal impulse vector field are given by

$$J_x = \sin x \sin y \quad \text{and} \quad J_y = \cos x \cos y \quad (5.12)$$

As a result, the curl of the thermal impulse vector assumes a prescribed non-zero initial value given by $R_z = -2\sin x \cos y$. Similar to eqn. (1.2), an additional PDE is then defined which accounts for the time evolution of the Burger's vector field $\mathbf{R} = \nabla \times \mathbf{J}$. In two-dimensions, the only non-zero component of $R_z$ whose evolution equation can be explicitly written as

$$\frac{\partial R_z}{\partial t} + \frac{\partial}{\partial x}\left(v_x R_z - S_y\right) + \frac{\partial}{\partial y}\left(v_y R_z + S_x\right) = 0 \qquad (5.13)$$

The constant $c_0$ in eqn. (5.2) was set to zero. In eqn. (5.5), we set $\gamma = 2$ so that we have a constitutive relation between the density and the pressure. Notice that in addition to having initial conditions given by $R_z = -2\sin x \cos y$, because the velocity field is zero, eqn. (5.13) becomes

$$\frac{\partial R_z}{\partial t} = -2\sin x \cos y \qquad (5.14)$$

In other words, the Burger's vector is analytically integrable! The goal of this test problem is to verify that Burger's vector, which is the curl of the thermal impulse vector, is correctly evolved by the scheme.

**Place Fig. 11 here**

The computational domain is taken to be a 100×100 Cartesian zone mesh in two-dimensions that spans $[-\pi, 3\pi] \times [-\pi, 3\pi]$. Using the third order curl-preserving scheme described in this paper (see Sub-Section III.5), we have run this test problem to a final time of unity. A CFL of 0.6 was used for the edge-centered curl-preserving scheme. (The reason it can sustain large CFL numbers is because it uses the multidimensional Riemann solver.) We also built a straightforward third order, zone-centered Godunov scheme that evolves eqns. (5.1), (5.2), (5.11) and (5.13) so that we can compare and contrast the curl constraint-preserving scheme with one that does not preserve the constraints. Please recall that because of eqn. (5.14), the Burger's vector is analytically integrable. Fig. 11a shows the analytically evaluated Burger's vector at a final time of unity. Fig. 11b shows the Burger's vector at the same final time when we used the third order accurate edge-centered curl-preserving scheme. Fig. 11c shows the Burger's vector at the same final time when we used a plain-vanilla, zone-centered, third order accurate Godunov scheme and directly finite differenced the zone-centered vector components $J_x$ and $J_y$ to obtain their curl. Although Figs. 11a, 11b and 11c look similar, please examine the quantitative values by using the color bars. We see that the third order accurate edge-centered curl-preserving scheme in Fig. 11b has closely tracked the analytical result in Fig. 11a. However, the numerical values in Fig. 11c are substantially different from those in Figs. 11a or 11b. This shows that a scheme that does not accurately account for the constraints in a curl-constrained PDE will indeed generate some noticeable errors even after a modest amount of time-integration.

**Place Fig. 12 here**

We can also make the errors more quantitative when we realize that eqn. (5.13) is analytically integrable. We can, therefore, evaluate the error in the curl at each timestep as

$$Error = \iint_{(x,y)\in[-\pi,3\pi]^2} |(\nabla \times \mathbf{J}) - R_z| dA \qquad (5.15)$$

Fig. 12 shows the error in the time-evolution of the Burger's vector when we use the third order accurate edge-centered curl-preserving scheme and when we use a plain-vanilla, zone-centered, third order accurate Godunov scheme. If the methods in this paper are not used, the errors in the curl are bigger by 2 to 4 orders of magnitude. This is significant because the curl indeed drives the evolution of the thermal impulse in eqn. (1.1d).

## VI) Results from Further Test Problems Involving the full GPR System from Eqn. (1.1)

While the previous Section showed several useful results with the help of a toy problem, we now turn our attention to the full GPR system from eqn. (1.1). Our test consists of a blast problem that is allowed to evolve with various extents of thermal conduction. The thermal conduction is controlled by the relaxation time $\tau$ in eqn. (1.1d) with the result that smaller values of $\tau$ result in the thermal conduction having an increasingly larger effect. Because we have run all our simulations with a time-explicit formulation of the source terms, we have restricted ourselves to $\tau = 100$, $\tau = 10$, $\tau = 1$ and $\tau = 0.1$. As for the other coefficients that regulate the PDE system in eqn. (1.1), we have used $\gamma = 1.4$, $c_h = 1$ and $c_v = 2.5$. Our objective is to carry out four simulations with the four different values of $\tau$ using the third order accurate, edge-centered, curl-preserving formulation with one-dimensional and multidimensional Riemann solvers that produces constrained evolution. We will then repeat the same four simulations with a plain-vanilla, third order accurate, zone-centered higher order Godunov scheme, with one-dimensional Riemann solver technology, where nothing special is done to account for the constrained evolution of the thermal impulse vector field. We will then intercompare the two methods and bring to the forefront the simulations where the maximum differences have been found.

The blast problem with thermal conduction that we consider here is set up on a 200×200 zone mesh that spans $(x, y) \in [-1,1]^2$. Continuative or periodic boundary conditions may be used, because the simulation is always stopped before the outgoing blast wave reaches the boundary. Let "$r$" be the radius, measured from the origin. The problem consists of initializing the primitive variables $(\rho, v_x, v_y, P, j_x, j_y)$ so that the central circular region is at higher density and pressure compared to the ambient region, as follows

$$(\rho, v_x, v_y, P, j_x, j_y) = \begin{cases} (2,0,0,1,0,0) & \text{for } r \leq 0.2 \\ (0.5,0,0,0.5,0,0) & \text{for } r > 0.2 \end{cases}$$

The simulation is run to a time of 0.7, which is just before the outer shock encounters the boundary. The edge-centered, curl-preserving scheme was run with a CFL of 0.6, whereas the plain-vanilla Godunov scheme was run with a CFL of 0.45. Notice that the interior of this blast wave set-up is cooler than the exterior, even if it is at higher pressure.

**Place Fig. 13 here**

**Place Fig. 14 here**

**Place Fig. 15 here**

Figs. 13a, 13b, 13c and 13d show the density, temperature, x-velocity and x-component of the thermal impulse at a final time of 0.7 for a simulation with $\tau = 1$ that was run with the third order accurate, edge-centered, curl-preserving formulation with one-dimensional and multidimensional Riemann solvers. Most of our other simulations look similar to Fig. 13, indicating that when the relaxation time is relatively large, the fluxes are dominated by the fluid fluxes and thermal conduction does not play a dominant role. However, Fig. 14 shows the same

variables as in Fig. 13, but with $\tau = 0.1$. The simulation in Fig. 14 was run with the same algorithm as Fig. 13. Fig. 15 shows another simulation with $\tau = 0.1$, but this time it was run with a plain-vanilla, third order accurate, zone-centered higher order Godunov scheme, with one-dimensional Riemann solver technology. All densities across Figs. 13 to 15 are shown with the same color table that spans $[0.45, 2]$; this makes it easy to inter-compare results across the three figures. All temperatures (defined by $P/\rho$) across Figs. 13 to 15 are shown with the same color table that spans $[0.5, 1.05]$. Likewise, the x-velocity across the three figures uses a color table that has a range of $[-0.18, 0.18]$. Similarly, the x-component of the thermal impulse uses a color table that has a range of $[-0.75, 0.75]$ across all three figures. The simulation rapidly sets up a system of outward-going circular shock and an inward-propagating reverse shock. This two-shock structure is well-known even for a purely hydrodynamical blast wave problem. The final time corresponds to a time when the inward-propagating reverse shock converges on to the origin. Because the interior of the blast is cooler than the exterior, a heat flux is set up that carries heat from the exterior to the interior of the blast. In other words, while the blast expands outwards, the heat propagates inwards. As a result, one can observe that the sign of the outward-going x-velocity is opposite to the sign of the inward-going x-component of the thermal impulse vector field in these three figures. Decreasing values of $\tau$ set up increasing amounts of heat flux. That heat flux has to compete with the other fluid dynamical fluxes in the problem with the result that we only expect small values of $\tau$ to produce heat fluxes that produce an appreciable difference relative to the fluid dynamical fluxes. Because the color bars are the same across flow variables in these three figures, let us now compare Fig. 14 with Fig. 15. We see that the density, temperature and the x-component of the thermal impulse are quite different in Fig. 14 compared to Fig. 15. Specifically, the x-component of the thermal impulse, which tracks the heat flux, is indeed very different. This shows that measurable differences have revealed themselves which can only be attributed to the difference in algorithms. We see that the curl-preserving formulation does make a substantial difference in the physical result!

**Place Fig. 16 here**

Now let us try to understand why the choice of algorithm plays a very important role in differentiating the results. To make this apparent, we plot $\left|\nabla\times\mathbf{J}\right|_{\max}$ as a function of simulation time for all our simulations of the blast problem in this Section. The results are shown in Fig. 16 in color-coded format. The solid curves show the results of the edge-centered curl-preserving algorithm. The dashed curves show the results of the plain-vanilla higher order Godunov algorithm. The dramatic difference is immediately evident. With diminishing values of $\tau$, i.e. with increasingly stronger source terms in eqn. (1.1d), we see that the edge-centered curl-preserving algorithm produces a progressively larger maximum in the curl of the thermal impulse. In other words, the outputs are proportional to, and regulated by, the input value of $\tau$! This is exactly what we would desire in a well-designed numerical experiment! Please also note from eqn. (1.1d) that the $\mathbf{v}\times\left(\nabla\times\mathbf{J}\right)$ term drives the build-up of the thermal impulse vector as long as there is a non-zero source term on the right-hand side of that equation. The non-zero source term in eqn. (1.1d) allows the evolution to cease being curl-free and makes it curl-preserving. But it is the non-zero value of the $\mathbf{v}\times\left(\nabla\times\mathbf{J}\right)$ term in eqn. (1.1d) that drives the further build-up of circulation in the thermal impulse field. The edge-centered curl-preserving algorithm does a perfect job of modulating all aspects of that build-up so that decreasing values of $\tau$ give us a well-modulated increase in the thermal flux! Now let us look at the dashed curves in Fig. 16. We see that regardless of the value of $\tau$, the schemes that give rise to the dashed curves all produce the about same amount of circulation. Besides, that circulation builds up very quickly and in a completely unregulated fashion. This is because the circulation always tracks the discretization errors which can be quite large due to the zone-centered collocation that is inherent in a plain-vanilla higher order Godunov scheme. The need for a curl-preserving formulation has been driven home to us thanks to Fig. 16.

**VII) Conclusions**

Structure-preserving PDEs are becoming increasingly important for several applications in physics and engineering. Of particular interest in this paper are PDEs that preserve the curl of a vector field. Examples of such PDEs include hyperelasticity, compressible multiphase flows with and without surface tension, first order reductions of the Einstein field equations as well as the novel first order hyperbolic reformulation of the Schrödinger equation, to name a few examples.

We have designed methods in this paper for increasingly high order numerical treatment of such PDEs. The essential building blocks for such methods are shown to be a non-linearly hybridized curl constraint-preserving, high order accurate reconstruction strategy and the use of multidimensional Riemann solvers that are needed for properly upwinded constraint-preserving time update. These two building blocks are then coupled to the third building block which is the SSP-RK family of time-stepping strategies.

Sections II and III show how curl-preserving reconstruction is carried out in two and three dimensions. The starting point is a one-dimensional WENO reconstruction along the edges. However, careful attention has to be paid to the curl constraints in order to get the reconstructed vector field at all locations within the zone of interest. This is why we think of the reconstruction as a WENO-like reconstruction. Notice that we use non-linearly hybridized WENO schemes to build the higher order moments in the edges and only then carry out the volumetric reconstruction according to principles of curl-preservation. As a result, our curl-preserving reconstruction is non-linearly hybridized and, therefore, suitable for integration with higher order Godunov scheme philosophy. In those sections we also demonstrate that when the reconstruction is combined with multidimensional Riemann solvers, we get numerical schemes that are multidimensionally upwinded and multidimensionally stable.

We refer to the reconstruction described above as edge-centered curl-preserving (ECCP) reconstruction. If all the higher moments in each edge are obtained with WENO reconstruction then they can be referred to as WENO-ECCP schemes. If all the moments in each edge are evolved according to DG-like principles, then they can be referred to as DG-ECCP schemes. If only some of the lower moments in each edge are evolved, while higher moments are reconstructed, then they can be referred to as PNPM-ECCP schemes.

Section IV presents some of the results of a von Neumann stability analysis of globally structure preserving WENO-like schemes in multiple space dimensions. The results in Section IV point to a class of high order WENO-like mimetic schemes for involution-constrained PDEs that have superior amplitude preservation and phase accuracy even in multiple dimensions.

In Section V a test problem is constructed that produces steady-state analytically exact solutions. In that Section we show that the mimetic finite volume schemes that use our methods indeed preserve order of accuracy while simultaneously satisfying the curl-free constraint. In some limits, these schemes also preserve the quadratic energy on the mesh. The utility of our mimetic schemes with high accuracy is also illustrated by the fact that additional quadratic energy terms are preserved with superlative precision as one goes to higher order. As a result, we have presented high order accurate mimetic finite volume-type schemes which have long time stability and excellent preservation of quadratic energy. Furthermore, we have also presented examples that show that our curl constraint-preserving schemes retain fidelity with known analytical results, while zone-centered higher order Godunov schemes lose that fidelity. Section VI shows further test problems involving thermal conductivity for the full GPR system in eqn. (1.1). We show that if the numerical scheme is not curl-preserving, eventually some measurable deficiencies will reveal themselves in the simulations.

In summary, in this paper we have shown that curl-constraint preserving WENO-like reconstruction provides a rich set of insights for mimetic scheme design when dealing with PDEs that have a curl-preserving involution.

**Acknowledgements**

DSB acknowledges support via NSF grants NSF-19-04774, NSF-AST-2009776 and NASA-2020-1241. MD acknowledges the financial support received from the Italian Ministry of Education, University and Research (MIUR) in the frame of the Departments of Excellence Initiative 2018-2022 attributed to DICAM of the University of Trento (grant L. 232/2016) and in the frame of the PRIN 2017 project *Innovative numerical methods for evolutionary partial differential equations and applications*. M.D. has also received funding from the University of Trento via the Strategic Initiative *Modeling and Simulation*. WB acknowledges funding from the Istituto Nazionale di Alta Matematica (INdAM) through the GNCS group and the program Young Researchers Funding 2018 via the research project *Semi-implicit structure-preserving schemes for continuum mechanics*.

**Figure Captions**

*Fig. 1 shows the collocation of vector components along the edges of a two-dimensional control volume. As evaluated over the edges of the square element, the discrete circulation is fully specified. The mean value and its linear variation are shown along each edge, in anticipation of a second order accurate reconstruction scheme. The reconstruction problem for a curl-free reconstruction consists of obtaining a polynomial-based vector field that is globally curl-free within this two-dimensional control volume. The reconstruction problem for a curl-preserving reconstruction consists of obtaining a polynomial-based vector field that matches the specified mean circulation in the zone.*

*Fig. 2 shows the components of the curl-free vector field around the four zones (i,j), (i-1,j), (i-1,j-1) and (i,j-1). A first order curl-free reconstruction is used. The multidimensionally upwinded potentials at the vertices of the zone (i,j) are also shown. The red dashed rectangle shows the effective control volume that is used for the update of $J_{i,j-1/2}^{x;n}$ . The blue dashed rectangle shows the effective control volume that is used for the update of $J_{i-1/2,j}^{y;n}$ .*

*Fig. 3 shows the collocation of vector components along the edges of the control volume. Within each of the two x-faces, the two y-faces and the two z-faces, the discrete circulation (evaluated over those faces) is either exactly zero or specified. The mean value and its linear variation are shown along each edge, in anticipation of a second order accurate reconstruction scheme. The reconstruction problem consists of obtaining a polynomial that is globally curl-free/curl-preserving within this control volume.*

*Fig. 4 shows the collocation of curl-free vector components along the edges of a three-dimensional zone. The zone center is indexed by (i,j,k) and the edges are indexed suitably, consistent with the zone center's indexing. We take all the velocity components to be constant and positive. The*

*upwinded potentials at two of the vertices of the mesh are also shown. The potentials at other vertices can be obtained by suitable shifts in the indexing. The purpose of this figure is to make it easy for us to understand how a curl-free reconstruction that is based on edge-centered vector components, in conjunction with a three dimensional Riemann solver, can give us a stable, globally curl-free scheme.*

*Fig. 5 shows the wave propagation characteristics for curl-preserving second order WENO-like schemes. Figs. 5a to 5d show one minus the absolute value of the amplification factor when the velocity vector makes angles of 0º , 15º , 30º and 45º relative to the x-direction of the 2D mesh. Figs. 5e to 5h show the phase error, again for the same angles. The 2D wave vector can make any angle relative to the 2D direction of velocity propagation, therefore, the amplitude and phase information are shown w.r.t. the angle made between the velocity direction and the direction of the wave vector. In each plot, the blue curve refers to waves that span 5 cells per wavelength; the green curve refers to waves that span 10 cells per wavelength; the red curve refers to waves that span 15 waves per wavelength.*

*Fig. 6 shows the wave propagation characteristics for curl-preserving third order WENO-like schemes. Figs. 6a to 6d show one minus the absolute value of the amplification factor when the velocity vector makes angles of 0º , 15º , 30º and 45º relative to the x-direction of the 2D mesh. Figs. 6e to 6h show the phase error, again for the same angles. The 2D wave vector can make any angle relative to the 2D direction of velocity propagation, therefore, the amplitude and phase information are shown w.r.t. the angle made between the velocity direction and the direction of the wave vector. In each plot, the blue curve refers to waves that span 5 cells per wavelength; the green curve refers to waves that span 10 cells per wavelength; the red curve refers to waves that span 15 waves per wavelength.*

*Fig. 7 shows the wave propagation characteristics for curl-preserving fourth order WENO-like schemes. Figs. 7a to 7d show one minus the absolute value of the amplification factor when the velocity vector makes angles of 0º , 15º , 30º and 45º relative to the x-direction of the 2D mesh.*

*Figs. 7e to 7h show the phase error, again for the same angles. The 2D wave vector can make any angle relative to the 2D direction of velocity propagation, therefore, the amplitude and phase information are shown w.r.t. the angle made between the velocity direction and the direction of the wave vector. In each plot, the blue curve refers to waves that span 5 cells per wavelength; the green curve refers to waves that span 10 cells per wavelength; the red curve refers to waves that span 15 waves per wavelength.*

*Fig. 7 shows the wave propagation characteristics for curl-preserving fourth order WENO-like schemes. Figs. 7a to 7d show one minus the absolute value of the amplification factor when the velocity vector makes angles of 0º , 15º , 30º and 45º relative to the x-direction of the 2D mesh. Figs. 7e to 7h show the phase error, again for the same angles. The 2D wave vector can make any angle relative to the 2D direction of velocity propagation, therefore, the amplitude and phase information are shown w.r.t. the angle made between the velocity direction and the direction of the wave vector. In each plot, the blue curve refers to waves that span 5 cells per wavelength; the green curve refers to waves that span 10 cells per wavelength; the red curve refers to waves that span 15 waves per wavelength.*

*Fig 8 shows four different simulations of our model test problem at 50X50 zone resolution. All figures show the y-component of the curl-free vector field after a very long simulation time; i.e. stopping time of 100. The color scheme is the same across all figures and the intervals between contours is also the same. Fig. 8a shows the result of the second-order curl free scheme from Boscheri et al. (2020) which frequently averages the solution to different locations of the mesh in the course of a timestep. Figs. 8b, 8c and 8d show the results from our new scheme that uses curl-free reconstruction at second, third and fourth order, along with the multidimensional Riemann solver. The color saturation is seen to increase with increasing order and the solution is also less diffusive as one progresses from Fig. 8a to 8d. Very different second order schemes were used for Figs. 8a and 8b; we see that Fig. 8b shows appreciably lower diffusion.*

*Fig. 9 shows the evolution of the quadratic energy in the curl-free vector field as a function of time for the four schemes shown in Fig. 8. The model problem was run with 50X50 zone resolution to a time of 100. In principle, the quadratic energy of the constrained vector field integrated over the mesh should stay unchanged. The plot shows how much the quadratic energy decays with time for various schemes. Because of the frequent averaging over the course of a timestep, the second order SCF(O2) scheme from Dumbser et al. [32] shows appreciable decay of energy. The second order ECCP(O2) scheme designed here improves on the SCF(O2) scheme. The third order ECCP(O3) scheme shows substantial improvement over the second order schemes. The fourth order ECCP(O4) scheme provides the best preservation of quadratic energy, thus highlighting the value of well-designed higher order curl-preserving schemes for the evolution of curl-constrained vector fields.*

*Fig. 10 shows the time series of the maximum pointwise error of the curl of **J** for the stationary curl free test problem until t=50. The computational domain is discretized with 32x32 uniform Cartesian cells. The results for second (black solid line), third (blue solid line) and fourth (red solid line) order curl-free WENO schemes are shown.*

*Figs 11a, b, c show the colorized plot of Rz as a function of position in the test problem at a time of 1. Because the problem lends itself to analytical treatment, Fig. 11a was obtained analytically. Fig. 11b was obtained numerically. Using the color bar, please note the close concordance in values between Figs. 11a and 11b. Fig. 11c was obtained from direct differentiation of a zone-centered higher order Godunov scheme that is not curl-preserving. Please use the color bars to compare the numerical values in Fig. 11c to those in Fig. 11a to see the significant differences. The color bar for each panel is different and scaled to the min and max of the data shown.*

*Fig. 12 shows the maximum error in Rz as a function of time for the edge-centered curl-preserving scheme versus a higher order Godunov scheme that does not preserve the curl. The latter shows errors that are 2 to 4 magnitudes higher.*

*Fig. 13a, 13b, 13c, 13d show density, temperature, x-velocity and x-component of the thermal impulse when t=1 was used along with the ECCP scheme at third order. The solution is shown at a final time of 0.7. All densities across Figs. 13 to 15 are shown on the same scale to facilitate inter-comparison. Likewise, for all pressures across Figs. 13 to 15 . Similarly for all x-velocities and all x-components of the thermal impulse across Figs. 13 to 15.*

*Fig. 14a, 14b, 14c, 14d show density, temperature, x-velocity and x-component of the thermal impulse when t=0.1 was used along with the ECCP scheme at third order. The solution is shown at a final time of 0.7. All densities across Figs. 13 to 15 are shown on the same scale to facilitate inter-comparison. Likewise, for all pressures across Figs. 13 to 15 . Similarly for all x-velocities and all x-components of the thermal impulse across Figs. 13 to 15.*

*Fig. 15a, 15b, 15c, 15d show density, temperature, x-velocity and x-component of the thermal impulse when t=0.1 was used along with the plain-vanilla, zone-centered higher order Godunov scheme at third order. The solution is shown at a final time of 0.7. All densities across Figs. 13 to 15 are shown on the same scale to facilitate inter-comparison. Likewise, for all pressures across Figs. 13 to 15 . Similarly for all x-velocities and all x-components of the thermal impulse across Figs. 13 to 15.*

*Fig. 16 shows the evolution of the maximum absolute value of the curl of the thermal impulse vector as a function of time. The results are color-coded for the four different values of t used in the thermally conducting blast problem. Solid curves show the results from the edge-centered curl-preserving (ECCP) algorithm that uses the multidimensional Riemann solver. Dashed curves show the results from the plain-vanilla zone-centered higher order Godunov (HOG) schemes.*

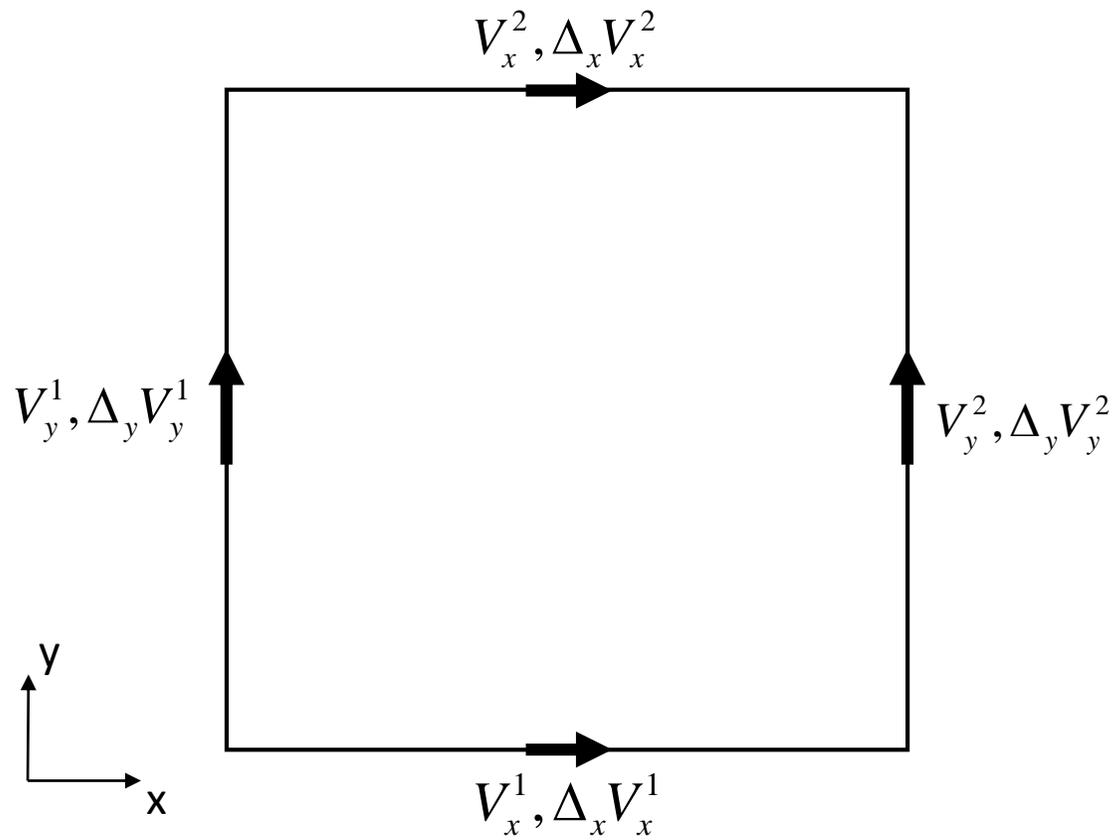

Fig. 1 shows the collocation of vector components along the edges of a two-dimensional control volume. As evaluated over the edges of the square element, the discrete circulation is fully specified. The mean value and its linear variation are shown along each edge, in anticipation of a second order accurate reconstruction scheme. The reconstruction problem for a curl-free reconstruction consists of obtaining a polynomial-based vector field that is globally curl-free within this two-dimensional control volume. The reconstruction problem for a curl-preserving reconstruction consists of obtaining a polynomial-based vector field that matches the specified mean circulation in the zone.

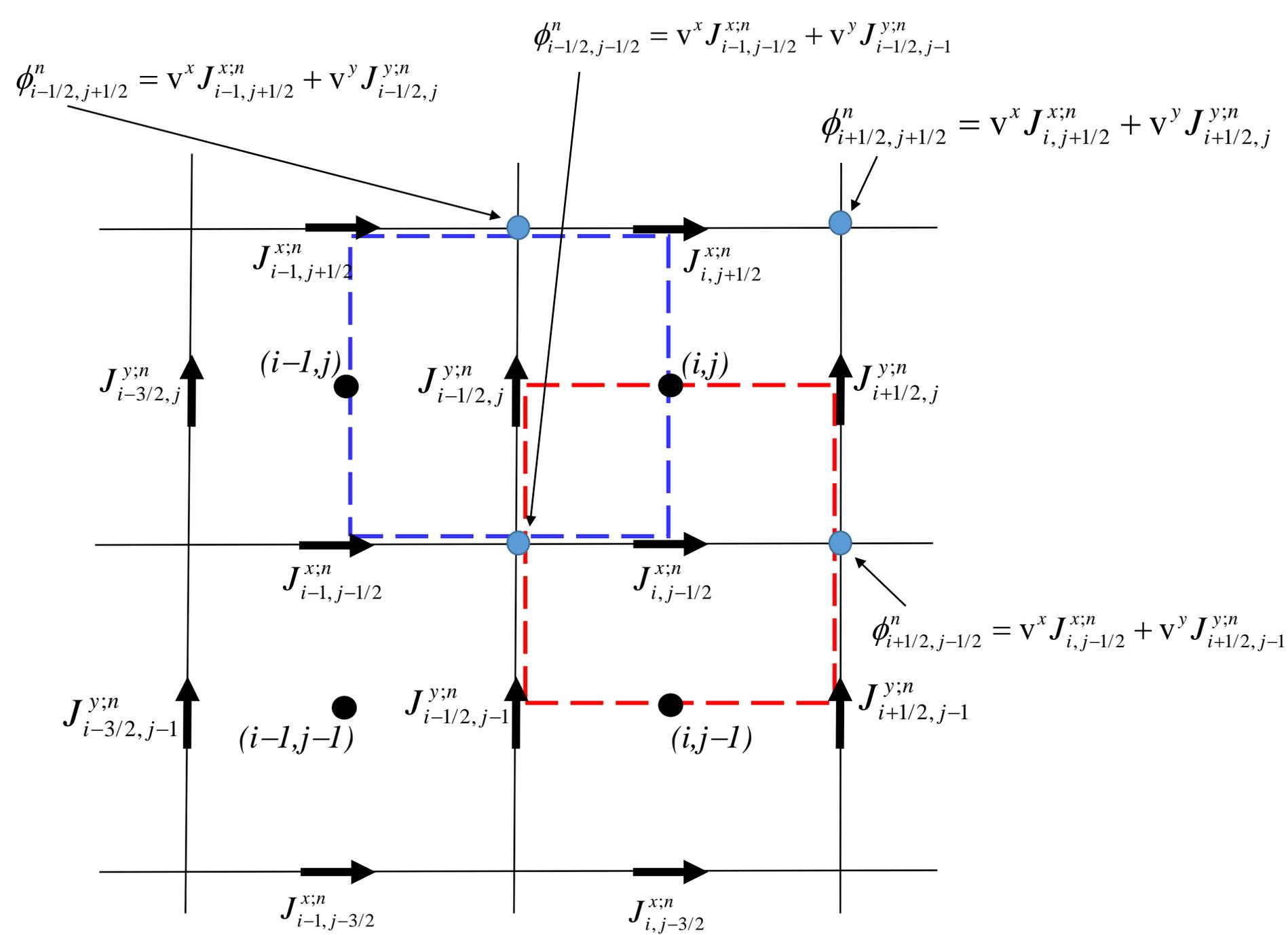

Fig. 2 shows the components of the curl-free vector field around the four zones (i,j), (i-1,j), (i-1,j-1) and (i,j-1). A first order curl-free reconstruction is used. The multidimensionally upwinded potentials at the vertices of the zone (i,j) are also shown. The red dashed rectangle shows the effective control volume that is used for the update of $J^{x;n}_{i,j-1/2}$. The blue dashed rectangle shows the effective control volume that is used for the update of $J^{y;n}_{i-1/2,j}$.

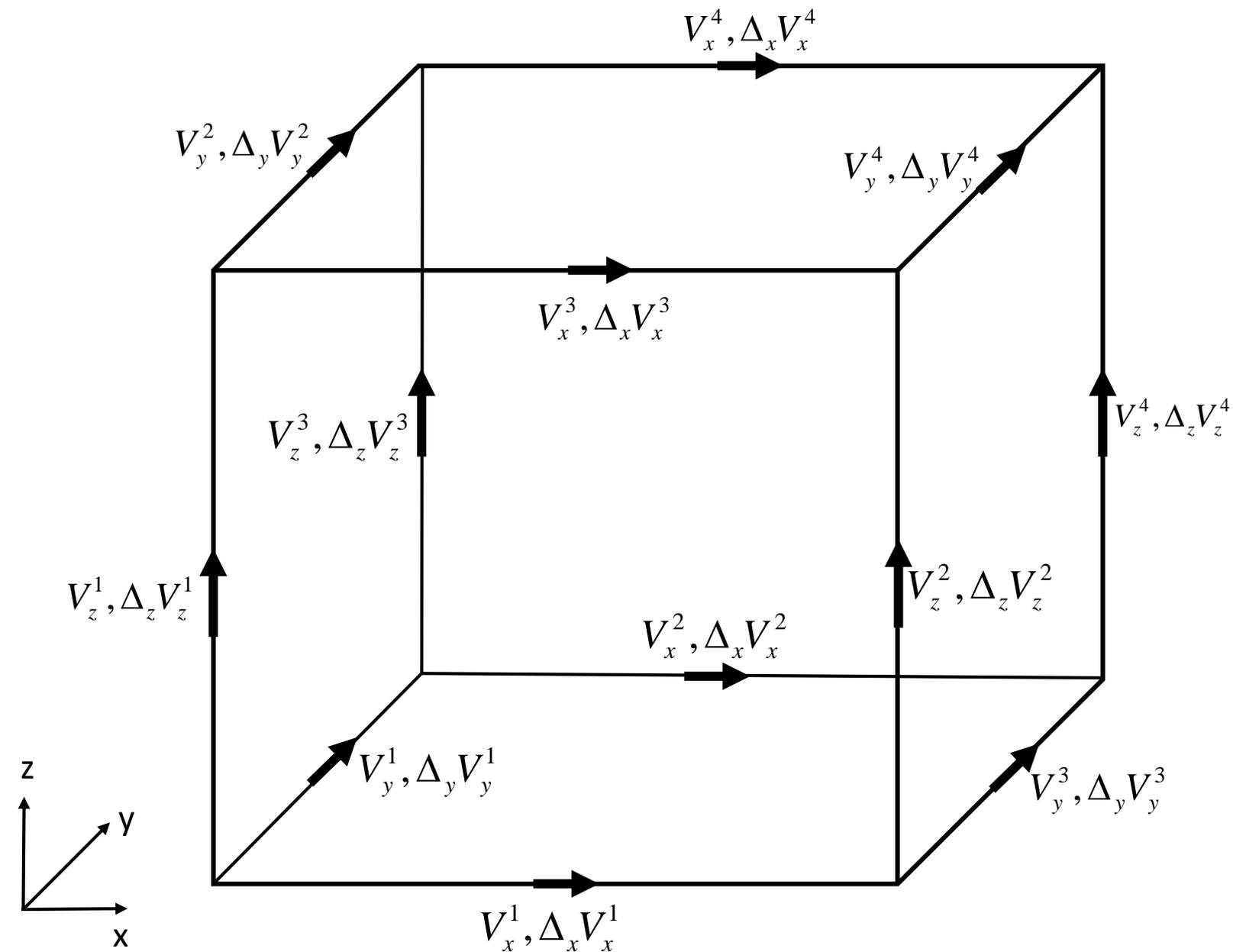

Fig. 3 shows the collocation of vector components along the edges of the control volume. Within each of the two x-faces, the two y-faces and the two z-faces, the discrete circulation (evaluated over those faces) is either exactly zero or specified. The mean value and its linear variation are shown along each edge, in anticipation of a second order accurate reconstruction scheme. The reconstruction problem consists of obtaining a polynomial that is globally curl-free/curl-preserving within this control volume.

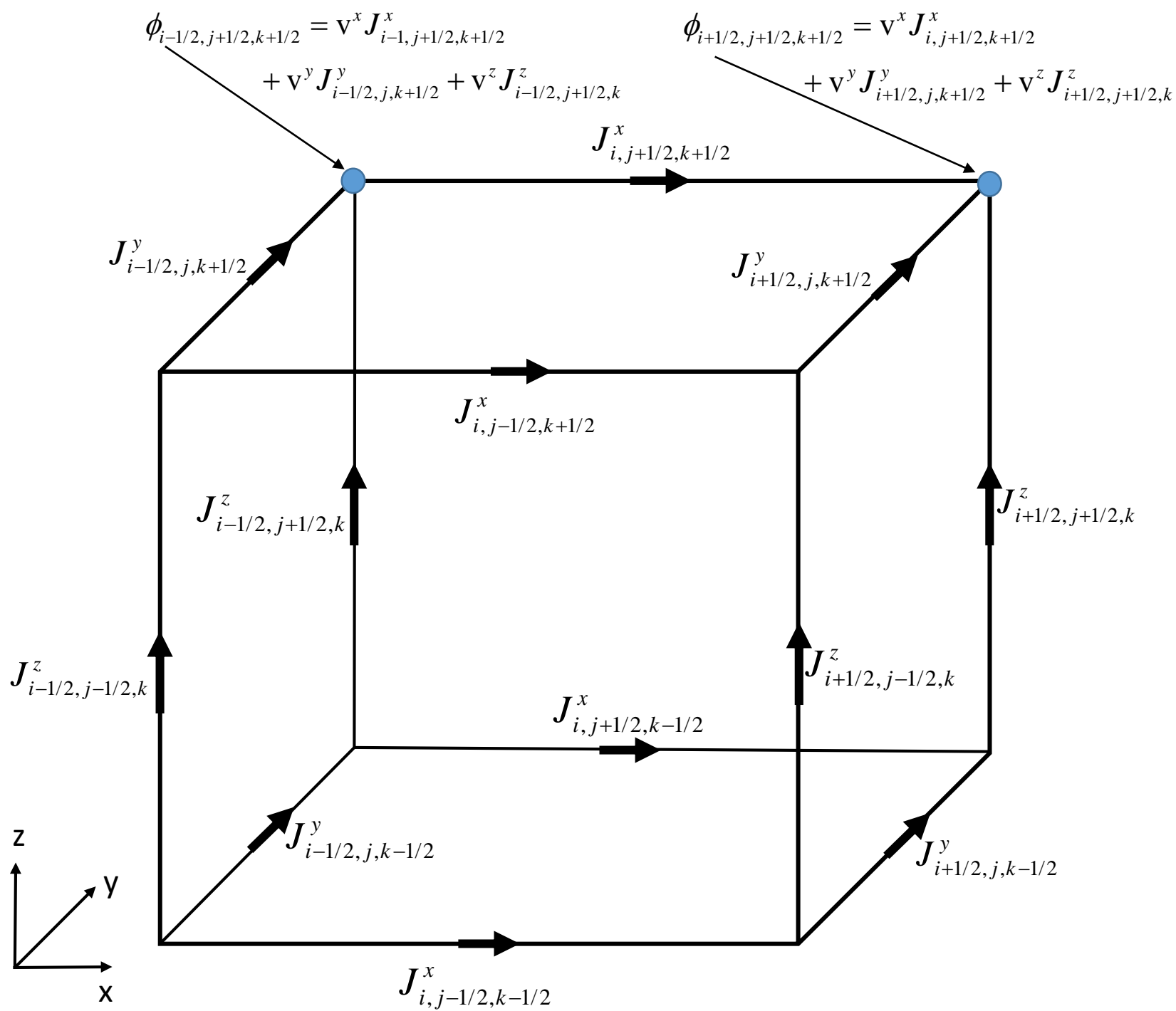

Fig. 4 shows the collocation of curl-free vector components along the edges of a three-dimensional zone. The zone center is indexed by (i,j,k) and the edges are indexed suitably, consistent with the zone center's indexing. We take all the velocity components to be constant and positive. The upwinded potentials at two of the vertices of the mesh are also shown. The potentials at other vertices can be obtained by suitable shifts in the indexing. The purpose of this figure is to make it easy for us to understand how a curl-free reconstruction that is based on edge-centered vector components, in conjunction with a three dimensional Riemann solver, can give us a stable, globally curl-free scheme.

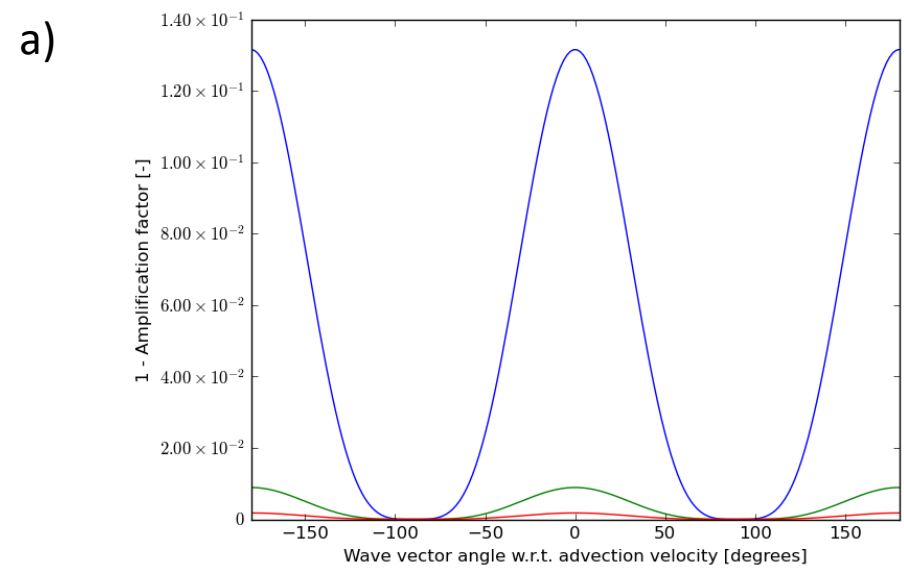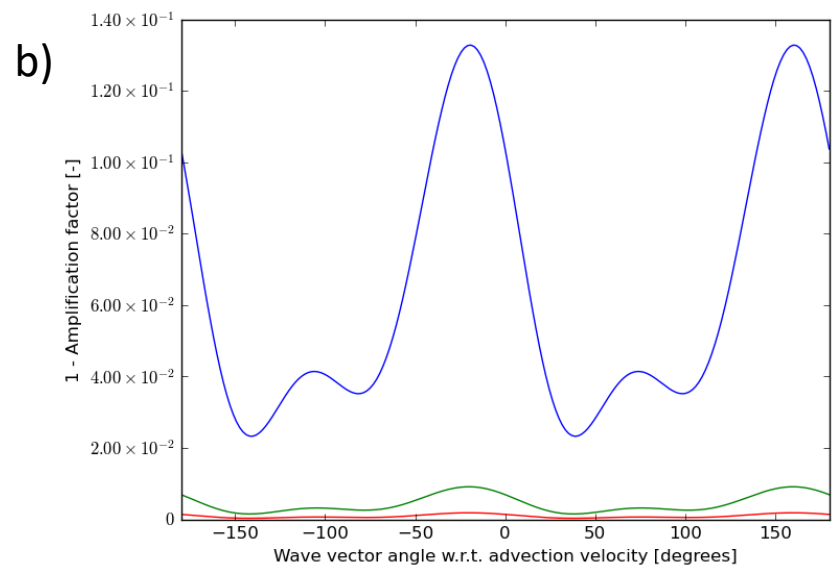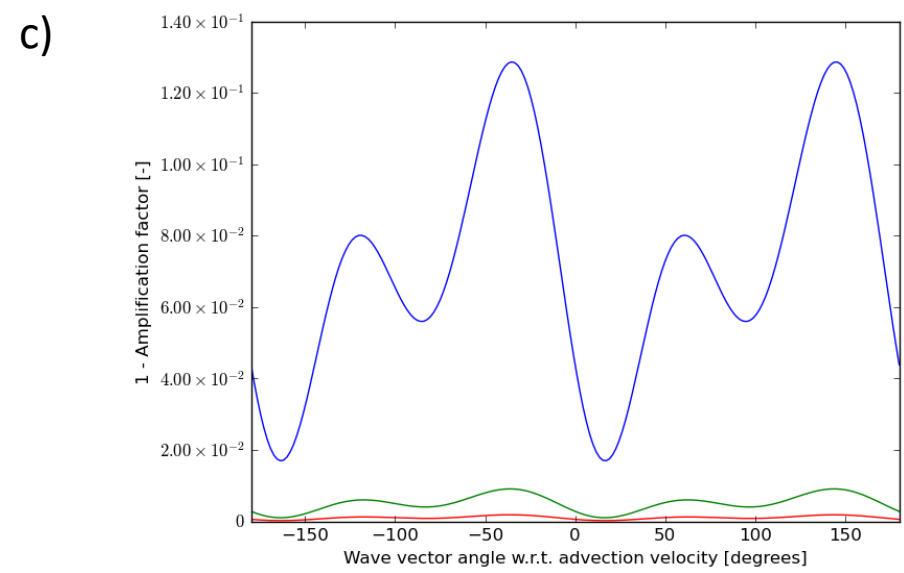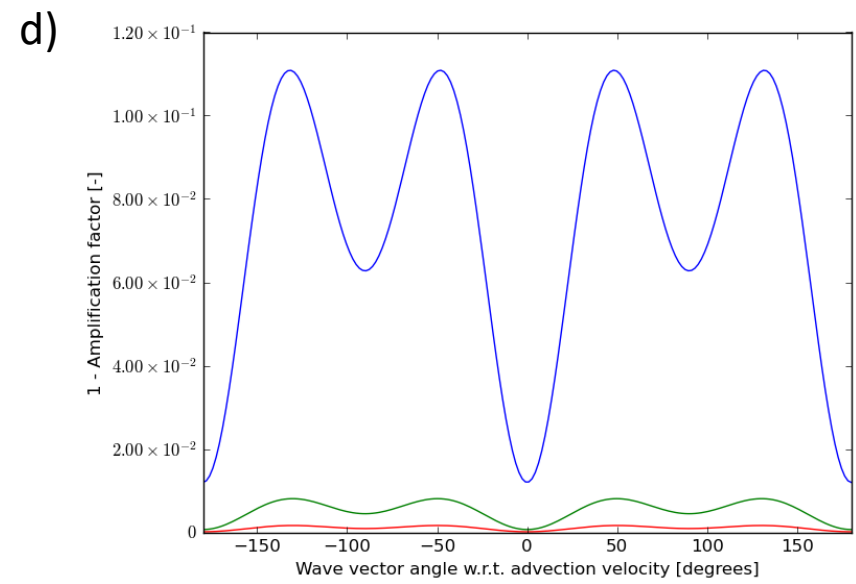

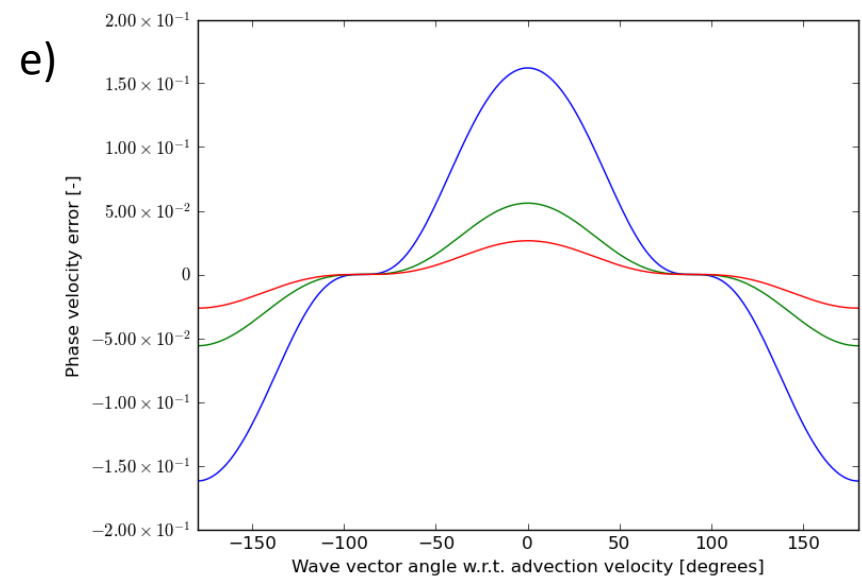
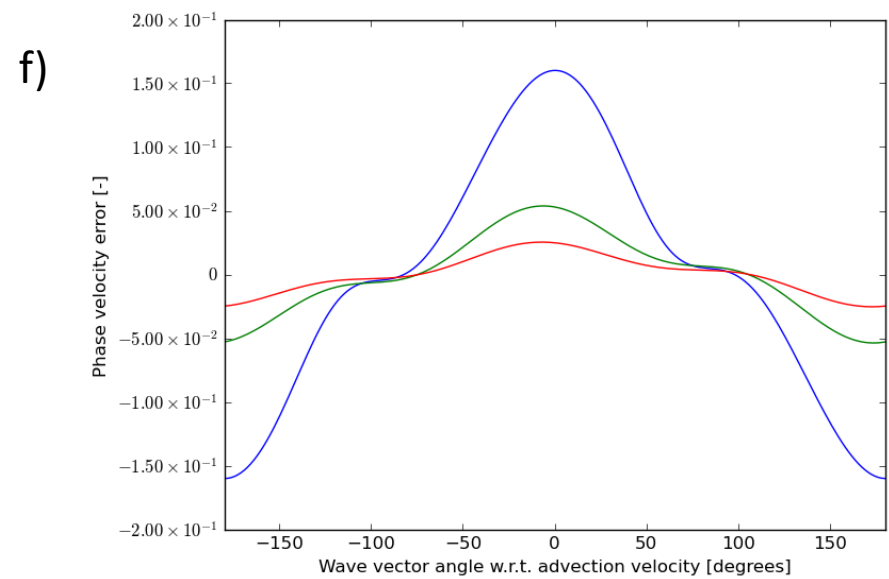
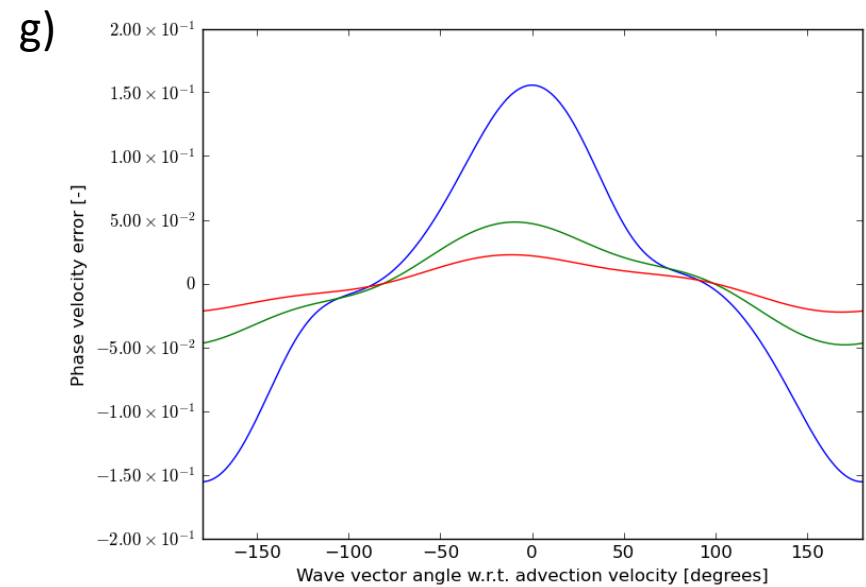
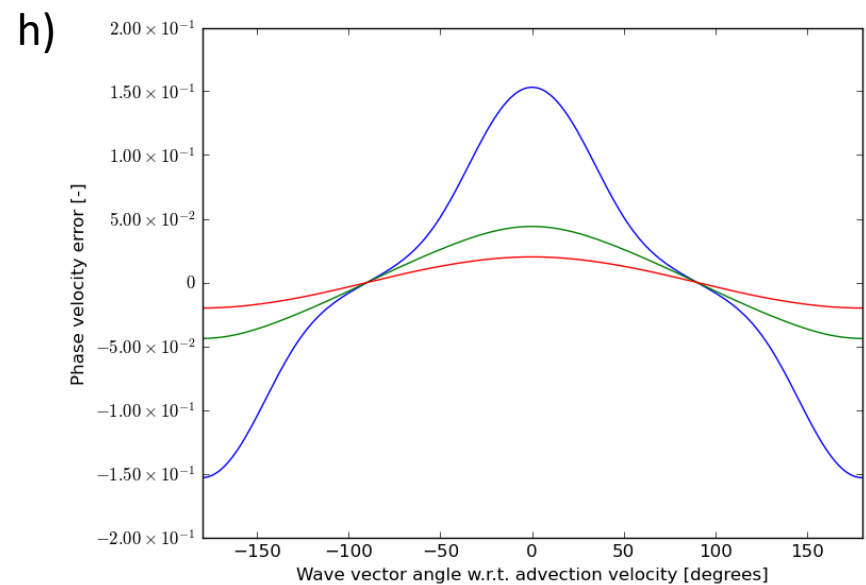

*Fig. 5 shows the wave propagation characteristics for curl-preserving second order WENO-like schemes. Figs. 5a to 5d show one minus the absolute value of the amplification factor when the velocity vector makes angles of 0º , 15º , 30º and 45º relative to the x-direction of the 2D mesh. Figs. 5e to 5h show the phase error, again for the same angles. The 2D wave vector can make any angle relative to the 2D direction of velocity propagation, therefore, the amplitude and phase information are shown w.r.t. the angle made between the velocity direction and the direction of the wave vector. In each plot, the blue curve refers to waves that span 5 cells per wavelength; the green curve refers to waves that span 10 cells per wavelength; the red curve refers to waves that span 15 waves per wavelength.*

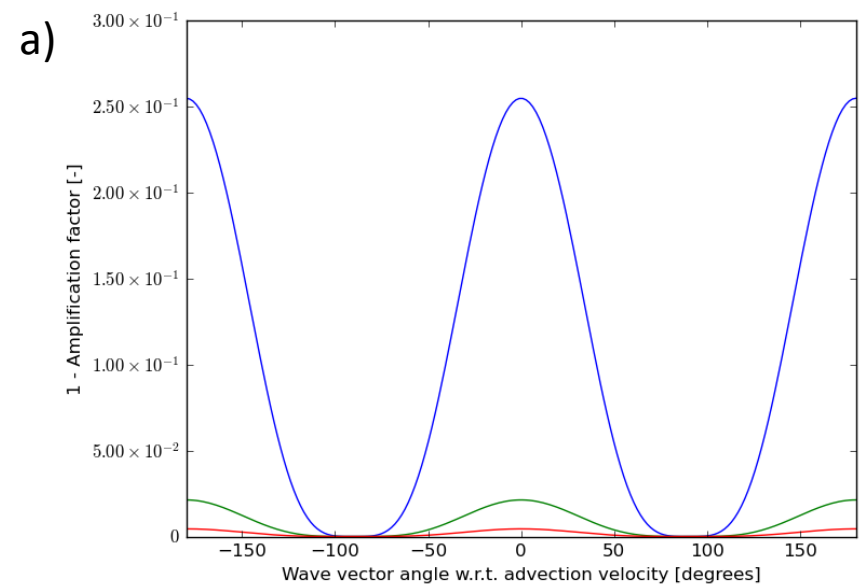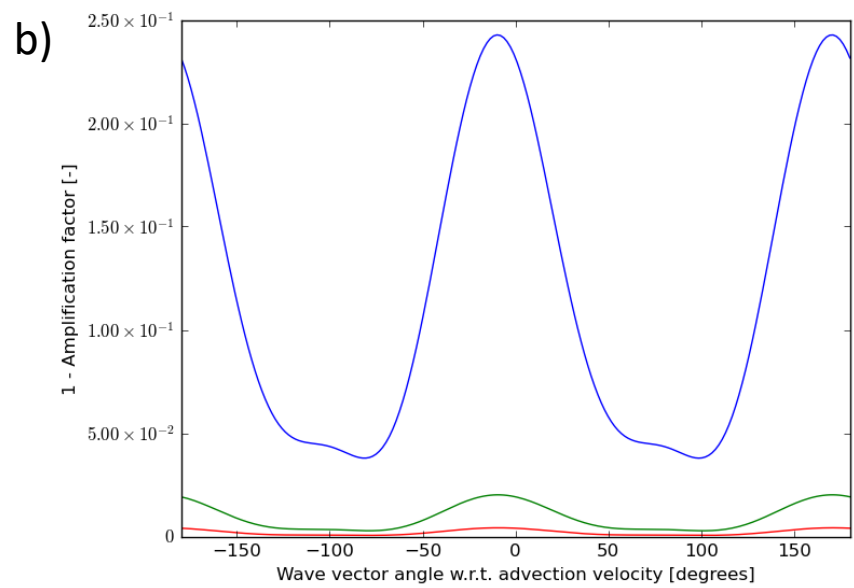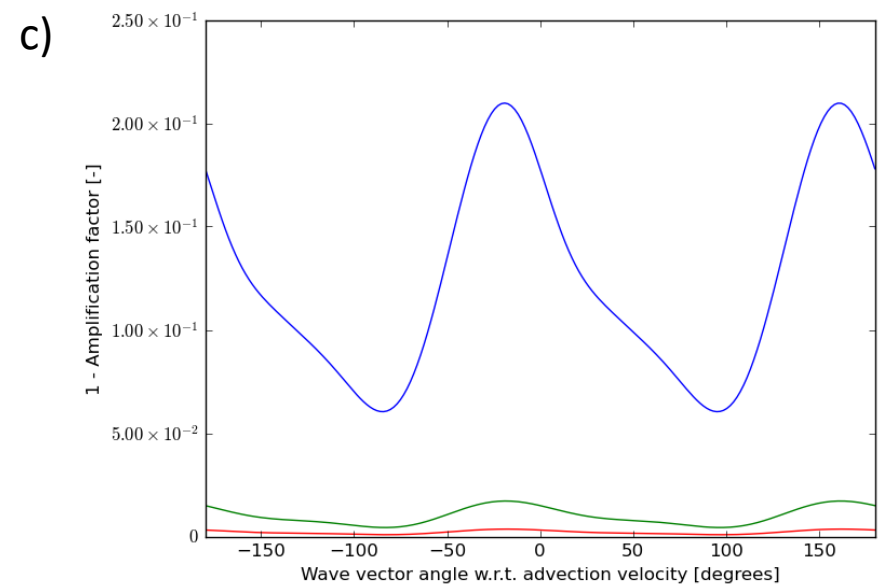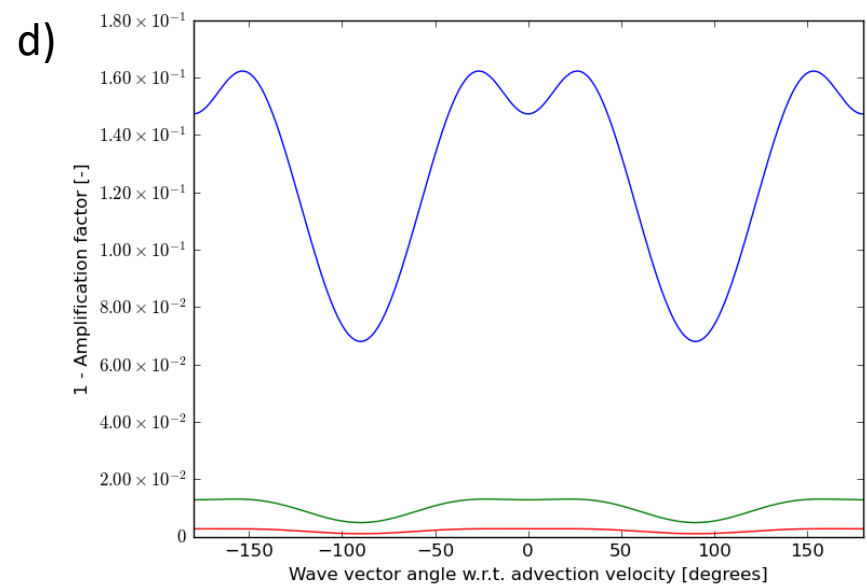

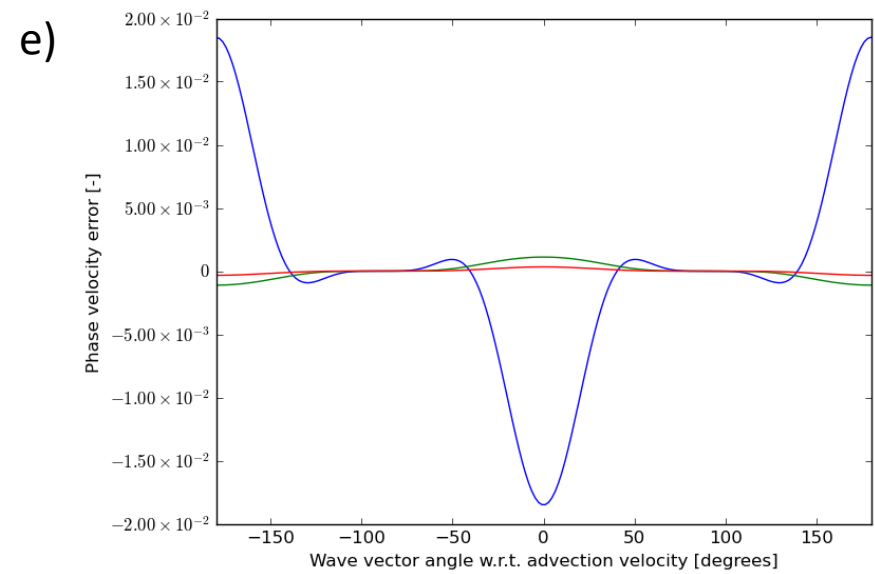
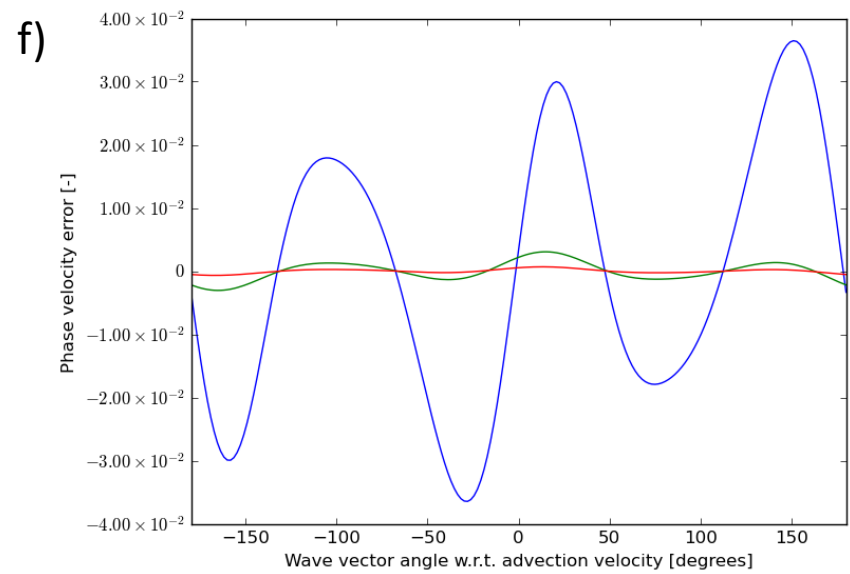
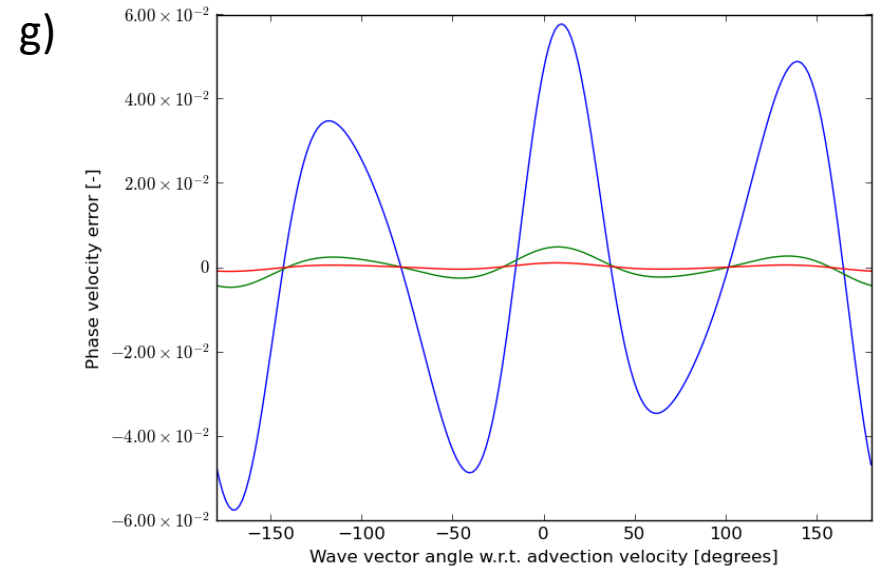
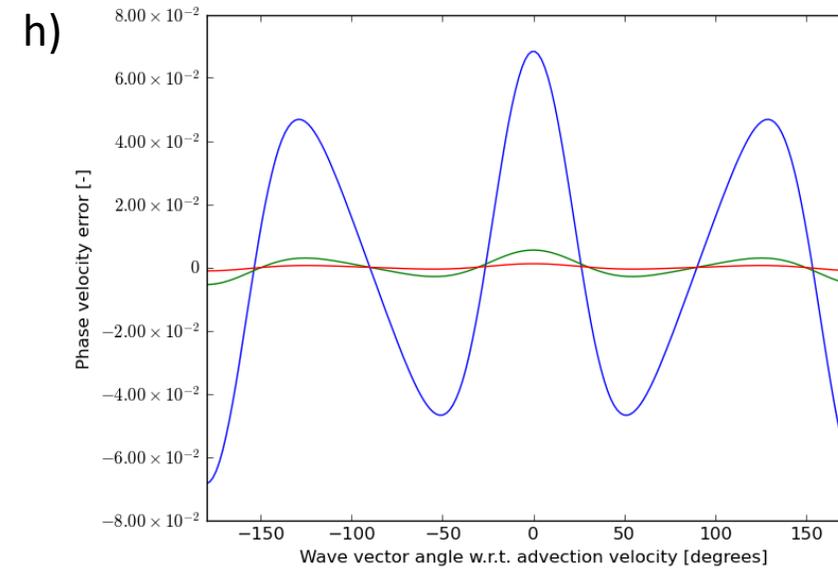

*Fig. 6 shows the wave propagation characteristics for curl-preserving third order WENO-like schemes. Figs. 6a to 6d show one minus the absolute value of the amplification factor when the velocity vector makes angles of 0º, 15º, 30º and 45º relative to the x-direction of the 2D mesh. Figs. 6e to 6h show the phase error, again for the same angles. The 2D wave vector can make any angle relative to the 2D direction of velocity propagation, therefore, the amplitude and phase information are shown w.r.t. the angle made between the velocity direction and the direction of the wave vector. In each plot, the blue curve refers to waves that span 5 cells per wavelength; the green curve refers to waves that span 10 cells per wavelength; the red curve refers to waves that span 15 waves per wavelength.*

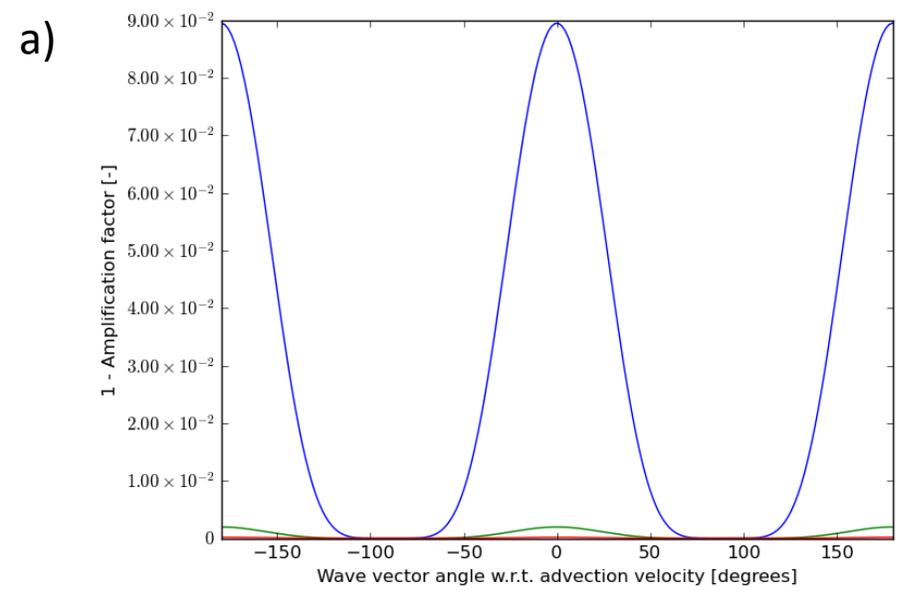
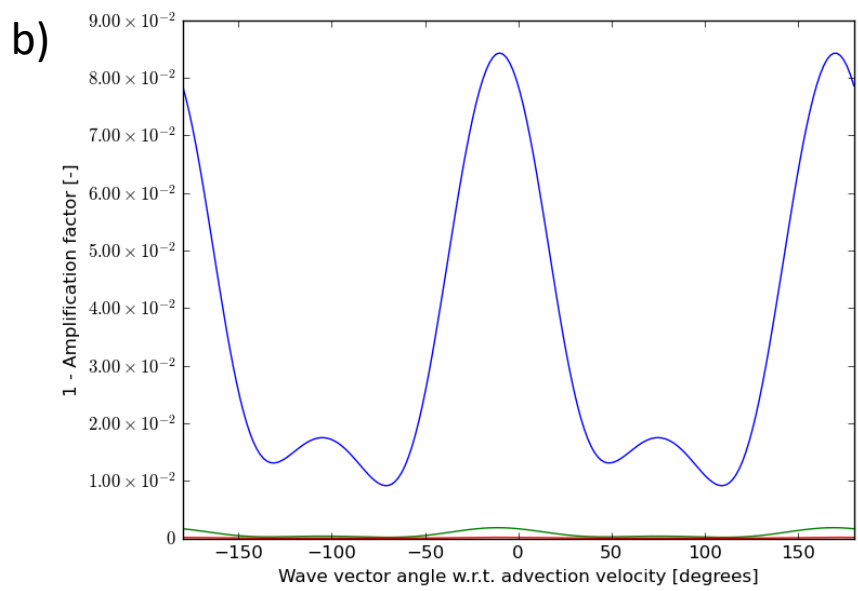
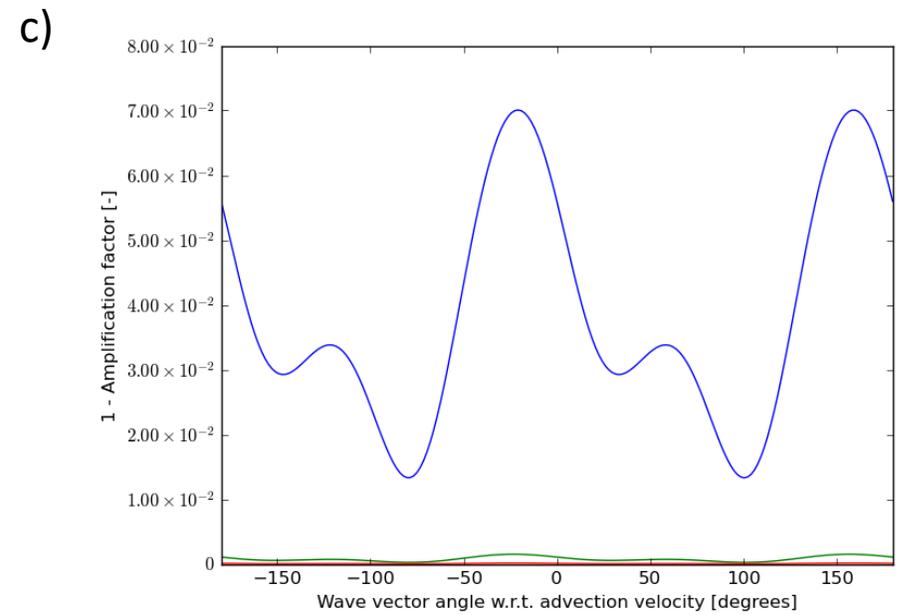
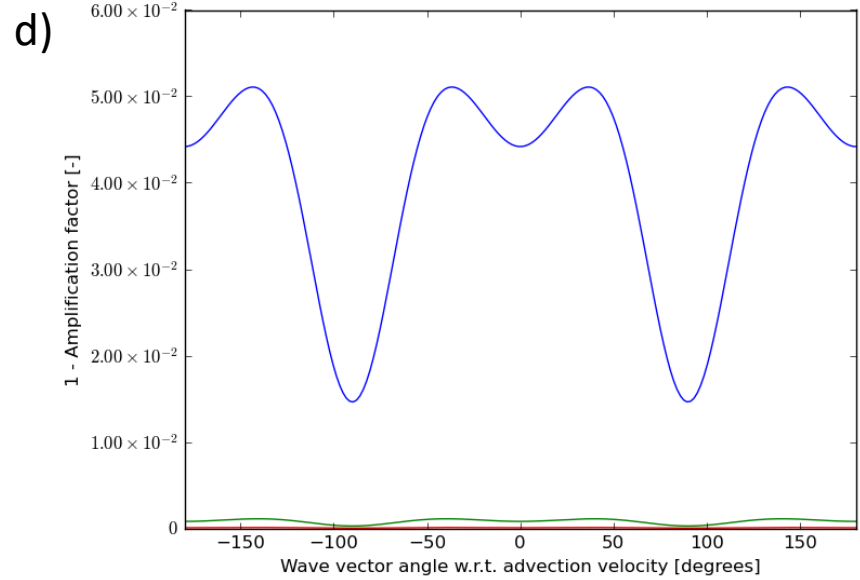

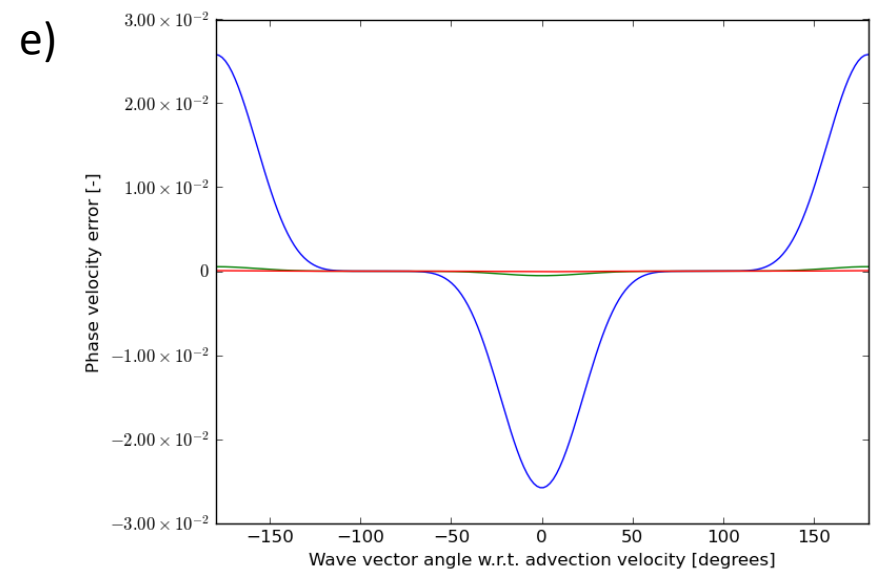
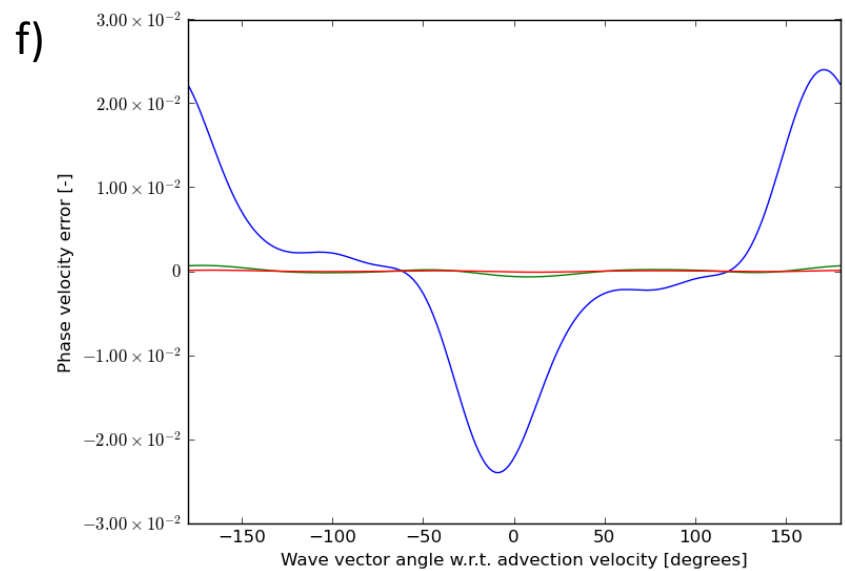
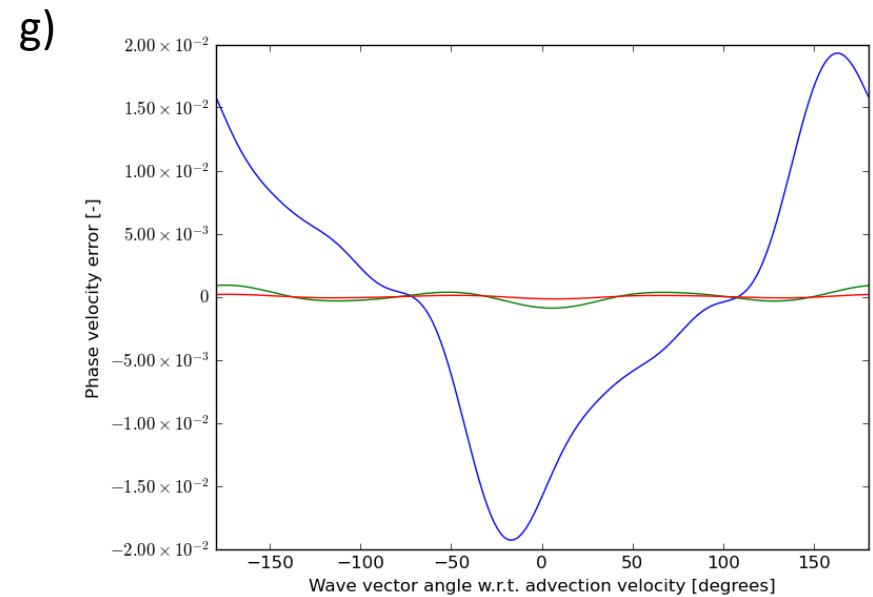
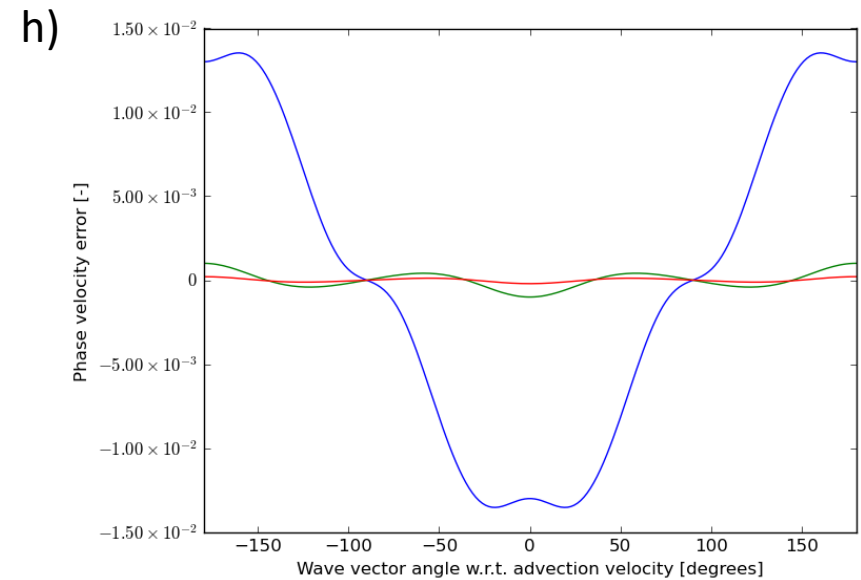

*Fig. 7 shows the wave propagation characteristics for curl-preserving fourth order WENO-like schemes. Figs. 7a to 7d show one minus the absolute value of the amplification factor when the velocity vector makes angles of 0º, 15º, 30º and 45º relative to the x-direction of the 2D mesh. Figs. 7e to 7h show the phase error, again for the same angles. The 2D wave vector can make any angle relative to the 2D direction of velocity propagation, therefore, the amplitude and phase information are shown w.r.t. the angle made between the velocity direction and the direction of the wave vector. In each plot, the blue curve refers to waves that span 5 cells per wavelength; the green curve refers to waves that span 10 cells per wavelength; the red curve refers to waves that span 15 waves per wavelength.*

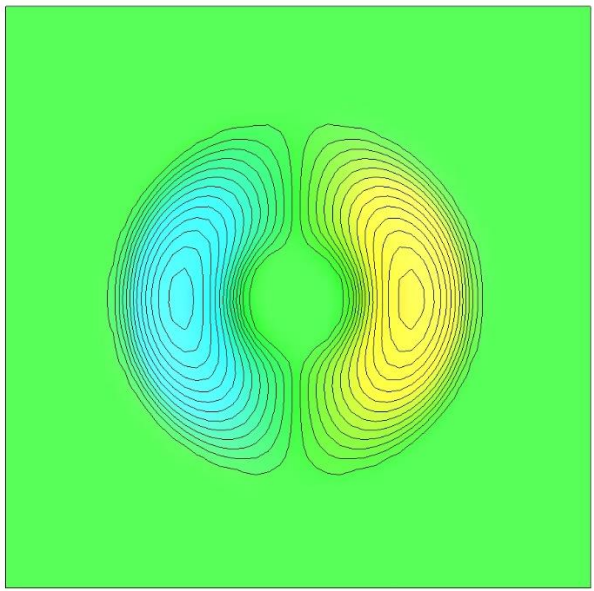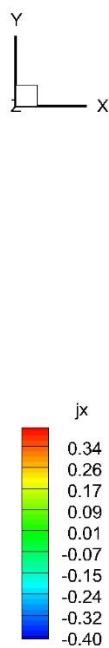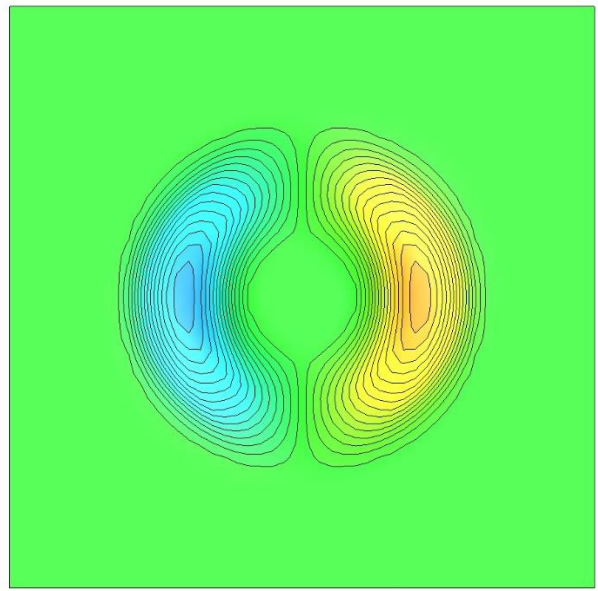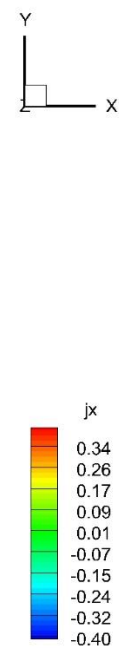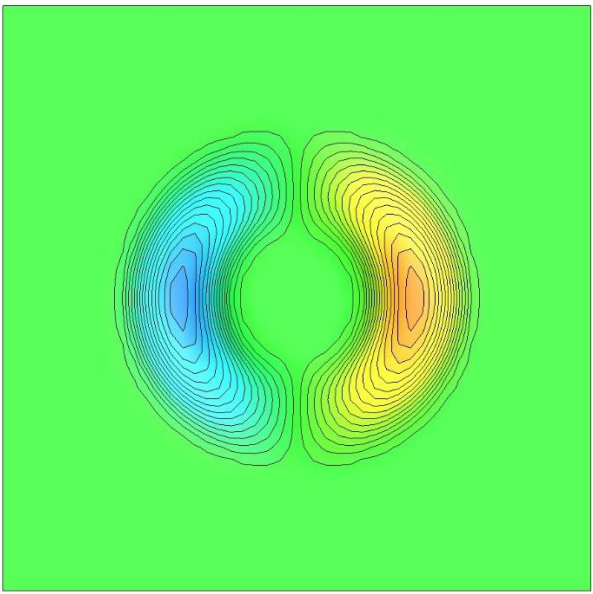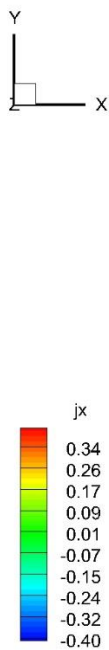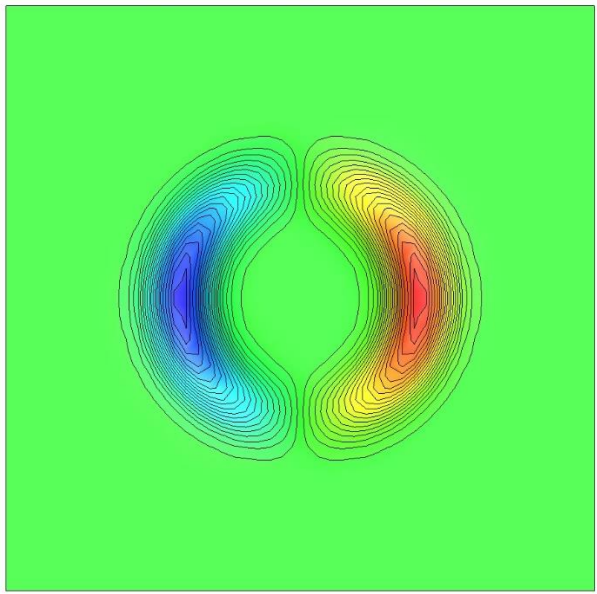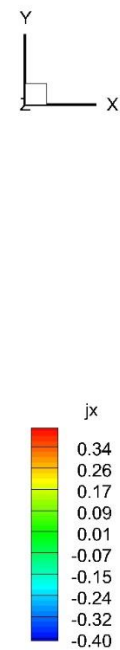

*Fig 8 shows four different simulations of our model test problem at 50X50 zone resolution. All figures show the y-component of the curl-free vector field after a very long simulation time; i.e. stopping time of 100. The color scheme is the same across all figures and the intervals between contours is also the same. Fig. 8a shows the result of the second-order curl free scheme from Boscheri et al. (2020) which frequently averages the solution to different locations of the mesh in the course of a timestep. Figs. 8b, 8c and 8d show the results from our new scheme that uses curl-free reconstruction at second, third and fourth order, along with the multidimensional Riemann solver. The color saturation is seen to increase with increasing order and the solution is also less diffusive as one progresses from Fig. 8a to 8d. Very different second order schemes were used for Figs. 8a and 8b; we see that Fig. 8b shows appreciably lower diffusion.*

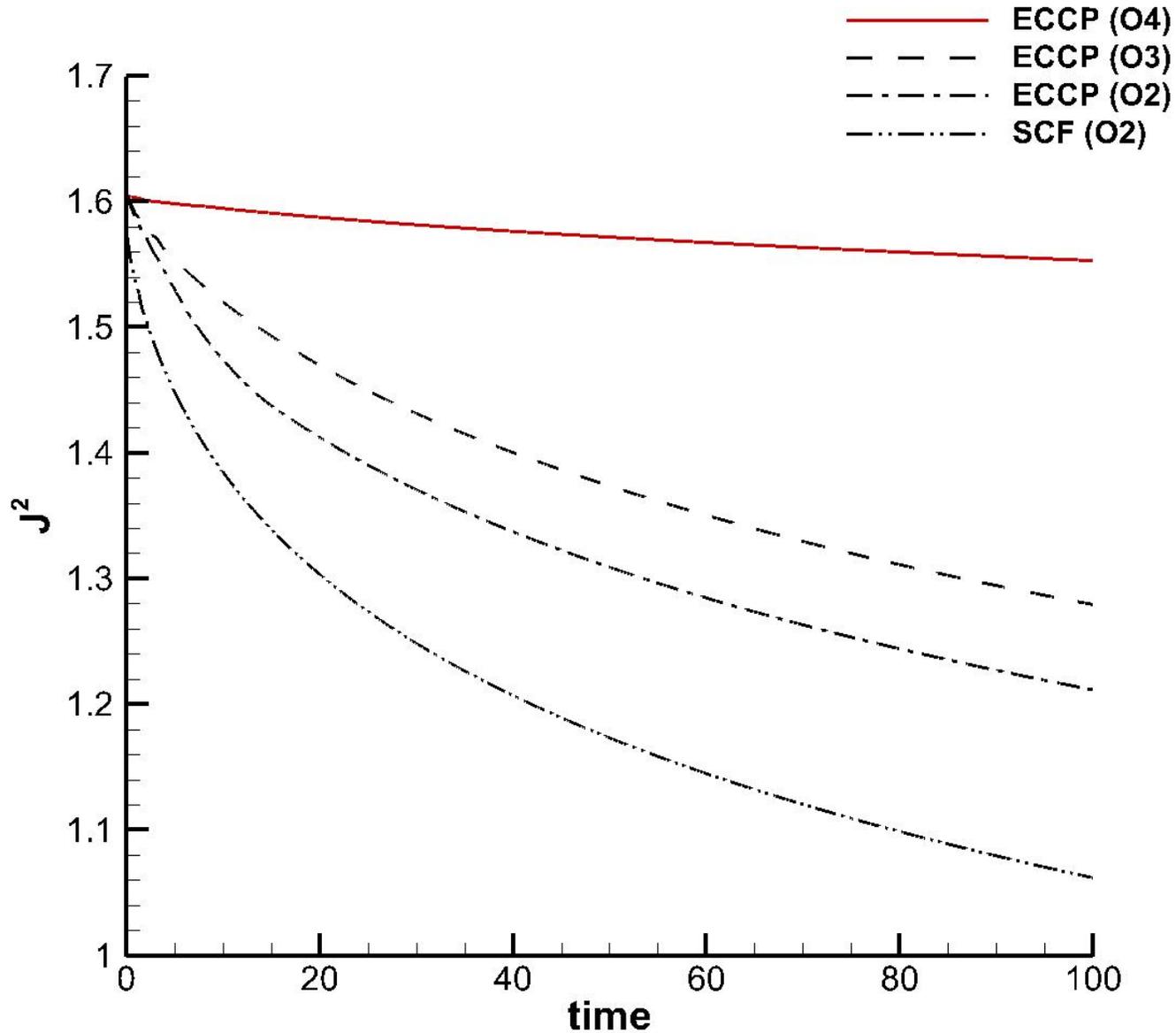

*Fig. 9 shows the evolution of the quadratic energy in the curl-free vector field as a function of time for the four schemes shown in Fig. 8. The model problem was run with 50X50 zone resolution to a time of 100. In principle, the quadratic energy of the constrained vector field integrated over the mesh should stay unchanged. The plot shows how much the quadratic energy decays with time for various schemes. Because of the frequent averaging over the course of a timestep, the second order SCF(O2) scheme from Dumbser et al. [35] shows appreciable decay of energy. The second order ECCP(O2) scheme designed here improves on the SCF(O2) scheme. The third order ECCP(O3) scheme shows substantial improvement over the second order schemes. The fourth order ECCP(O4) scheme provides the best preservation of quadratic energy, thus highlighting the value of well-designed higher order curl-preserving schemes for the evolution of curl-constrained vector fields.*

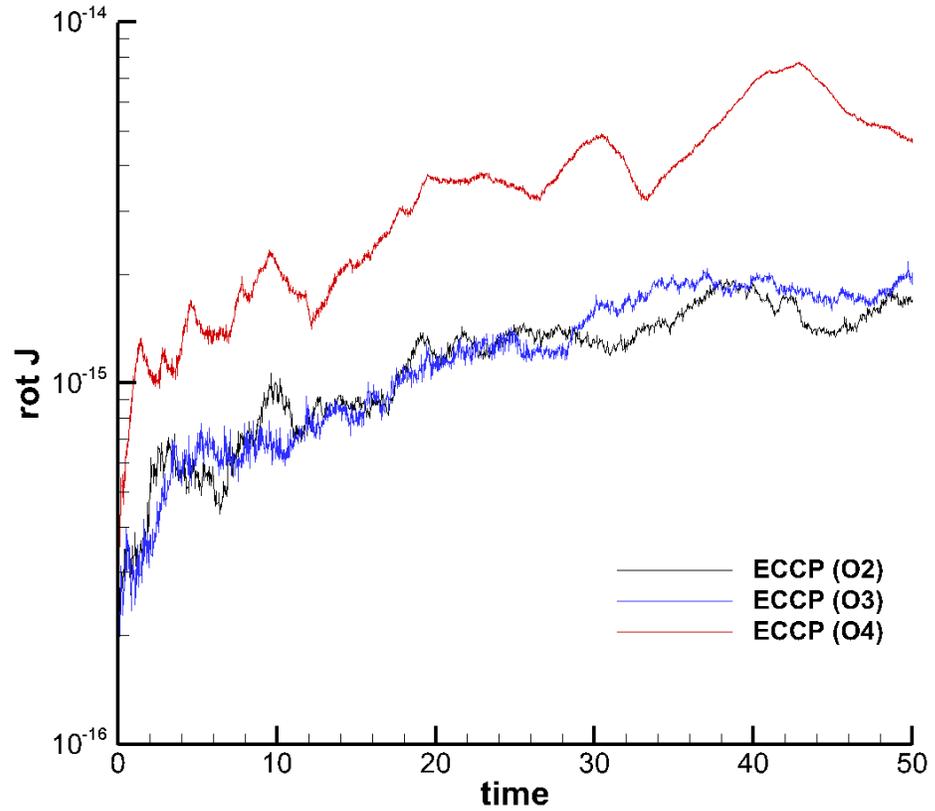

Fig. 10 shows the time series of the maximum pointwise error of the curl of ***J*** for the stationary curl free test problem until t=50. The computational domain is discretized with 32x32 uniform Cartesian cells. The results for second (black solid line), third (blue solid line) and fourth (red solid line) order curl-free WENO schemes are shown.

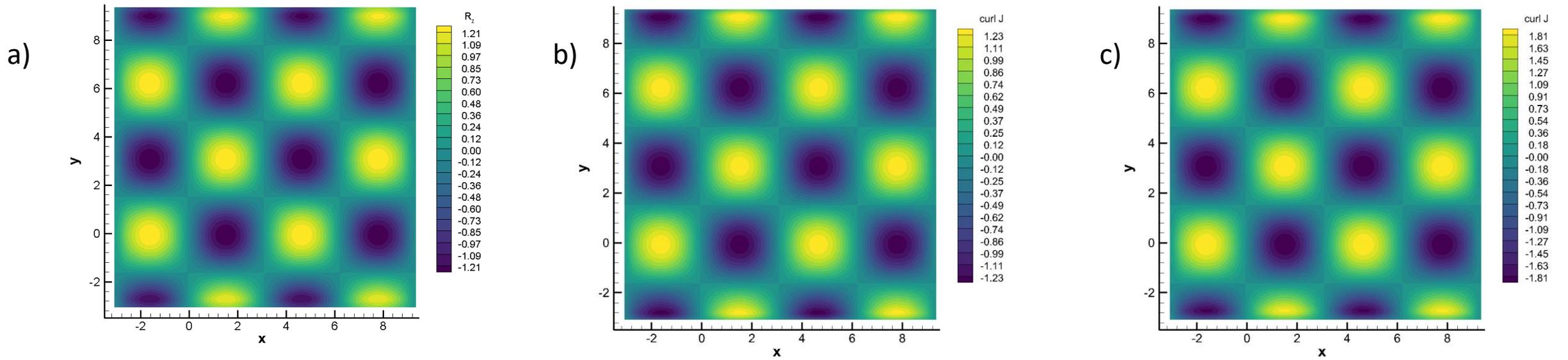

Figs 11a, b, c show the colorized plot of Rz as a function of position in the test problem at a time of 1. Because the problem lends itself to analytical treatment, Fig. 11a was obtained analytically. Fig. 11b was obtained numerically. Using the color bar, please note the close concordance in values between Figs. 11a and 11b. Fig. 11c was obtained from direct differentiation of a zone-centered higher order Godunov scheme that is not curl-preserving. Please use the color bars to compare the numerical values in Fig. 11c to those in Fig. 11a to see the significant differences. The color bar for each panel is different and scaled to the min and max of the data shown.

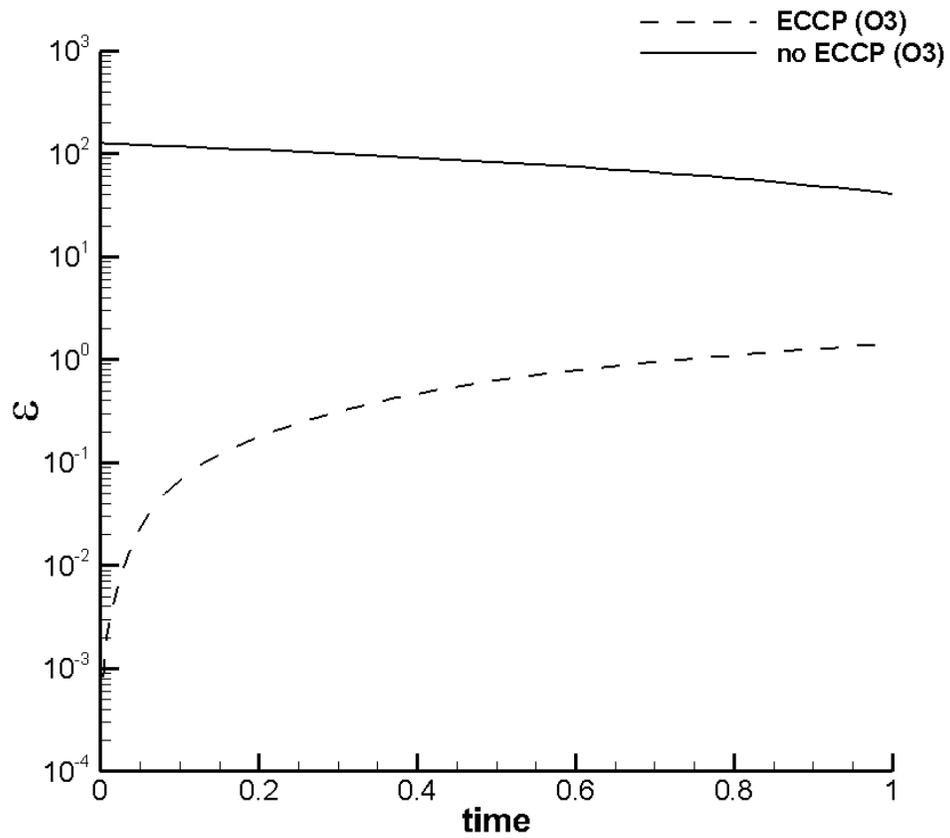

*Fig. 12 shows the maximum error in Rz as a function of time for the edge-centered curl-preserving scheme versus a higher order Godunov scheme that does not preserve the curl. The latter shows errors that are 2 to 4 magnitudes higher.*

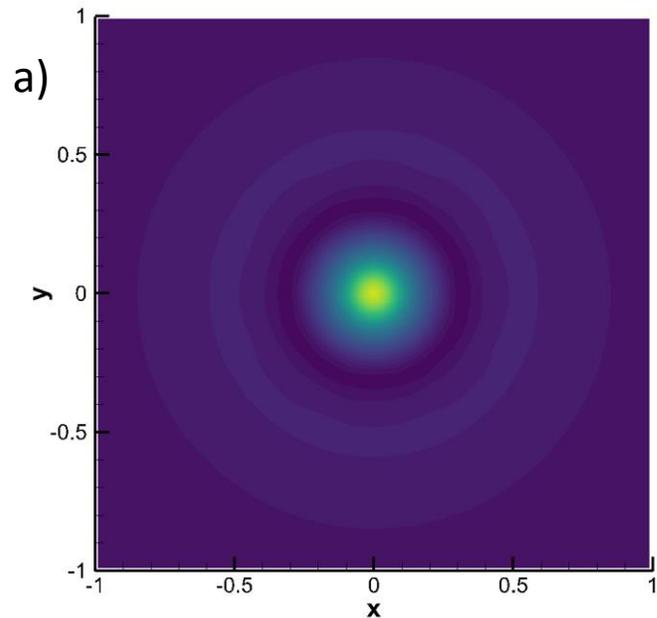
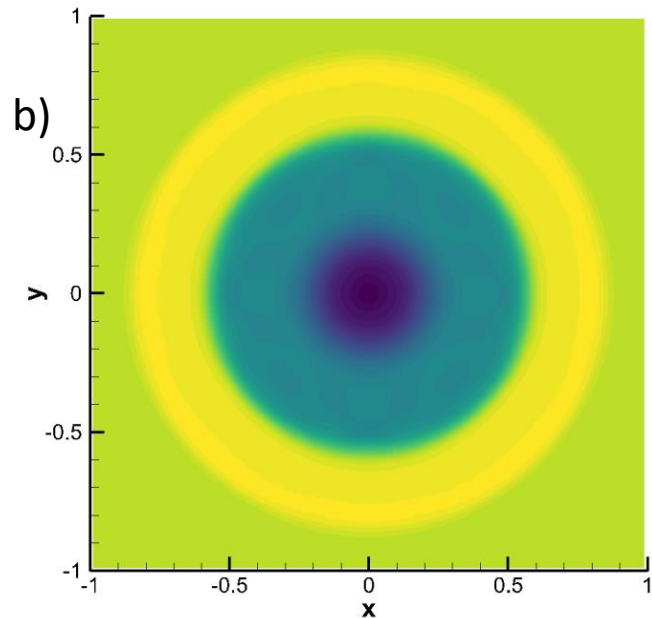
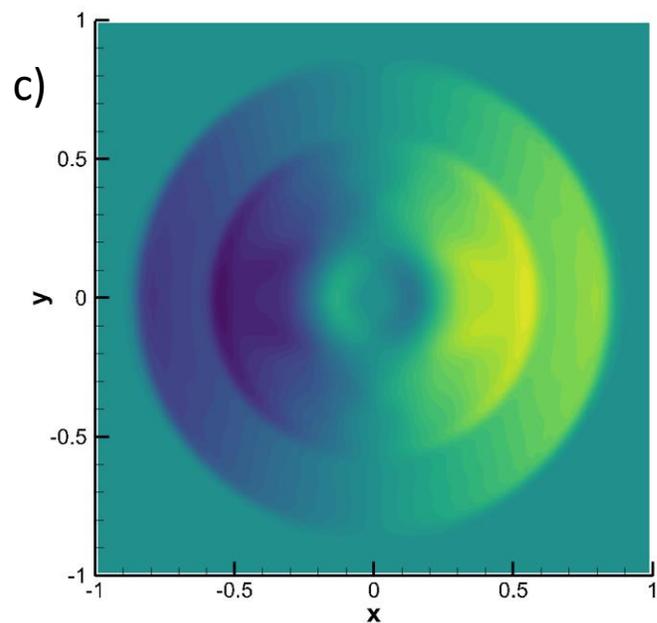
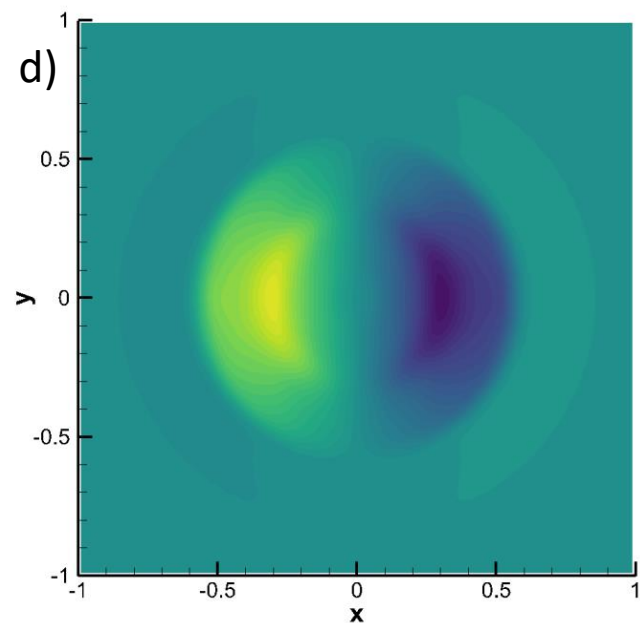

*Fig. 13a, 13b, 13c, 13d show density, temperature, x-velocity and x-component of the thermal impulse when t=1 was used along with the ECCP scheme at third order. The solution is shown at a final time of 0.7. All densities across Figs. 13 to 15 are shown on the same scale to facilitate inter-comparison. Likewise, for all pressures across Figs. 13 to 15. Similarly for all x-velocities and all x-components of the thermal impulse across Figs. 13 to 15.*

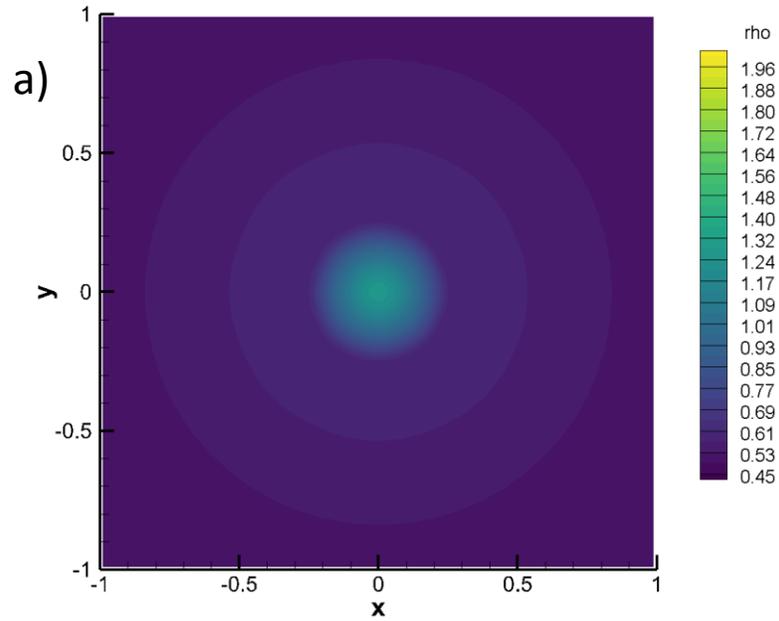
a)

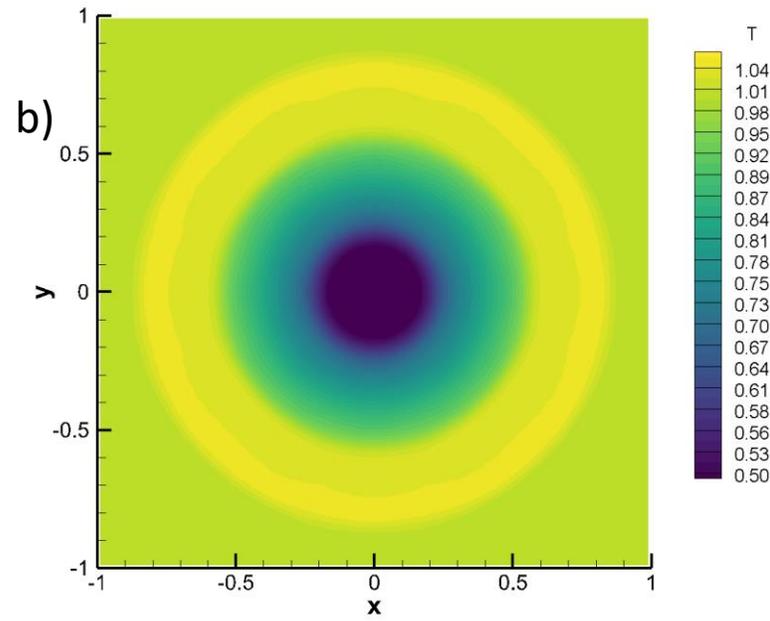
b)

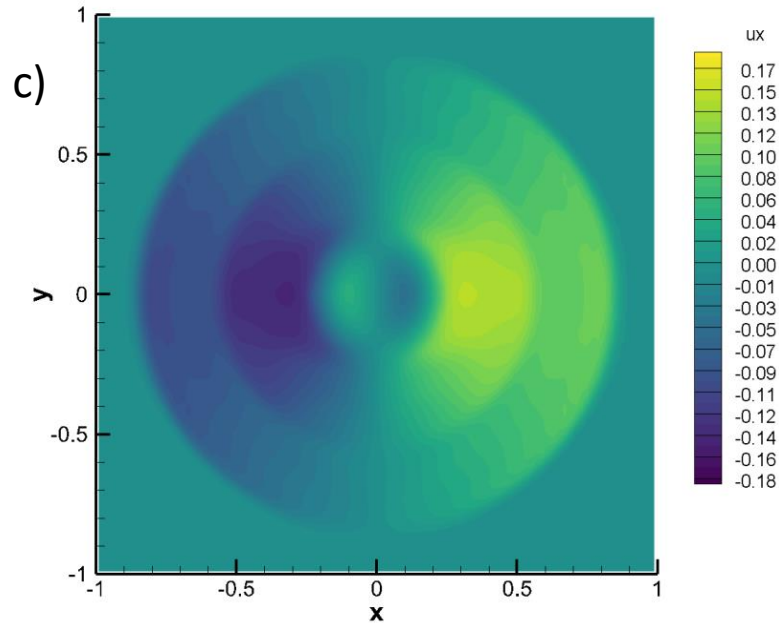
c)

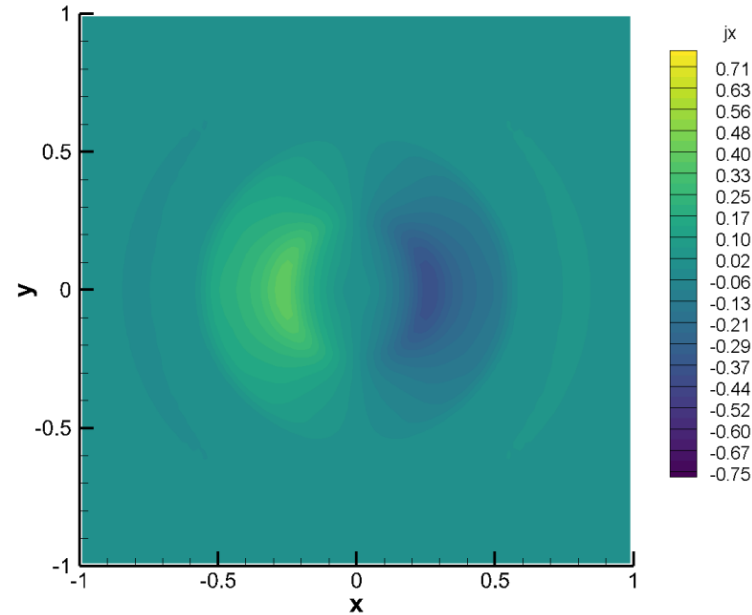

*Fig. 14a, 14b, 14c, 14d show density, temperature, x-velocity and x-component of the thermal impulse when t=0.1 was used along with the ECCP scheme at third order. The solution is shown at a final time of 0.7. All densities across Figs. 13 to 15 are shown on the same scale to facilitate inter-comparison. Likewise, for all pressures across Figs. 13 to 15. Similarly for all x-velocities and all x-components of the thermal impulse across Figs. 13 to 15.*

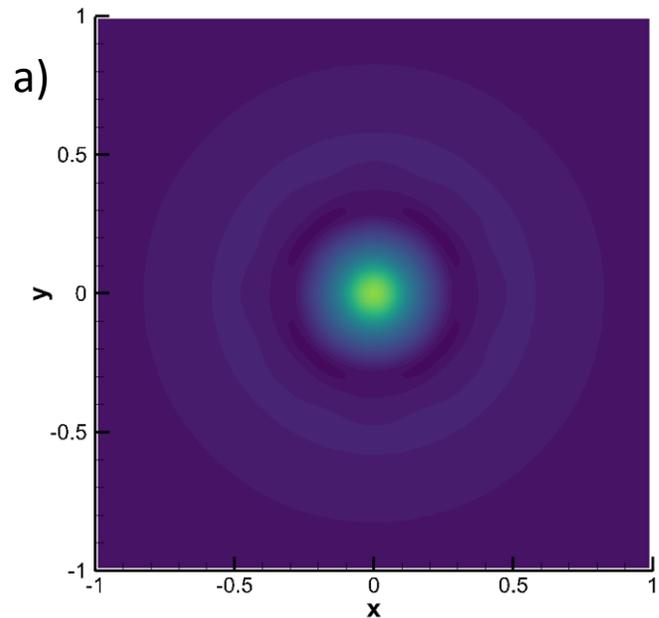
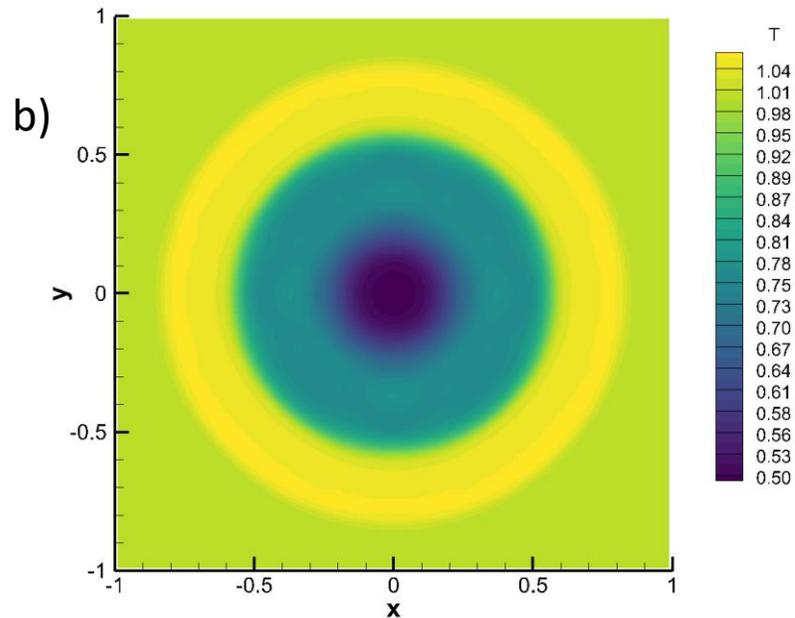
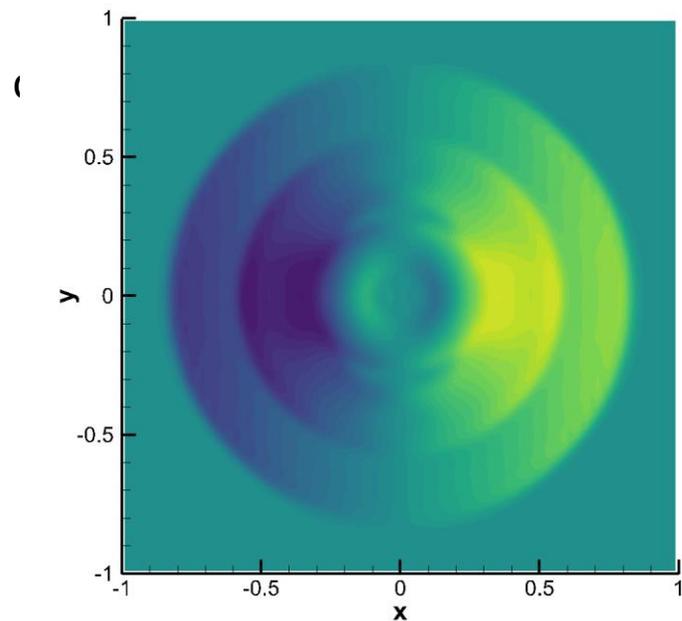
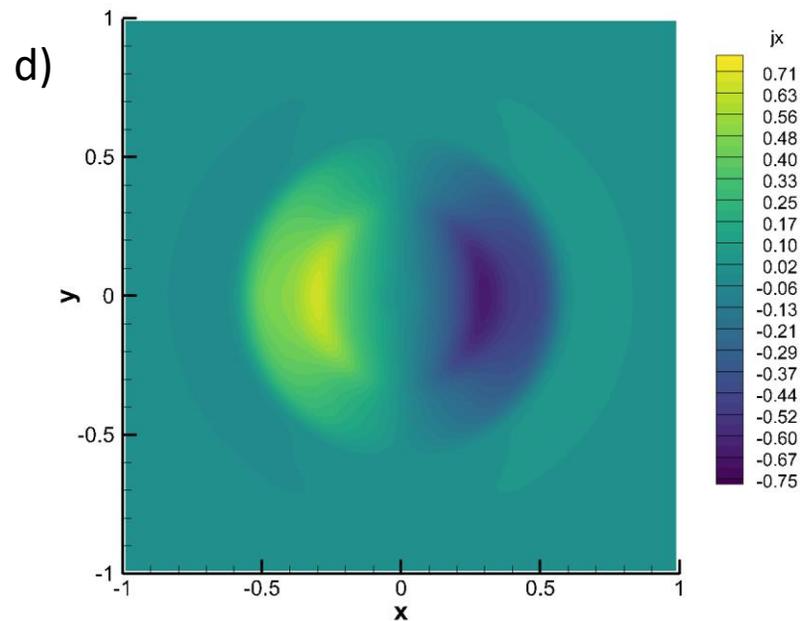

*Fig. 15a, 15b, 15c, 15d show density, temperature, x-velocity and x-component of the thermal impulse when t=0.1 was used along with the plain-vanilla, zone-centered higher order Godunov scheme at third order. The solution is shown at a final time of 0.7. All densities across Figs. 13 to 15 are shown on the same scale to facilitate inter-comparison. Likewise, for all pressures across Figs. 13 to 15 . Similarly for all x-velocities and all x-components of the thermal impulse across Figs. 13 to 15.*

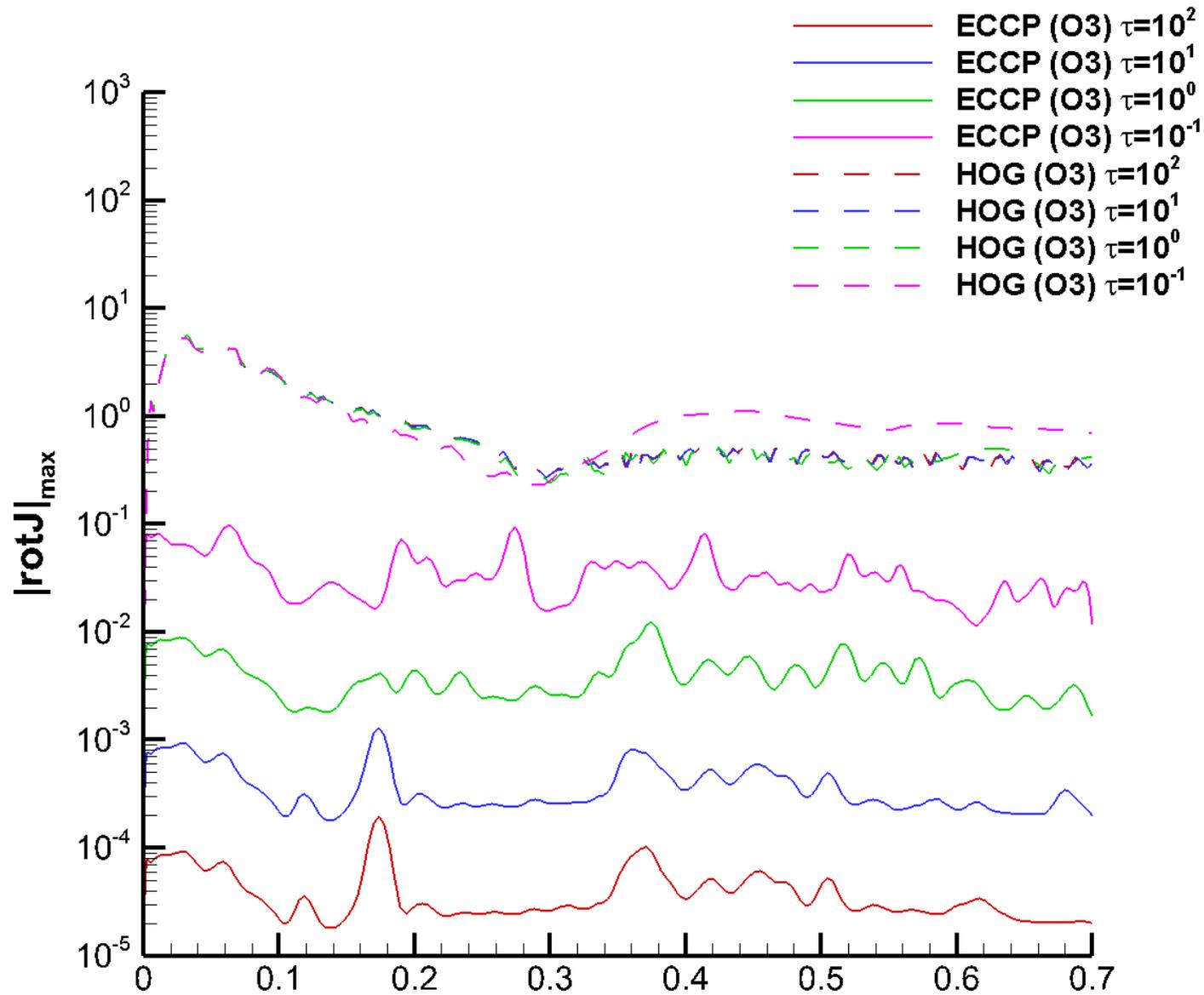

Fig. 16 shows the evolution of the maximum absolute value of the curl of the thermal impulse vector as a function of time. The results are color-coded for the four different values of $\tau$ used in the thermally conducting blast problem. Solid curves show the results from the edge-centered curl-preserving (ECCP) algorithm that uses the multidimensional Riemann solver. Dashed curves show the results from the plain-vanilla zone-centered higher order Godunov (HOG) schemes.